\documentclass[11pt]{amsart}
\usepackage{amssymb}
\usepackage{hyperref}

\newtheorem{theorem}{Theorem}[section]
\newtheorem{splittingtheorem}[theorem]{Splitting Theorem}
\newtheorem{slicetheorem}[theorem]{Slice Theorem}
\newtheorem{maintheorem}{Theorem}
\newtheorem{maincorollary}{Corollary}
\newtheorem{lemma}[theorem]{Lemma}
\theoremstyle{definition}
\newtheorem{definition}[theorem]{Definition}
\theoremstyle{remark}
\newtheorem{remark}[theorem]{\bf Remark}
\numberwithin{equation}{section} \theoremstyle{corollary}
\newtheorem{corollary}[theorem]{Corollary}
\theoremstyle{proposition}
\newtheorem{proposition}[theorem]{Proposition}

\newfont{\EUL}{eufm10 scaled 1000}

\def\Sp{{\rm Sp}}

\newcommand\Ad{{\rm Ad}}
\newcommand\Aut{{\rm Aut}}
\newcommand\Ca{\mathbb{O}}

\newcommand\C{\mathbb{C}}

\newcommand\Exp{{\rm exp}}

\newcommand\Gl{{\rm G \ell}}
\newcommand\G{{\mathbb G}}

\newcommand\Isom{{\rm Isom}}
\newcommand\K{\mathbb{K}}
\newcommand\LA{{\rm A}}
\newcommand\LB{{\rm B}}
\newcommand\LC{{\rm C}}

\newcommand\LE{{\rm E}}
\newcommand\LF{{\rm F}}
\newcommand\LG{{\rm G}}
\newcommand\La{\mathfrak a}

\newcommand\Le{\mathfrak e}
\newcommand\Lf{\mathfrak f}

\newcommand\N{{\rm N}}
\newcommand\Nat{{\mathbb N}}

\newcommand\R{\mathbb{R}}
\newcommand\Rho{\rm P}
\newcommand\SO{{\rm SO}}
\newcommand\SUxU[2]{ {\rm S(U(}#1) \times {\rm U(}#2{\rm ))} }
\newcommand\SU{{\rm SU}}

\newcommand\Spin{{\rm Spin}}
\newcommand\Sym{\rm Sym}
\newcommand\T{{\rm T}}
\newcommand\U{{\rm U}}

\newcommand\bl{\vspace{1em}\hfill\\}

\newcommand\e{{\rm e}}
\newcommand\eK{{\e K}}
\newcommand\eS{{\rm S}}

\newcommand\g{\mbox{\EUL g}}
\newcommand\hl{\rule[-.48em]{0em}{1.5em}}
\newcommand\h{\mbox{\EUL h}}
\newcommand\id{{\rm id}}

\newcommand\m{\mbox{\EUL m}}
\newcommand\pp{\mbox{\EUL p}}
\newcommand\pr{{\rm pr}}
\newcommand\p{\varrho}

\newcommand\rk{{\rm rk}}
\newcommand\so{\mathfrak s \mathfrak o}

\newcommand\spin{\mathfrak s \mathfrak p \mathfrak i \mathfrak n}
\newcommand\str{\rule[-.48em]{0em}{1.7em}}
\newcommand\strh{\rule[-.7em]{0em}{2.2em}}
\newcommand\stru{\rule[-.65em]{0em}{1.7em}}
\newcommand\strs{\rule[-.48em]{0em}{1.7em}}
\newcommand\suxu[2]{\mathfrak s(\mathfrak u({#1}) {+} \mathfrak u({#2}))}

\newcommand\s{\sigma}

\newcommand\ts{\textstyle}
\newcommand\x{{\times}}

\renewcommand\H{\mathbb{H}}
\renewcommand\O{{\rm O}}
\renewcommand\P{{\rm P}}
\renewcommand\a{\alpha}
\renewcommand\aa{{\mathfrak a}}

\renewcommand\eS{{\rm S}}
\renewcommand\k{\mbox{\EUL k}}

\renewcommand\ss{\mathfrak s}
\renewcommand\t{\tau}

\begin{document}
\setcounter{tocdepth}{1}

%
%
%

\pagenumbering{arabic}\pagestyle{plain}

\title{Polar actions on symmetric spaces}

\author{Andreas Kollross}

\address{Institut f\"{u}r Mathematik\\Universit\"{a}t Augsburg\\86135
Augsburg\\Germany} \email{kollross@math.uni-augsburg.de}




\begin{abstract}
We study isometric Lie group actions on symmetric spaces admitting a section, i.e.\ a submanifold which
meets all orbits orthogonally at every intersection point. We classify such actions on the compact
symmetric spaces with simple isometry group and rank greater than one. In particular we show that these
actions are hyperpolar, i.e.\ the sections are flat.
\end{abstract}

\maketitle


\section{Introduction and main results}
\pagenumbering{arabic}

An isometric action of a compact Lie group on a Riemannian manifold is called {\em polar} if there
exists a connected immersed submanifold~$\Sigma$ which intersects the orbits orthogonally and meets
every orbit. Such a submanifold~$\Sigma$ is then called a {\em section} of the group action. If the
section is flat in the induced metric, the action is called {\em hyperpolar}. Our main result is a
classification of polar actions on compact symmetric spaces with simple isometry group and rank greater
that one. This classification shows that these actions are in fact all hyperpolar.


One may think of the elements in a section as being {\em canonical forms}, representing the orbits of
the group action uniquely up to the action of a finite group, the {\em Weyl group}. This point of view
may be illustrated by the example of the orthogonal group~$\O(n)$ acting on the space of real symmetric
$n \times n$-matrices by conjugation, where the subspace of diagonal matrices is a section. Another
motivation comes from submanifold geometry, in particular from the theory of isoparametric submanifolds
and their generalizations~\cite{terng}, \cite{tt1}. The orbits of polar actions have many remarkable
geometric properties, for instance, the principal orbits of polar representations are isoparametric
submanifolds of Euclidean space.


However, the history of the subject probably starts with an application in topo\-logy. Bott~\cite{bott}
and Bott and Samelson~\cite{bs} considered the adjoint action of a compact Lie group on itself and on
its Lie algebra~\cite{bott} and more generally, the isotropy action of a compact symmetric
space~\cite{bs}. The motivation of Bott and Samelson to consider these actions was that they are
``variationally complete'', which made it possible to apply Morse theory to the space of loops in the
symmetric space. Conlon~\cite{conlon1} proved that hyperpolar actions on Riemannian manifolds are
variationally complete, referring to the sections as {\em $K$-transversal domains}.
Hermann~\cite{hermann} found another class of examples, namely if $H$ and $K$ are both symmetric
subgroups of a simple compact Lie group~$G$, then the action of $H$ on the symmetric space~$G / K$ is
hyperpolar. It was shown much later~\cite{gt} that actions on compact symmetric spaces are variationally
complete if and only if they are hyperpolar.

Conlon~\cite{conlon2} observed that s-representations are hyperpolar and later on Dadok~\cite{dadok}
obtained a classification of irreducible polar representations. The classification shows that the
connected components of the orbits of a polar representation agree with the orbits of an
s-representation after a suitable identification of the representation spaces. Reducible polar
representations were classified by Bergmann~\cite{bergmann}.

Cohomogeneity one actions, i.e.\ actions whose principal orbits are hypersurfaces, are a special case of
independent interest. Cohomogeneity one actions on spheres were classified by Hsiang and
Lawson~\cite{hsl}. Later Takagi~\cite{takagi}, D'Atri~\cite{datri}, and Iwata~\cite{iwata} classified
cohomogeneity one actions on $\C \P^n$, $\H \P^n$ and $\Ca \P^2$, respectively.

Szenthe~\cite{szenthe}, Palais and Terng~\cite{pt} investigated fundamental properties of polar actions
on Riemannian manifolds. Heintze, Palais, Terng and Thorbergsson~\cite{hptt}, \cite{hptt2} obtained
structural results for hyperpolar actions on compact symmetric spaces, studied relations to polar
actions on infinite dimensional Hilbert space and involutions of affine Kac-Moody algebras. They showed
in particular~\cite{hptt2} that compact Riemannian homogeneous spaces admitting a hyperpolar action with
a fixed point are symmetric.

In~\cite{kollross}, the author gave a classification of hyperpolar actions on the irreducible compact
symmetric spaces, the main result being that these actions are orbit equivalent to the examples found by
Hermann if the cohomogeneity is~$\ge 2$.

Podest\`{a} and Thorbergsson~\cite{pth1} classified polar actions on the compact symmetric spaces of
rank one. The first result on polar actions on irreducible symmetric spaces of higher rank without
assuming flatness of the sections was obtained by Br{\"u}ck~\cite{brueck}, who showed that polar actions
with a fixed point on these spaces are hyperpolar. Podest\`{a} and Thorbergsson~\cite{pth2} proved that
polar actions on compact irreducible homogeneous K{\"a}hler manifolds are co\-iso\-tro\-pic and classified
co\-iso\-tro\-pic and polar actions on the real Grassmannians~$\G_2(\R^n)$ of rank two. It turned out
that all polar actions on these spaces are hyperpolar.

This approach was further pursued by Biliotti and Gori~\cite{bg}, who classified co\-iso\-tro\-pic and
polar actions on the complex Grassmannians~$\G_k(\C^n)$. The classification of co\-iso\-tro\-pic actions
on the compact irreducible Hermitian symmetric spaces was recently completed by
Biliotti~\cite{biliotti}, showing in particular that polar actions on these spaces are hyperpolar, which
led Biliotti to conjecture that this holds for all compact irreducible symmetric spaces.

The present work extends the classification of polar actions to all irreducible symmetric spaces of
type~I, i.e.\ to the compact symmetric spaces with simple isometry group, confirming the conjecture of
Biliotti for these spaces. We show that polar actions on the symmetric spaces of type~I and higher rank
are hyperpolar. That is, they are of cohomogeneity one or orbit equivalent to the examples found by
Hermann. Our main result can be stated as follows.

\begin{maintheorem}\label{MainTheorem}
Let $M$ be a compact symmetric space of rank greater than one whose isometry group $G$ is simple. Let $H
\subset G$ be a closed connected  non-trivial subgroup acting polarly on~$M$. Then the action of $H$ on
$M$ is hyperpolar, that is, the sections are flat in the induced metric. Moreover, the sections are
embedded submanifolds.
\end{maintheorem}

\begin{table}[h]\rm
\begin{tabular}{|c|c|c|}
\hline \str Type & $\pi^{-1}(H)$ & $\tilde M$ \\ \hline \hline

A\,III-II & $\SU(2n{-}2k{-}1)\times\SU(2k{+}1)$ & $\SU(2n) / \Sp(n)$  \\
\hline

A\,III-II & $\eS(\U(2n-2) \times \U(1) \times \U(1))$ & $\SU(2n) / \Sp(n)$  \\
\hline

A\,III-II & $\eS(\U(2n-2) \times \U(1))$ & $\SU(2n) / \Sp(n)$  \\
\hline

A\,III-III & $\SU(k) \times \SU(n-k)$ & $\SU(n) /
\SUxU{\ell}{n-\ell},\,(k,\ell) \neq \left(\frac{n}{2},\frac{n}{2}\right)$ \\
\hline \hline

BD\,I-I& $\LG_{2} \times \SO(n-7)$ & $\SO(n) / \SO(2) \times \SO(n-2),\,n\ge7$ \\ \hline

BD\,I-I& $\Spin(7) \times \SO(n-8)$ & $\SO(n) / \SO(2) \times \SO(n-2),\,n\ge8$ \\ \hline

BD\,I-I& $\Spin(7) \times \SO(n-8)$ & $\SO(n) / \SO(3) \times \SO(n-3),\,n\ge8$ \\ \hline

BD\,I-I& $\LG_{2} \times \LG_{2}$ & $\SO(14) / \SO(2) \times \SO(12)$ \\
\hline

BD\,I-I& $\LG_{2} \times \Spin(7)$ & $\SO(15) / \SO(2) \times \SO(13)$ \\
\hline

BD\,I-I& $\Spin(7) \times \Spin(7)$ & $\SO(16) / \SO(2) \times \SO(14)$
\\ \hline

BD\,I-I& $\Spin(7) \times \Spin(7)$ & $\SO(16) / \SO(3) \times \SO(13)$
\\ \hline\hline

C\,I-II& $\SU(n)$ & $\Sp(n) / \Sp(k) \times \Sp(n - k)$ \\
\hline

C\,II-II& $\Sp(a + b)$ & $\Sp(a + b + 1) / \Sp(a) \times \Sp(b + 1)$ \\
\hline

C\,II-II& $\Sp(a + b) \times \U(1)$ & $\Sp(a + b + 1) / \Sp(a) \times \Sp(b + 1)$ \\
\hline \hline

D$_4$\,I-I'& $\LG_2$ & $\SO(8) / \U(4)$ \\
\hline

D\,I-III & $\SO(2n-2)$ & $\SO(2n) / \U(n)$ \\ \hline

D\,III-I & $\SU(n)$ & $\SO(2n) / \SO(k) \times \SO(2n-k),\, k < n$
\\ \hline

D\,III-III' & $\SU(n)$ & $\SO(2n) / \a(\U(n))$ \\ \hline \hline

E\,II-IV& $\SU(6)$ & $\LE_6 / \LF_4$ \\ \hline

E\,II-IV& $\SU(6) \cdot \U(1)$ & $\LE_6 / \LF_4$ \\ \hline

E\,III-II& $\Spin(10)$ & $\LE_6 / \SU(6) \cdot \Sp(1)$ \\ \hline

E\,VII-VI& $\LE_6$ & $\LE_7 / \SO'(12) \cdot \Sp(1)$ \\ \hline

\end{tabular}
\bl\caption{Hyperpolar subactions of Hermann actions}\label{TOrbEqSubgrHermann}
\end{table}

In \cite{kollross}, the hyperpolar actions on irreducible compact symmetric spaces were only determined
up to orbit equivalence. In the present work we obtain the complete classification of connected Lie
groups acting polarly without fixed points on the symmetric spaces of higher rank with simple isometry
group. For actions with fixed points  the complete classification follows immediately from
Corollary~\ref{PolFixOnHomSp} and Lemma~\ref{OrbitEqPolarSubGr}.

\begin{maintheorem}\label{MainTheorem2}
Let $M = G / K$ be a connected compact symmetric space of rank greater than one whose isometry group is
simple. Let $H \subset G$ be a closed connected proper subgroup such that the $H$-action on $G / K$ is
polar, non-trivial, non-transitive, and without fixed point. Then
\begin{enumerate}

\item either $H \subset G$ is maximal connected (and as described
in Theorem~A of~\cite{kollross})

\item or the universal cover of the symmetric space~$\tilde M$ and
the conjugacy class of the subgroup~$H \subset G$ are as given by Table~\ref{TOrbEqSubgrHermann}, where
$\pi \colon \Isom(\tilde M) \to G$ is the covering map, and there exists a connected subgroup $H_0
\subset G$ whose Lie algebra $\h_0 \subset \g$ is the fixed point set of an involution of~$\g$ and such
that the $H_0$-action on~$G / K$ has the same orbits as the $H$-action.
\end{enumerate}
\end{maintheorem}

The first column of Table~\ref{TOrbEqSubgrHermann} indicates a connected subgroup~$H_0$ of~$G$
containing~$H$, see~Table~\ref{THermannActions} and the remarks there. By $\alpha$ we denote a
non-trivial outer automorphism of~$\SO(2n)$ of order two given by conjugation with an element
from~$\O(2n) \setminus \SO(2n)$. The proofs of Theorems~\ref{MainTheorem} and \ref{MainTheorem2} are
completed and summarized on pages~\pageref{MainProofBegin}--\pageref{MainProofEnd}.

Combining Theorem~\ref{MainTheorem} with the results of~\cite{pth1} and Corollary~D of~\cite{hl}, we
obtain the following result on sections and Weyl groups of polar actions. The Weyl group $W_{\Sigma} =
\N_H(\Sigma) / Z_H(\Sigma)$ is a quotient group of the group $\hat W_{\Sigma}$ as defined in
Lemma~\ref{TotGeodHyperWeyl}.

\begin{maincorollary}\label{MainCorollary2}
Let $H$ be a connected compact Lie group acting polarly on a compact symmetric space~$M$ with simple
isometry group. Then a section~$\Sigma$ of the $H$-action on~$M$ is isometric to a flat torus, a sphere
or a real projective space. The group $\hat W_{\Sigma}$ acting on the universal cover of~$\Sigma$ is an
irreducible affine Coxeter group in case $\Sigma$ is flat or a finite Coxeter group of Euclidean space
restricted to a sphere in case $\Sigma$ is non-flat.
\end{maincorollary}

In particular, the Weyl groups of such polar actions can be described by connected Dynkin diagrams of
affine type (in the hyperpolar case) or Dynkin diagrams of the finite type (in the polar, non-hyperpolar
case).

This article is organized as follows. We start by setting up terminology and notation. We then review
examples and known results on polar actions. In Section~\ref{TotGeodSubmf} we recall some facts about
symmetric spaces and their totally geodesic submanifolds; in particular, we give a characterization of
maximal totally geodesic submanifolds and obtain an upper bound on the dimension of totally geodesic
submanifolds locally isometric to a product of spheres.

In Section~\ref{CriteriaPolarity}, we recall a criterion which reduces the problem of deciding whether
an action on a symmetric spaces is polar or not to a problem on the Lie algebra level. In
Section~\ref{SectionsWeylAct} we prove the Splitting Theorem~\ref{SplittingTheorem} which says that if a
section~$\Sigma$ of a polar action admits a local splitting $\tilde \Sigma = \tilde \Sigma_1 \times
\tilde \Sigma_2$ such that the Weyl group acts trivially on one factor~$\tilde \Sigma_2$, then the
symmetric space is locally a Riemannian product~$M \times \Tilde \Sigma_2$. As a consequence, we show
that the section of a polar action on a compact irreducible symmetric space is locally isometric to a
product of spaces of constant curvature. This observation is crucial for our classification since it
implies an upper bound on the cohomogeneity, reducing the classification problem to a finite number of
cases.

In Section~\ref{PolSub} we introduce another main tool by collecting various sufficient conditions for
actions to be {\em polarity minimal}, which means that the restriction to a closed connected subgroup
with orbits of lower dimension is either non-polar or trivial. This is of essential importance since it
enables us to restrict our attention at first to maximal subgroups of the isometry group. In many cases
we are able to show that the action of a maximal connected subgroup is non-polar and polarity minimal,
thereby excluding all of its subgroups.

In the remaining part of the paper, the classification is carried out. We start with the maximal
connected subgroups in the isometry group of a symmetric space. In Section~\ref{SubHermannHighRk}, we
consider Hermann actions, i.e.\ actions of symmetric subgroups of the isometry group. We show that
actions of cohomogeneity~$\ge 2$ are polarity minimal and determine orbit equivalent subactions. We then
consider maximal connected subgroups in the isometry group of classical symmetric spaces which are given
by irreducible representations of non-simple groups. It turns out that they are either non-polar and
polarity minimal or of cohomogeneity one. In Section~\ref{SubSimple} we study actions of simple
irreducible subgroups in the classical groups. In Section~\ref{ClassExcept} we consider actions on the
exceptional symmetric spaces. It turns out that the actions of non-symmetric maximal subgroups are
non-polar and polarity minimal.

It then remains to study subactions of cohomogeneity one and transitive actions. Since we do not have an
{\em a priori} proof that these actions are polarity minimal, it is necessary to descend from maximal
connected subgroups~$H_1 \subset G$ acting with cohomogeneity~$\le 1$ to further subgroups~$$G \supset
H_1 \supset H_2 \supset \ldots,$$ where $H_{n+1} \subset H_n$ is maximal connected, until we arrive at
an action which is polarity minimal.

\paragraph{\em Acknowledgements} I would like to thank
Ernst Heintze, Mamoru Mimura, Chuu-Lian Terng and Gudlaugur Thorbergsson for helpful discussions and
comments; I am especially indebted to Burkhard Wilking for providing a crucial step in the proof of the
Splitting Theorem~\ref{SplittingTheorem} and for pointing out an error in an earlier version of this
paper.

\hfill\\
\begin{center}\begin{minipage}{33em}{\small
\tableofcontents}\end{minipage}\end{center}


\section{Preliminaries and examples}
\label{PrelimEg}

An isometric action of a compact Lie group~$G$ on a Riemannian manifold is called {\em polar} if there
exists a connected immersed submanifold~$\Sigma$ such that $\Sigma$ meets all $G$-orbits and the
intersection of~$\Sigma$ with any $G$-orbit is orthogonal at all intersection points. Such a submanifold
$\Sigma$ is called a {\em section} for the $G$-action on~$M$. In particular, the actions of finite
groups and transitive actions are special cases of polar actions, the section being the whole space or a
point, respectively.

Note that we do not require the section to be an embedded submanifold, generalizing the definitions
of~\cite{pt} and \cite{conlon1}. However, it turns out by our classification that on symmetric spaces of
type~I the sections are flat and therefore closed embedded submanifolds by Corollary~2.12
of~\cite{hptt}.
\par
The dimension of $\Sigma$ equals the cohomogeneity of the $G$-action and hence the tangent space $\T_p
\Sigma$ at a regular point~$p \in \Sigma$ coincides with the normal space~$\N_p (G \cdot p)$ at $p$ of
the orbit through~$p$. From this, it follows that any two sections are mapped isometrically onto each
other by some group element. It has been proved in~\cite{pt} that sections are totally geodesic.
\par
In the special case where the sections are flat in the induced metric, the action is called {\em
hyperpolar}. Examples for hyperpolar actions are given by the action of a compact Lie group on itself by
conjugation, where the sections are the maximal tori. More generally, the action of an isotropy group of
a symmetric space is hyperpolar, the sections being the flats of the symmetric space.
\par
For a polar action one can define the Weyl group by considering the normalizer of a section, i.e.\ all
group elements which map the section onto itself, this group acts on the section by isometries and the
Weyl group is defined by factoring out the kernel of this action.

\begin{definition}
Let $M$ be a Riemannian manifold on which the compact Lie group $G$ acts polarly with section $\Sigma$.
The {\em (generalized) Weyl group} $W_{\Sigma} = W_{\Sigma}(M,G)$ is the group
$N_G(\Sigma)/Z_G(\Sigma)$, where $N_G(\Sigma) = \{ g \in G \mid g \cdot \Sigma  = \Sigma \}$ and
$Z_G(\Sigma) = \{ g \in G \mid g \cdot s  = s \mbox{ for all} s \in \Sigma \}$ are the {\em normalizer}
and {\em centralizer} of $\Sigma$ in $G$, respectively.
\end{definition}

Two Riemannian $G$-manifolds are called {\em conjugate} if there exists an equivariant isometry between
them. In particular, the actions of two conjugate subgroups of the isometry group of a Riemannian
manifold are conjugate. To study isometric actions on a Riemannian manifold it suffices to consider
conjugacy classes of subgroups in the isometry group.

\par

Two isometric actions of two Lie groups $G$ and $G'$ on a Riemannian manifold $M$ are called {\em orbit
equivalent} if there exists an isometry of $M$ which maps $G$-orbits onto $G'$-orbits; they are called
{\em locally orbit equivalent} if there is an isometry mapping connected components of $G$-orbits onto
connected components of $G'$-orbits. Obvious examples of orbit equivalent actions are given by various
groups acting transitively on spheres, e.g.\ the actions of $\SO(4n)$, $\SU(2n)$, $\U(2n)$, and $\Sp(n)$
on $\R^{4n}$ are all orbit equivalent.

\par

We use the term {\em subaction} for the restriction of an action of a group~$G$ to a subgroup $H
\subseteq G$; in case the $H$-orbits coincide with the $G$-orbits, the $H$-action is called {\em orbit
equivalent subaction}.

\par

A normal subgroup $N$ of a compact Lie group $G = G' \cdot N$ acting isometrically on a Riemannian
manifold is called {\em inessential} if the $G$-action restricted to~$G'$ is orbit equivalent to the
$G$-action. An isometric action of a compact connected Lie group $G$ on a Riemannian manifold $M$ is
called {\em orbit maximal} if any other isometric action of any other compact connected Lie group $G'$
such that every $G$-orbit is contained in a $G'$-orbit is either orbit equivalent or transitive on $M$.

An immersed submanifold $M$ in a symmetric space $N$ is said to have {\em parallel focal structure} if
the normal bundle $\nu(M)$ is globally flat and the focal data is invariant under normal parallel
translation, that is, for every parallel normal field $v$ on $M$ the rank of $d\eta_{v(x)}$ is locally
constant on $M$, where the end point map $\eta \colon \N M \rightarrow N, v \mapsto \exp(v)$ is defined
to be the restriction of the exponential map to the normal bundle $\nu(M)$, see~\cite{tt1}. The
principal orbits of a polar action on a symmetric space have parallel focal structure \cite{ewert}.
\par
A submanifold with parallel focal structure is called {\em equifocal} if the normal bundle $\nu(M)$ is
{\em abelian}, that is, $\exp(\nu(M))$ is contained in some totally geodesic flat subspace of~$N$ for
each point $x \in M$. The principal orbits of hyperpolar actions on symmetric spaces of compact type are
equifocal submanifolds, see~\cite{tt1}, Theorem 2.1. Our Theorem~\ref{MainTheorem} shows that
submanifolds with parallel focal structure which arise as principal orbits of polar actions on symmetric
spaces of higher rank with simple compact isometry group are in fact equifocal. We conjecture that more
generally submanifolds with parallel focal structure in irreducible compact symmetric spaces of higher
rank are equifocal, hence of codimension one or homogeneous by the result of Christ~\cite{christ}.

\subsection{Notation}

We will frequently use the following notational conventions for compact Lie groups and their
representations. We view the classical Lie groups~$\SO(n)$, $\SU(n)$, and $\Sp(n)$ as matrix Lie groups
as described in~\cite{helgason}, Ch.~X, \S~2.1. We assume that reducible subgroups of the classical
groups are standardly embedded, e.g.\ by $\SO(m)\times\SO(n)$ we denote the subgroup
\begin{equation}\label{StdEmbedd}
\left\{\left.\left(\begin{array}{c|c}A&\\\hline&B\end{array}\right)
 \right|A\in\SO(m),\,B\in\SO(n)\right\}\subset\SO(m+n).
\end{equation}
We write $H_1 \otimes H_2$ for the Kronecker product of two matrix Lie groups. When we write $\LG_2$, we
refer to an irreducible representation by orthogonal $7 \times 7$-matrices; similarly, $\Spin(7)$ stands
for a matrix Lie group which is the image of the $8$-dimensional spin representation of~$\Spin(7)$. By
$\R^{n}$, $\C^{n}$, $\H^{n}$ we will denote the standard representation of $\O(n)$, $\U(n)$ or $\Sp(n)$,
respectively.


\subsection{Polar representations}\label{PolRep}

An important class of examples for polar actions is given by polar representations on Euclidean space.
Since the sections of polar actions are totally geodesic, they are linear subspaces in the case of polar
representations and polar representations are therefore automatically hyperpolar. Polar representations
are of importance for our classification since they occur as slice representations of polar actions. Let
$M$ be a Riemannian $G$-manifold and let $G_p$ be the isotropy subgroup at $p$. The restriction of the
isotropy representation to $\N_p(G\cdot{p})$ is called the {\it slice representation} at $p$. Slice
representations are a fundamental tool for the study of Lie group actions since they provide a means to
describe the local behavior of an action in a tubular neighborhood of an orbit by a linear
representation.

\begin{slicetheorem}\label{SliceTheorem}
Let $M$ be a Riemannian $G$-manifold, let $p \in M$ and $V = \N_p(G \cdot p)$ the normal space at $p$ to
the $G$-orbit through $p$. Then there is an equivariant diffeomorphism $\Psi$ of a $G$-invariant open
neighborhood around the zero section in the normal bundle $G \times_{G_p} V \rightarrow G/G_p$ onto a
$G$-invariant open neighborhood around the orbit $G\cdot p$ such that the zero section in $G
\times_{G_p} V$ is mapped to the orbit $G\cdot p$. The diffeomorphism $\Psi$ is given by the end point
map which maps any normal vector $v_q \in N_q(G \cdot p)$ to its image under the exponential map
$\exp_q(v_q)$.
\end{slicetheorem}


\begin{proof}
See e.g.~\cite{jaenich}, p.~3.
\end{proof}


It is an immediate consequence of the Slice Theorem~\ref{SliceTheorem} that the slice representation and
the $G$-action on $M$ have the same cohomogeneity.  Slice representations are in particular useful for
our classification since the polarity of an action is inherited by its slice representations.
\begin{proposition}\label{PolarSlice}
Let $M$ be a Riemannian $G$-manifold. If the action on $N$ is polar then for all $p \in M$ the slice
representation at $p$ is polar with $\T_p\Sigma$ as a section, where $\Sigma$ is the section of the
$G$-action on $M$ containing $p$.
\end{proposition}

\begin{proof}
This was proved in~\cite{pt}, Theorem~4.6. Although in \cite{pt} the sections are assumed to be embedded
submanifolds, the proof is still valid if one requires the sections only to be immersed.
\end{proof}


We use the term {\em effectivized slice representation} to describe the representation of the isotropy
group with the effectivity kernel factored out. Let us recall some known results about polar
representations.

\begin{definition}
Let $G$ be a compact Lie group and $K$ be a closed subgroup. By
$$\chi(G,K) = \Ad_G|_K \ominus \Ad_K$$ we denote the equivalence class
of the isotropy representation of the homogeneous space $G/K$, i.e.\ the restriction of the adjoint
representation of~$G$ to~$K$ acting on a $K$-invariant complement of $\k$ in $\g$. In the special case
of a symmetric pair~$(G,K)$, see below, the (equivalence class of the) representation~$\chi(G/K)$ is
called an {\em s-representation}.
\end{definition}

A compact subgroup~$K$ of a Lie group~$G$ is called {\em symmetric subgroup} if there exists an
involutive automorphism of~$G$ such that $G_0^{\s} \subseteq K \subseteq G^{\s}$, where $G^{\s}$ and
$G_0^{\s}$ denote the fixed point set of~$\s$ and its connected component, respectively. A pair $(G,K)$,
where $G$ is a Lie group and $K$ a symmetric subgroup is called a {\em symmetric pair}. Any Riemannian
globally symmetric space~$M$ has a homogeneous presentation $G /K$, where $G$ is the isometry group
of~$M$, such that $K$ is a symmetric subgroup of~$G$. Conversely, if $(G,K)$ is a symmetric pair, then
$G / K$ endowed with a $G$-invariant metric is a Riemannian globally symmetric space,
see~\cite{helgason}.

It is well known that the adjoint representations of compact Lie groups and more generally
s-representations are polar. As far as concerns the geometry of the orbits, also the converse is true.


\begin{theorem}[Dadok]\label{DadoksTheorem}
A representation $\rho \colon G \to \O(n)$ of a compact Lie group is polar if and only if it is locally
orbit equivalent to an s-representation, i.e.\ the connected components of its orbit agree with the
orbits of an s-representation after a suitable isometric identification of the representation spaces.
\end{theorem}


\begin{proof}
The proof given in \cite{dadok} relies on a classification of the irreducible polar representations.
See~\cite{eh1} for a conceptual proof in case the cohomogeneity is $\ge 3$.  See~\cite{lowcoh} for an
alternative proof, where a similar classification strategy as in the present work is used.
\end{proof}

It is shown in Theorem~3.12 of \cite{hptt} that irreducible polar representations of cohomogeneity~$\ge
2$ are orbit maximal when restricted to a sphere around the origin. For irreducible s-representations of
cohomogeneity $\ge 2$, orbit equi\-valent subgroups were determined in~\cite{eh2}. We state the result
below.


\begin{lemma}\label{OrbitEqPolarSubGr}
Let $G$ be a connected simple compact Lie group and let $K$ be a connected symmetric subgroup such that
$\rk(G / K) \ge 2$.  Let $H \subseteq K$ be a closed connected subgroup. Then $\chi(G,K)$ and
$\chi(G,K)|_H$ are orbit equivalent if and only if either $H = K$ or the triple $(G,K,H)$ is as given in
Table~\ref{TOrbitEqPolarSubGr}.
\begin{table}[h]\rm
  \begin{tabular}{|c|c|c|c|}\hline
    \hl \str $G$ & $K$ & $H$ & Range \\\hline\hline
    $\SO(9)$ & $\SO(2)\times \SO(7)$ & $\SO(2)\times \LG_{2}$ & \\\hline
    $\SO(10)$ & $\SO(2)\times \SO(8)$ & $\SO(2)\times \Spin(7)$ & \\\hline
    $\SO(11)$ & $\SO(3)\times \SO(8)$ & $\SO(3)\times \Spin(7)$ & \\\hline
    $\SU(p+q)$ & $\eS(\U(p)\times \U(q))$ & $\SU(p)\times \SU(q)$ & $p \neq q$ \\\hline
    $\SO(2n)$ & $\U(n)$ & $\SU(n)$ & $n$ odd \\\hline
    $\LE_6$ & $\U(1)\cdot \Spin(10)$ & $\Spin(10)$ & \\\hline
  \end{tabular}
\bl\caption{Orbit equivalent subactions of polar representations} \label{TOrbitEqPolarSubGr}
\end{table}
\end{lemma}


\subsection{Hyperpolar actions on symmetric spaces}

If $H$, $K$ are two symmetric subgroups of the compact Lie group~$G$, then the action of $H$ on $G / K$
is hyperpolar~\cite{hermann}. Slightly more generally, if $H$ is a subgroup of $G$ such that its Lie
algebra $\h \subset \g$ is the fixed point set of an involution of~$\g$, then the action of $H$ on the
symmetric space $G / K$ is hyperpolar and we call such actions {\em Hermann actions}. In the special
case $H = K$ we have the isotropy action of the symmetric space and the sections are just the flats of
the symmetric space. It was shown in~\cite{kollross} that all hyperpolar actions on irreducible
symmetric spaces of compact type are of cohomogeneity one or orbit equivalent to Hermann actions.

All fixed-point free Hermann actions on the symmetric spaces of type~I are given by
Table~\ref{THermannActions}. Here $\a$ denotes the non-trivial diagram automorphism of $\SO(2n)$ given
by conjugation with a matrix from~$\O(2n) \setminus \SO(2n)$ and $\tau$ an order three diagram
automorphism of $\Spin(8)$. The type of the Hermann action indicated in the first column refers to the
type of the symmetric subgroups involved as given in Table~\ref{TSymmSpaces}, e.g.\ the symbol A\,I-II
refers to the action of $H$ on $G /K$, where $G / H = \SU(2n) / \SO(2n)$ is a symmetric space of
type~A\,I and $G / K = \SU(2n) / \Sp(n)$ is a symmetric space of type~A\,II; whereas for the action $K$
on $G / H$ we use the notation~A\,II-I. For the conjugacy classes of connected symmetric subgroups in
simple compact Lie groups see~\cite{kollross}, 3.1.1 and 3.1.2. The cohomogeneity of the actions is
given in the last column.


\begin{table}[h]
\begin{tabular}{|l*{4}{|c}|}
\hline \hl \str Type & $H$ & $G$ & $K$ & Coh.\  \\\hline \hline \hl
A\,I-II & $\SO(2n)$ & $\SU(2n)$ & $\Sp(n)$ & $n-1$ \\ \hline \hl
A\,I-III ($k\le\lfloor\frac{n}{2}\rfloor$) & $\SO(n)$ & $\SU(n)$ & $\SUxU{k}{n-k}$ & $k$   \\ \hline \hl
A\,II-III ($k \le n)$& $\Sp(n)$ & $\SU(2n)$ & $\eS(\U(k){\times}\U(2n-k))$ & $\left\lfloor\frac{k}{2}\right\rfloor$ \\
\hline \hl
A\,III-III ($k\le\ell\le\lfloor\frac{n}{2}\rfloor$)& $\eS(\U(k){\times}\U(n{-}k))$ & $\SU(n)$ &
$\eS(\U(\ell){\times}\U(n{-}\ell))$ & $k$ \\ \hline \hl
BD\,I- I ($k\le\ell\le\lfloor\frac{n}{2}\rfloor$) & $\SO(k){\times}\SO(n{-}k)$ & $\SO(n)$ &
$\SO(\ell){\times}\SO(n{-}\ell)$ & $k$
\\ \hline \hl
C\,I-II ($k\le\lfloor\frac{n}{2}\rfloor$) & $\U(n)$ & $\Sp(n)$ & $\Sp(k){\times}\Sp(n{-}k)$ & $k$ \\ \hline \hl
C\,II-II ($k\le\ell\le\lfloor\frac{n}{2}\lfloor$)  & $\Sp(k){\times}\Sp(n{-}k)$ & $\Sp(n)$ & $\Sp(\ell){\times}\Sp(n{-}\ell)$ & $k$ \\
\hline \hl
D\,I-III ($k \le n$)  & $\SO(k){\times}\SO(2n{-}k)$ & $\SO(2n)$ & $\U(n)$ &
$\left\lfloor\frac{k}{2}\right\rfloor$
\\ \hline \hl
D\,III-III' & $\U(2n)$ & $\SO(4n)$ & $\a(\U(2n))$ & $n-1$
\\ \hline \hl
D$_4$ I-I' ($k{\le}\ell{\le}3$) & $\Spin(k){\cdot}\Spin(8{-}k)$ & $\Spin(8)$ & $\tau(\Spin(\ell){\cdot}\Spin(8{-}\ell))$ & $k-1$ \\
\hline \hl
E\,I-II & $\Sp(4)/\{\pm1\}$&$\LE_6$&$\SU(6){\cdot}\Sp(1)$&$4$\\
\hline \hl
E\,I-III& $\Sp(4)/\{\pm1\}$&$\LE_6$&$\Spin(10){\cdot}\U(1)$&$2$\\
\hline \hl
E\,I-IV & $\Sp(4)/\{\pm1\}$&$\LE_6$&$\LF_4$&$2$\\ \hline \hl
E\,II-III & $\SU(6){\cdot}\Sp(1)$&$\LE_6$&$\Spin(10){\cdot}\U(1)$&$2$\\
\hline \hl
E\,II-IV & $\SU(6){\cdot}\Sp(1)$&$\LE_6$&$\LF_4$&$1$\\ \hline \hl
E\,III-IV & $\Spin(10){\cdot}\U(1)$&$\LE_6$&$\LF_4$&$1$\\ \hline \hl
E\,V-VI & $\SU(8)/\{\pm1\}$&$\LE_7$&$\SO'(12){\cdot}\Sp(1)$&$4$\\
\hline \hl
E\,V-VII & $\SU(8)/\{\pm1\}$&$\LE_7$&$\LE_6{\cdot}\U(1)$&$3$\\
\hline \hl
E\,VI-VII & $\SO'(12){\cdot}\Sp(1)$&$\LE_7$&$\LE_6{\cdot}\U(1)$&$2$\\ \hline \hl
E\,VIII-IX & $\SO'(16)$&$\LE_8$&$\LE_7{\cdot}\Sp(1)$&$4$\\\hline \hl
F\,I-II & $\Sp(3){\cdot}\Sp(1)$&$\LF_4$&$\Spin(9)$&$1$\\ \hline
\end{tabular}
\bl\caption{Hermann actions}\label{THermannActions}
\end{table}


Hyperpolar actions on compact symmetric spaces have the remarkable property that they lift under certain
Riemannian submersions to actions which are again hyperpolar, cf.~\cite{hptt}.


\begin{proposition}\label{HyperPLiftToGroup}
Let $G$ be a compact simple Lie group and let $K \subset G$ be a symmetric subgroup, $M = G/K$ the
corresponding symmetric space and $H$ a closed subgroup of $G$. Then the $H$-action on $M$ is hyperpolar
if and only if the $H \times K$-action on $G$ is hyperpolar.
\end{proposition}


\begin{proof}
See \cite{hptt}, Proposition 2.11
\end{proof}


It can be shown using Proposition~\ref{PolCrit} that polar actions on compact symmetric spaces have this
lifting property only if they are hyperpolar. In particular, if we lift the known polar actions on
symmetric spaces to the groups, we do {\em not} obtain any examples of polar actions besides the
hyperpolar ones. Another remarkable property of hyperpolar actions is that they are orbit maximal on
irreducible symmetric spaces of compact type.


\begin{proposition}\label{HyperPMax}
Let $M = G/K$ be a connected irreducible symmetric space of compact type and $H \subset L \subset G$
closed connected subgroups. If the $H$-action on $M$ is hyperpolar, then the $L$-action on $M$ is
transitive or orbit equivalent to the $H$-action.
\end{proposition}
\begin{proof}
See \cite{hl}, Corollary D.
\end{proof}

The non-orbit maximal examples of polar actions found by Podest\`{a} and Thorbergsson~\cite{pth1} show
that Proposition~\ref{HyperPMax} does not directly generalize to polar actions, see below. However, it
is a consequence of our classification that polar actions on symmetric spaces of rank~$\ge 2$ with
simple compact isometry group are orbit maximal.


\subsection{Polar actions on rank one symmetric spaces}
Polar actions on rank one symmetric spaces have been classified by Podest{\`a} and Thorbergsson \cite{pth1}.
The hyperpolar, i.e.\ cohomogeneity one, actions on these spaces had before been classified in
\cite{hsl}, \cite{takagi}, \cite{datri} and \cite{iwata}.
The results can be summarized as follows. The classification of polar actions on spheres (and real
projective spaces) follows from \cite{dadok}, since every polar action on the sphere is given as the
restriction of a polar representation to the sphere. The isotropy representations of Hermitian symmetric
spaces of real dimension $2n+2$ induce polar actions on $\C\P^n$ and all polar actions on $\C\P^n$ are
orbit equivalent to actions obtained in this fashion. Similarly, all polar actions on $\H\P^n$ come from
isotropy representations of products of quaternion-K{\"a}hler symmetric spaces, with the additional
restriction that all factors but one must be of rank one. While all these polar actions arise from polar
actions on the sphere, the actions on the Cayley plane $\Ca \P^2 = \LF_4/\Spin(9)$ do not have such an
interpretation. The maximal connected subgroup $\SU(3)\cdot\SU(3)\subset\LF_4$ acts polarly on the
Cayley plane with cohomogeneity two. The groups $\Sp(3) \cdot \Sp(1)$, $\Sp(3) \cdot \U(1)$, $\Sp(3)$
and $\Spin(9)$ act with cohomogeneity one. In addition there are three polar actions of cohomogeneity
two with a fixed point of the following subgroups of $\Spin(9)$:
$$
  \Spin(8),
 \qquad \SO(2)\cdot\Spin(7),
 \qquad \Spin(3)\cdot\Spin(6).
$$
In particular, polar actions on rank one symmetric spaces are not orbit maximal in general.


\section{Symmetric spaces and their totally geodesic submanifolds}\label{TotGeodSubmf}

In the following we will collect some useful facts about symmetric spaces and their totally geodesic
submanifolds. Sections of polar actions are totally geodesic submanifolds and it will be shown in
Theorem~\ref{ProductOfSpheres} that the sections of a non-trivial polar action on an irreducible compact
symmetric space are locally isometric to Riemannian products whose factors are spaces of constant
curvature. We give an upper bound on the dimension of such submanifolds in Lemma~\ref{ProdSphDim}. For
the proof of Theorem~\ref{ProductOfSpheres}, which is essentially a consequence of the Splitting
Theorem~\ref{SplittingTheorem}, we will need the characterization of totally geodesic hypersurfaces in
reducible symmetric spaces given in Corollary~\ref{TotGeodHyperSurface}, because the Weyl group of a
polar action is generated by reflections in totally geodesic hypersurfaces. We will conclude this
section by recalling a well known characterization of maximal subgroups in the classical groups.


Every symmetric space $M$ may be presented as $G/K$, where $G$ is the isometry group of $M$ and $K$ is a
symmetric subgroup of $G$. Conversely, if $(G,\,K)$ is a symmetric pair, then $G/K$ is a symmetric space
if it is equipped with an appropriate metric. A Riemannian manifold is an irreducible Riemannian
symmetric space of compact type if and only if it is isometric to
\begin{itemize}

\item either $G/K$, where $G$ is a simple, compact, connected Lie
group and $K$ a symmetric subgroup of G {\it (symmetric space of type~I)}

\item or a simple, compact, connected Lie group equipped with a
biinvariant metric {\em (symmetric space of type~II)}

\end{itemize}

The local isometry classes of the symmetric spaces of type~I are given by Table~\ref{TSymmSpaces}. By
$\SO'(2n)$ we denote the image of a half-spin representation of~$\Spin(2n)$.
\begin{table}[h]
\begin{tabular}{|l|l|c|c|}
\hline \str
 Type & $G/K$ & Rank & Dimension \\
\hline\hline
 A\,I & $\SU(n)/\SO(n)$ & $n-1$ & $\frac{1}{2}(n-1)(n+2)$ \\
\hline
 A\,II & $\SU(2n)/\Sp(n)$ & $n-1$ & $(n-1)(2n+1)$ \\
\hline
 A\,III & $\G_p(\C^{p+q}) = \SU(p+q)/\SUxU{p}{q}$ & $\min(p,\,q)$ & $2pq$ \\
\hline
 BD\,I & $\G_p(\R^{p+q}) = \SO(p+q)/\SO(p)\times \SO(q)$ & $\min(p,\,q)$ & $pq$ \\
\hline
 C\,I & $\Sp(n)/\U(n)$ & $n$ & $n(n+1)$ \\
\hline
 C\,II & $\G_p(\H^{p+q}) = \Sp(p+q)/\Sp(p) \times \Sp(q)$ & $\min(p,\,q)$ & $4pq$ \\
\hline
 D\,III & $\SO(2n)/\U(n)$ & $\lfloor\frac{n}{2}\rfloor$ & $n(n-1)$ \\
\hline
 E\,I & $\LE_6/ \left( \Sp(4)/\left\{\pm 1\right\} \right)$ & $6$ & $42$ \\
\hline
 E\,II & $\LE_6/\SU(6){\cdot}\Sp(1)$ & $4$ & $40$ \\
\hline
 E\,III & $\LE_6/\Spin(10){\cdot}\U(1)$ & $2$ & $32$ \\
\hline
 E\,IV & $\LE_6/\LF_4$ & $2$ & $26$ \\
\hline
 E\,V & $\LE_7/\left(\SU(8)/\left\{\pm 1\right\}\right)$ & $7$ & $70$ \\
\hline
 E\,VI & $\LE_7/\SO'(12){\cdot}\Sp(1)$ & $4$ & $64$ \\
\hline
 E\,VII & $\LE_7/\LE_6{\cdot}\U(1)$ & $3$ & $54$ \\
\hline
 E\,VIII & $\LE_8/\SO'(16)$ & $8$ & $128$ \\
\hline
 E\,IX & $\LE_8/\LE_7{\cdot}\Sp(1)$ & $4$ & $112$ \\
\hline
 F\,I & $\LF_4/\Sp(3){\cdot}\Sp(1)$ & $4$ & $28$ \\
\hline
 F\,II & $\LF_4/\Spin(9)$ & $1$ & $16$ \\
\hline
 G & $\LG_2/\SO(4)$ & $2$ & $8$ \\
\hline
\end{tabular}
\bl\caption{Symmetric spaces of type~I}\label{TSymmSpaces}
\end{table}


The global isometry classes of symmetric spaces are given by the following theorem, which will be needed
for the proof of the Splitting Theorem~\ref{SplittingTheorem}.

\begin{theorem}\label{Coverings}
Let $M$ be a simply connected Riemannian symmetric space with decomposition $M = M_0 \times M_1 \times
\ldots \times M_t$ into Euclidean and irreducible parts. Define $G = V \times I(M_1)_0 \times \ldots
\times I(M_t)_0$ where $V$ is the vector group of pure translations of the Euclidean space~$M_0$. Define
$\Delta = V \times \Delta_1 \times \ldots \times \Delta_t$ where $\Delta_i$ ($i>0$) is the centralizer
of~$I(M_i)_0$ in~$I(M_i)$. Then $G$ is the group generated by all transvections of $M$ and $\Delta$ is
the centralizer of~$G$ in~$I(M)$. In particular, the symmetric spaces covered by $M$ are just the
manifolds $M / \Gamma$ where $\Gamma$ is a discrete subgroup of $\Delta$.
\par
The group $\Delta_i$ is trivial if $M_i$ is noncompact, and is finite if $M_i$ is compact. In
particular, the discrete subgroups of $\Delta$ are just the subgroups $\Gamma\subseteq\Delta$ with
discrete projection on the vector group~$V$.
\end{theorem}


\begin{proof}
See \cite{wolf}, Ch.~8, Sec.~3.
\end{proof}


Totally geodesic submanifolds of symmetric spaces correspond to Lie triple systems.

\begin{proposition}\label{LieTripleSystems}

Let $M$ be a Riemannian globally symmetric space and let $p_0 \in M$. Let $G = I(M)_0$ and let $K =
G_{p_0}$. Let $\g = \k \oplus \pp$, where $\k$ is the Lie algebra of $K$ and where we identify $\pp =
\T_{p_0}M$ as usual. Let $\sigma_*\colon\g\rightarrow\g$ be the automorphism  of $\g$ which acts on $\k$
as $\id_{\k}$ and on $\pp$ as $-\id_{\pp}$. The totally geodesic submanifolds of $M$ containing $p_0$
are in one-to-one correspondence with the Lie triple systems $\ss \subseteq \pp$; i.e.\ if $\ss
\subseteq \pp$ is a Lie triple system, then $\Exp(\ss) \subseteq M$ is a totally geodesic submanifold
and, conversely, if $S \subseteq M$ is a totally geodesic submanifold, then $\T_{p_0} \Sigma \subseteq
\pp$ is a Lie triple system.
\par
Moreover, for any Lie triple system $\ss \subset \pp$, define $\g' = \ss+[\ss,\ss]$ and $\k' =
[\ss,\ss]$; then $\g'$ is the Lie subalgebra of $\g$ generated by $\ss$, $\g'$ is invariant under
$\sigma_*$ and $\k' = \g'\cap\k$. Let $G'$ and $K'$ be the connected Lie subgroups of $G$ with Lie
algebras $\g'$ and $\k'$, respectively. Then $(G',K')$ is a symmetric pair and $G'$ acts transitively on
$\Exp(\ss)$.
\end{proposition}


\begin{proof}
See \cite{helgason}, Ch.~IV, \S~7.
\end{proof}


\begin{lemma}\label{ProdSphDim}
Let $(G,K)$ be a symmetric pair such that $M = G / K$ is a Riemannian symmetric space of compact type
and let $\Sigma \subseteq M$ be a totally geodesic submanifold whose universal cover is a product of
spheres.  Then $\dim(\Sigma) \le \rk(G) + \rk(K)$.
\end{lemma}


\begin{proof}
Let$(G_{\Sigma},K_{\Sigma})$ be the symmetric pair corresponding to $\Sigma$, we have $\g_{\Sigma} =
\g_{\Sigma}^1 \oplus \ldots \oplus \g_{\Sigma}^m$, where $\g_{\Sigma}^i \cong \so( n_i + 1 )$. Let
$\k_{\Sigma}^i \subseteq\g_{\Sigma}^i$ such that $\k_{\Sigma}^i \cong \so( n_i )$. Let $\g = \k \oplus
\pp$ as usual. By Proposition~\ref{LieTripleSystems}, we may assume $\k_{\Sigma}^i \subseteq \k$. Now
choose maximal abelian subalgebras $\aa_i \subseteq \g_{\Sigma}^i$ as follows. If $n_i$ is even, then
$\rk (\k_{\Sigma}^i) = \rk (\g_{\Sigma}^i)$ and we may choose $\aa^i = \aa_i^{\k} \subseteq
\k_{\Sigma}^i$. If $n_i$ is odd, then we may choose $\aa^i \subseteq \g_{\Sigma}^i$ such that $\aa_i =
\aa_i^{\pp} \oplus \aa_i^{\k}$ where $\aa_i^{\pp} \subseteq \pp$ is one-dimensional and $\aa_i^{\k}
\subseteq \k_{\Sigma}^i$. Let $\aa^{\k} = \bigoplus_{i=1}^m \aa_i^{\k}$ and $\aa^{\pp} = \bigoplus_{n_i
\equiv 1(2)} \aa_i^{\pp}$. Then we have $\dim(\Sigma) = 2 \cdot \dim(\aa^{\k}) + \dim(\aa^{\pp})$. Since
$\aa^{\k} \oplus \aa^{\pp}$ is an abelian subalgebra of~$\g$ and $\aa^{\k}$ is an abelian subalgebra
of~$\k$, it follows that $\dim(\Sigma) \le \rk (G) + \rk (K)$.
\end{proof}


The estimate on the dimension given by the Lemma above is not optimal in all cases.  See \cite{cn} for
classifications of totally geodesic submanifolds in symmetric spaces.


The following Theorem, which characterizes maximal totally geodesic submanifolds of reducible symmetric
spaces, is an analogue of Theorem~15.1 in~\cite{dynkin1}, which characterizes maximal subalgebras of
semisimple Lie algebras; we give a proof which is similar to the proof in~\cite{dynkin1}.


\begin{theorem}\label{MaxTotGeodSubmf}
Let $S$ be a connected simply connected symmetric space with decomposition $S = S_0 \times S_1 \times
\ldots \times S_k$ such that $S_1,\ldots,S_k$ are irreducible and $S_0$ is of Euclidean type. Let $V$ be
a maximal totally geodesic submanifold of $S$, (i.e.\ if there is a totally geodesic submanifold $W$
such that $V \subseteq W \subseteq S$ then either $V = W$ or $W = S$). Let $p = (p_0,\ldots,p_k)\in V$.
Then either there is an index $i\in\{0,\ldots,k\}$ and a totally geodesic submanifold $\tilde{V}\subset
S_i$ such that
  $$V = S_0\times \ldots \times S_{i-1} \times \tilde{V}
        \times S_{i-1} \times \ldots \times S_k,$$
or there are two factors $S_i$ and $S_j$ ($i \neq j$) and a map $\phi\colon S_i \rightarrow S_j$  which
is an isometry up to scaling such that
  $$V = \prod_{\ell = 1 \atop \ell \neq i,j}^k S_{\ell}
        \times \{(x,\phi(x))|x \in S_i\}.$$
\end{theorem}


\begin{proof}
Let $G = I(S)$ and let $K = I(S)_p$ such that $\g = \k\oplus\pp$ is a Cartan decomposition associated
with the symmetric space $S = G/K$. Let $G_i = I(S_i)$ and let $K_i = I(S_i)_{p_i}$ such that $\g_i =
\k_i\oplus\pp_i$ are Cartan decompositions corresponding to the irreducible factors $S_i = G_i/K_i$.
Since $V$ is a totally geodesic submanifold, we have that $\nu = \T_p V\subseteq\pp$ is a Lie triple
system by Proposition~\ref{LieTripleSystems}. Obviously, the projection $\pi_i(\nu)$ onto each of the
summands $\pp_i$ is again a Lie triple system. Now there are two cases:
\par
Either there is an index $i \in \{0,\ldots,k\}$ such that $\pr_i(V)\neq S_i$, where $\pr_i \colon S
\rightarrow S_i$ denotes the canonical projection onto $S_i$. Then $\pr_i(V)$ is a totally geodesic
submanifold in $S_i$ and there is a maximal totally geodesic submanifold $\tilde{V} \subset \pp_i$
containing $\pr_i(V)$. Thus, $S_0\times \ldots \times S_{i-1} \times \tilde{V}
 \times S_{i-1} \times \ldots \times S_k$ is a totally geodesic submanifold
of $S$ containing $V$ and which is, by maximality, equal to $V$.
\par
Or $\pi_i(\nu) = \pp_i$ for all $i = 0,\ldots,k$, where the Lie algebra epimorphisms
$\pi_i\colon\g\rightarrow \g_i$ are given by the canonical projections. In this case, it follows that
there are at least two indices $i,j\in\{0,\ldots,k\}$ such that $\pp_i$ and $\pp_j$ are both not
contained in $\nu$. Define $\nu^* = \nu\cap(\pp_i\oplus\pp_j) \neq \nu$. This is a Lie triple system in
$\pp$, since it is the intersection of two Lie triple systems. Hence $\bigoplus_{\ell = 0 \atop \ell
\neq i,j}^{k}\pp_{\ell}\oplus\nu^*$ is a Lie triple system in $\pp$ which contains $\nu$ and is
different from $\pp$, thus, by maximality, is the tangent space $\T_pV$.
\par
It remains now to study the Lie triple system $\nu^*\subset\pp_i\oplus\pp_j$. By
Proposition~\ref{LieTripleSystems}, it follows that the Lie algebra $\g' = \nu^*\oplus[\nu^*,\nu^*]$
generated by $\nu^*$ is the Lie algebra of a group $G' \subset G_i \times G_j$ acting transitively on
the totally geodesic submanifold $V^*$ of $S_i \times S_j$ which is the exponential image of $\nu^*
\subset \T_{(p_i,p_j)}(S_i \times S_j)$.
\par
We show that $\g'\cap\g_i$ is an ideal in $\g_i$: Let $x \in \g_i$, $y \in \g'\cap\g_i$, there is a $z
\in \g'$ such that $\pi_i(z) = x$ and it follows that $[x,y] = [z,y]\in\g'\cap\g_i$. By the same
argument, $\g'\cap\g_j$ is an ideal in $\g_j$.
\par
Let us assume for the moment that $i,j\neq 0$. Since $\pi_i(\nu) = \pp_i$ and $\pi_j(\nu) = \pp_j$, we
have that $\pi_i(\g') = \g_i$ and $\pi_j(\g') = \g_j$ since $S_i$ and $S_j$ are irreducible symmetric
spaces. Since they are the Lie algebras of isometry groups of irreducible symmetric spaces, the $\g_i$
are either simple or the direct sum of two isomorphic simple ideals $\h_i\oplus\h_i$ (in case $S_i$ is
of type~II). Therefore the ideal $\g'\cap\g_i$ is either zero, equal to $\g_i$ or equal to $\h_i$. The
last case is impossible, since $\g'$ has to be invariant under the action of the Cartan involution of
$\h_i\oplus\h_i$, which is given by $(x,y)\mapsto(y,x)$; the case $\g'\cap\g_i = \g_i$ is also
impossible, since $\pp_i$ is not contained in $\nu^*$. Now we can show that $\g_i$ and $\g_j$ are
isomorphic: Let $x\in\g_i$ then there is a an element $y\in\g_j$ such that $(x,y)\in\g'$; but this
element is uniquely defined since otherwise, $\g_j$ would have a non-trivial intersection with $\g'$.
The map $\g_i\rightarrow\g_j$ we defined in this way is easily seen to be a Lie algebra isomorphism and
the subalgebra $\g'$ is given by the diagonal embedding of
$\g'\rightarrow\g'\oplus\g'\cong\g_i\oplus\g_j$. It remains to be shown that the spaces $S_i$ and $S_j$
are isometric up to scaling: This follows from the requirement that the Cartan involution corresponding
to $\g_i\oplus\g_j = (\k_i\oplus\k_j)\oplus(\pp_i\oplus\pp_j)$ must leave the diagonally embedded
subalgebra $\g'$ invariant and is thus of the form $(x,y)\mapsto(\sigma(x),\sigma(y))$, where $\sigma$
is an involution of $\g_i\cong\g_j$.
\par
Finally, assume $i = 0$, $j\in\{1,\ldots,k\}$. This case can be included in the above proof if we
further split $\pp_0$ into the direct sum of $\g'\cap\pp_0$ plus a complementary subspace. Then one is
again in the situation that $\g$ can be written as a direct sum of ideals all of which have either
trivial intersection with the Lie algebra generated by $\nu$ or are contained in this Lie algebra. Then
the same type of argument leads to the contradiction that an abelian Lie algebra is isomorphic to one of
$\g_1,\ldots,\g_k$.
\end{proof}


\begin{corollary}\label{TotGeodHyperSurface}
Let $S$ be a connected simply connected symmetric space with decomposition $S = S_0 \times S_1 \times
\ldots \times S_k$ such that $S_1, \ldots, S_k$ are irreducible and $S_0$ is of Euclidean type. Let $H$
be a totally geodesic hypersurface of $S$. Let $p = (p_0, \ldots, p_k) \in H$. Then there is an index $i
\in \{0, \ldots, k\}$ and a totally geodesic hypersurface $\tilde{H}\subset S_i$ such that
  $$H = S_0\times \ldots \times S_{i-1} \times \tilde{H}
        \times S_{i-1} \times \ldots \times S_k.$$
\end{corollary}


\begin{proof}
Obviously, a totally geodesic hypersurface is a maximal totally geodesic submanifold, so we may apply
Theorem~\ref{MaxTotGeodSubmf}. Irreducible non-flat symmetric spaces are at least of dimension two, thus
the second possibility in the assertion of Theorem~\ref{MaxTotGeodSubmf} does not occur here, since this
would lead to submanifolds of codimension at least two.
\end{proof}


The following facts on the maximal connected subgroups of the classical groups can be proven by standard
arguments from the representation theory of compact Lie groups, see e.g.\ ~\cite{dynkin2}. It should be
remarked that some of subgroups of~$\SO(n)$, $\SU(n)$ or $\Sp(n)$ given by irreducible representations
of simple groups of corresponding (real, complex or quaternionic) type are not maximal connected,
see~\cite{dynkin2} for complete lists of inclusions.


\begin{proposition}\label{MaxSubgrSOn}
Let $K$ be a connected proper subgroup of $\SO(n)$. Then there is an automorphism $\a$ of $\SO(n)$ such
that $\a(K)$ is contained in one of the following subgroups of $\SO(n)$
\begin{enumerate}

\item $\SO(k)\times\SO(n-k),\quad$ $1 \le k \le n-1$ .

\item $\SO(p)\otimes\SO(q),\quad$  $pq = n,\,3\le p \le q$.

\item $\U(k),\quad 2k = n$.

\item $\Sp(p)\cdot\Sp(q),\quad$ $4pq = n\neq4$ .

\end{enumerate}
or $K$ is a simple irreducible subgroup $K = \p(H) \subset \SO(n)$, where $H$ is a simple compact Lie
group and $\p$ is an irreducible representation of $H$ of real type such that $\deg\p = n$.
\end{proposition}


\begin{proposition}\label{MaxSubgrSUn}
Let $K$ be a connected proper subgroup of $\SU(n)$. Then there is an automorphism $\a$ of $\SU(n)$ such
that $\a(K)$ is contained in one of the following subgroups of $\SU(n)$
\begin{enumerate}

\item $\SO(n)$

\item $\Sp(m),\quad 2m = n$

\item $\SUxU{k}{n-k},\quad 1 \le k \le n-1$

\item $\SU(p)\otimes\SU(q),\quad pq = n,\, p \ge 3,\,q \ge 2$

\end{enumerate}
or $K$ is a simple irreducible subgroup $K = \p(H) \subset \SU(n)$, where $H$ is a simple compact Lie
group and $\p$ is an irreducible representation of $H$ of complex type such that $\deg\p = n$.
\end{proposition}


\begin{proposition}\label{MaxSubgrSpn}
Let $K$ be a connected proper subgroup of $\Sp(n)$. Then there is an automorphism $\a$ of $\Sp(n)$ such
that $\a(K)$ is contained in one of the following subgroups of $\Sp(n)$
\begin{enumerate}

\item $\U(n),$

\item $\Sp(k)\times\Sp(n-k),\quad 1 \le k \le n-1$

\item $\SO(p)\otimes\Sp(q),\quad pq = n,\, p\ge 3,\, q\ge 1$

\end{enumerate}
or $K$ is a simple irreducible subgroup $K = \p(H) \subset \Sp(n)$, where $H$ is a simple compact Lie
group and $\p$ is an irreducible representation of $H$ of quaternionic type such that $\deg\p = 2n$.
\end{proposition}


\section{Criteria for polarity}\label{CriteriaPolarity}

The following is a generalization of the criterion for hyperpolarity given in~\cite{hptt}. Note that we
do not require the sections to be embedded submanifolds here. Hyperpolar actions are characterized by
the property that the Lie triple system $\nu$ in Proposition~\ref{PolCrit} is abelian.

\begin{proposition}\label{PolCrit}
Let $G$ be a connected compact Lie group, $K \subset G$ a symmetric subgroup and let $\g = \k + \pp$ be
the Cartan decomposition. Let $H \subseteq G$ be a closed subgroup. Let $k$ be the cohomogeneity of the
$H$-action on $G$. Then the following are equivalent.
\begin{enumerate}

\item The $H$-action on $G / K$ is polar w.r.t\ some Riemannian
metric induced by an $\Ad(G)$-invariant scalar product on~$\g$.

\item For any $g \in G$ such that $gK$ lies in a principal orbit
of the $H$-action on~$G / K$ the subspace $\nu = g^{-1} N_{g K}(H \cdot gK) \subseteq \pp$ is a
$k$-dimensional Lie triple system such that the Lie algebra~$\ss = \nu \oplus [\nu, \nu]$ generated by
$\nu$ is orthogonal to $\Ad(g^{-1})\h$.

\item The normal space $\N_{\eK}(H \cdot \eK) \subseteq \pp$
contains a $k$-dimensional Lie triple system $\nu$ such that the Lie algebra~$\ss = \nu \oplus [\nu,
\nu]$ generated by $\nu$ is orthogonal to $\h$.

\end{enumerate}
\end{proposition}


\begin{proof}
Let $g \in G$ be such that $gK$ lies in a principal orbit of the $H$-action on~$G / K$. Then the action
of $g^{-1}Hg$ on~$G / K$ has a principal orbit containing~$\eK$ and the equivalence of (i) and (ii)
follows from \cite{gorodski}, Proposition, p.~193.
\par
Assume now condition (i) holds. Let $\Sigma$ be a section of the polar $H$-action on~$G / K$ such that
$\eK \in \Sigma$. Let $\nu = T_{\eK} \Sigma \subseteq \pp$, let $\ss = \nu \oplus [\nu, \nu]$, and let
$S$ be the connected subgroup of~$G$ corresponding to~$\ss$. Since $S$ acts transitively on~$\Sigma$,
there is an element $s \in S$ such that the point $sK$ lies in a principal orbit of the $H$-action on~$G
/ K$. Now it follows from (ii) that $\Ad(s^{-1})\h$ is orthogonal to the Lie algebra generated by
$s^{-1} N_{s K}(H \cdot sK)$, which coincides with~$\ss$. Since $\Ad(s^{-1})$ leaves~$\ss^{\perp}$
invariant, we have that $\h$ is orthogonal to~$\ss$ and (iii) follows.
\par
We will now show that if (iii) holds, then $\Sigma = \exp(\nu) \subseteq G / K$ meets the orbits
orthogonally. Let $s K \in \Sigma$, where $s$ is an arbitrary element of the Lie group $S$ corresponding
to the Lie algebra generated by $\nu$. The tangent space of the $H$-orbit through $sK$ is orthogonal to
$\T_{sK} \Sigma$ if and only if $s^{-1} \h {s} \perp \nu$. But since the adjoint representation of $G$
restricted to $S$ leaves the orthogonal complement of $\ss$ invariant, $s^{-1}\h{s}$ is perpendicular to
$\ss$. Thus the $H$-action on $G / K$ is polar.
\end{proof}

As an immediate consequence of this criterion, the problem of classifying polar actions on $G / K$ is
reduced to a problem on the Lie algebra level. We conclude this section with the simple observation that
a polar action restricted to an invariant totally geodesic submanifold is polar.

\begin{lemma}\label{PolarOnTotGeodSubmf}
Let $G$ be a compact Lie group acting polarly on a connected Riemannian manifold~$N$. Let $M \subseteq
N$ be a connected totally geodesic submanifold which is invariant under the $G$-action. Then the
$G$-action on $M$ is polar.
\end{lemma}

\begin{proof}
Let $\Sigma \subseteq N$ be a section of the $G$-action on~$N$. Let $\Sigma_0$ be a connected component
of $\Sigma \cap M$. Then the totally geodesic submanifold $\Sigma_0 \subseteq M$ obviously meets the
$G$-orbits in~$M$ orthogonally at every intersection point. Furthermore, since $M$ is connected, any two
orbits of the $G$-action on $M$ can be joined by a shortest geodesic which meets the principal
$G$-orbits orthogonally and is hence contained in $\Sigma$ after conjugation with a group element. This
geodesic is now also contained in $M$, since $M$ is totally geodesic. This proves that $\Sigma_0$ meets
all $G$-orbits in $M$.
\end{proof}


\section{Sections and Weyl group actions}\label{SectionsWeylAct}

Let us first recall some known properties of the Weyl group.

\begin{lemma}[Thorbergsson, Podest{\`a}]\label{TotGeodHyperWeyl}
Let $M$ be a simply connected symmetric space on which a compact, connected Lie group $G$ acts polarly
and nontrivially. Let $\Sigma$ be a section of the polar action and let $p \in \Sigma$ be such that the
orbit through $p$ is singular. Then there is a totally geodesic hypersurface $H$ in~$\Sigma$ passing
through~$p$ and consisting of singular points; moreover there exists a non-trivial element $g \in
W_{\Sigma}$ which fixes $H$ pointwisely.
\par
The set of singular points in $\Sigma$ is a union of finitely many totally geodesic hypersurfaces
$\{H_i\}_{i \in I}$ in~$\Sigma$; the Weyl group $W_{\Sigma}$ is generated by reflections in the
hypersurfaces $\{H_i\}_{i \in I}$. Let $\tilde\Sigma$ be the universal covering of $\Sigma$ and let
$\{P_j\}_{j \in J}$ be the collection of all lifts of all the totally geodesic hypersurfaces $\{H_i\}_{i
\in I}$ in~$\Sigma$. Let $\hat W_{\Sigma}$ be the subgroup of the isometry group of~$\tilde \Sigma$
which is generated by the reflections in the hypersurfaces $\{P_j\}_{j \in J}$. Then $\hat W_{\Sigma}$
is a Coxeter group and $W_{\Sigma}$ is a quotient group of~$\hat W_{\Sigma}$.
\end{lemma}


\begin{proof}
See \cite{pth1}, Lemma~1A.4 or \cite{ewert}, Section~2.3 for a more general statement.
\end{proof}


The following splitting theorem is a generalization of Lemma~1A.2 in~\cite{pth1}, where $\tilde
\Sigma_1$ is a point and the hypothesis is equivalent to a trivial Weyl group action. We consider the
weaker hypothesis that the section of a polar action is locally a product such that the Weyl group acts
trivially on one factor.


\begin{splittingtheorem}\label{SplittingTheorem}

Let $N$ be a compact connected Riemannian symmetric space on which a connected compact Lie group $G$
acts polarly. Assume the universal covering $\tilde\Sigma$ of a section $\Sigma$ decomposes as a
Riemannian product $\tilde\Sigma
 = \tilde\Sigma_1 \times \tilde\Sigma_2$ and the action of $\hat
W_{\Sigma}$ on $\tilde\Sigma$ descends to an action on $\tilde\Sigma_1$ such that
$$
w \cdot (p,q) = (w \cdot p, q )\quad \mbox{for all $w \in \hat W_{\Sigma}$, $p\in\tilde\Sigma_1$,
$q\in\tilde\Sigma_2$.}
$$
Then the universal cover of~$N$ is a Riemannian product isometric to~$\tilde M \times \tilde\Sigma_2$,
where $M = G \cdot \Sigma_{1}$ and $\Sigma_{1}$ is the image of $\tilde \Sigma_{1}$ under the covering
map $\tilde \Sigma \rightarrow \Sigma$.
\end{splittingtheorem}


\begin{proof}
Let $\Sigma$ be a section and $p\in\Sigma$ be an arbitrary point of this section. For $i = 1,\,2$, let
$\Sigma_i = \Sigma_{i}(p)$ be the totally geodesic submanifolds of $\Sigma$ corresponding to
$\tilde\Sigma_i$ (uniquely determined by Theorem~\ref{Coverings}) such that $p\in\Sigma_{i}$.

First we show that the isotropy group $G_p$ acts trivially on $\T_p\Sigma_2$. Consider the slice
representation of $G_{p}$ on $V = N_{p}(G \cdot p)$ which is polar by Proposition~\ref{PolarSlice}, with
section $\T_{p}\Sigma = \T_p\Sigma_1 \oplus \T_p\Sigma_2$. Now consider the Weyl group $W'$ of this
polar linear representation; it coincides with $(W_{\Sigma})_{p}$. Its representation space decomposes
into a sum of irreducible modules and one trivial module and the section $\T_{p}\Sigma$ decomposes
accordingly, see~\cite{dadok}. It follows from the hypothesis that $W'$ acts trivially on the linear
subspace $\T_p\Sigma_2$. Since irreducible polar representations have irreducible Weyl groups it follows
that $G_{p}$ acts trivially on $\T_p\Sigma_2$.

We will now show that the set $M(p) = G \cdot \Sigma_1$ is an embedded submanifold of $N$. By the Slice
Theorem~\ref{SliceTheorem} there is an equivariant diffeomorphism $\Psi$ of a $G$-invariant open
neighborhood around the zero section in the normal bundle $G \times_{G_p} V \rightarrow G/G_p$ onto a
$G$-invariant open neighborhood around the orbit $G\cdot p$ such that the zero section in $G
\times_{G_p} V$ is mapped to the orbit $G\cdot p$. The diffeomorphism $\Psi$ is given by the end point
map which maps any normal vector $v_q \in V = \N_q(G \cdot p)$ to its image under the exponential map
$\exp_q(v_q)$.

Since $\Sigma_2$ is a totally geodesic submanifold of $N$, we have $\exp_p(\T_p\Sigma_2) = \Sigma_2$.
The subspace $\T_p \Sigma_1 \subseteq V$ is fixed by the Weyl group $W'$ of the slice representation and
hence fixed by $G_{p}$. Hence the orthogonal complement $S$ of $\T_{p}\Sigma_2$ in $V$ is a linear
subspace invariant under the polar representation of $G_{p}$ on $V$. Therefore, $S$ defines a smooth
subbundle of the normal bundle $G \times_{G_p} V$.

But since $\T_{p}\Sigma_1$ is the section of the $G_{p}$-representation on $S$ we have $S = G_p \cdot
\T_p\Sigma_1$. From the fact that $\Psi$ is an equivariant diffeomorphism it follows now that the
elements of the subbundle defined by $S$ are mapped into the set $G \cdot \Sigma_1$. This shows that, in
a neighborhood of $p$, the subset $M(p) = G \cdot \Sigma_1 \subseteq N$ is a smooth submanifold of
codimension $\dim(\Sigma_2)$. Thus we see that the symmetric space $N$ is foliated by the totally
geodesic submanifolds $\{g \cdot \Sigma_{2}(p)\}_{g \in G,\,p \in \Sigma}$ with integrable normal bundle
whose integral manifolds are given by $\{M(p)\}_{p \in \Sigma}$.
\par
It follows from Theorem~A of~\cite{bh} that the universal cover~$\tilde N$ of~$N$ is topologically a
product diffeomorphic to $\tilde M \times \tilde \Sigma_2$ such that the projection of $\tilde N$ on the
factor $\tilde \Sigma_2$ is a Riemannian submersion. We have just shown that the horizontal distribution
of this Riemannian submersion is integrable. Since the sectional curvature of $N$ is nonnegative, it
follows from Theorem~1.3 of~\cite{walschap} that the fibers of this Riemannian submersion are totally
geodesic. We conclude that $\tilde N$ is a Riemannian product isometric to $\tilde M \times \tilde
\Sigma_2$.
\end{proof}


\begin{corollary}\label{SingularOrbit}
Let $N$ be an irreducible Riemannian symmetric space of compact type on which a compact Lie group $G$
acts polarly and nontrivially. Then the $G$-action on~$N$ has a singular orbit.
\end{corollary}

\begin{proof}
Assume there is no singular orbit. Then by Lemma~\ref{TotGeodHyperWeyl} the Weyl group $W_{\Sigma}$ acts
trivially on $\Sigma$. Hence it follows from the Splitting Theorem~\ref{SplittingTheorem} that $\tilde
N$ is a Riemannian product $\tilde M \times \tilde \Sigma$, where $M$ is a $G$-orbit. But this is a
contradiction to the irreducibility of~$N$.
\end{proof}


The following theorem is a generalization of Proposition~1B.1 of~\cite{pth1}, where it was proved that
the section of a polar action on a compact rank one symmetric space has constant curvature.


\begin{theorem}\label{ProductOfSpheres}
Let $N$ be an irreducible compact simply connected symmetric space on which a compact Lie group $G$ acts
polarly and non-trivially with section $\Sigma$. Then $\Sigma$ is covered by a Riemannian product of
spaces which have constant curvature.
\end{theorem}


\begin{proof}
Let $\tilde\Sigma = \tilde\Sigma_1\times\tilde\Sigma_2$ be a decomposition of the universal covering
$\tilde\Sigma$ of $\Sigma$ such that $\Sigma_1$ is a Riemannian product of spaces of constant curvature
and $\Sigma_2$ is either a point or a Riemannian product of irreducible symmetric spaces of non-constant
curvature. The section $\Sigma$ contains a union of finitely many totally geodesic hypersurfaces
$\{H_i\}_{i \in I}$ such that the Weyl group $W_{\Sigma}$ is generated by the reflections in the
hypersurfaces $\{H_i\}_{i \in I}$. In view of Corollary~\ref{TotGeodHyperSurface} and since it is well
known that the only irreducible symmetric spaces containing totally geodesic hypersurfaces are those of
constant curvature, it is clear that the hypothesis of the Splitting Theorem~\ref{SplittingTheorem} is
fulfilled and we conclude that $\Sigma_2$ is a point.
\end{proof}


It follows from Theorem~\ref{ProductOfSpheres} and Lemma~\ref{ProdSphDim} that the cohomogeneity of a
polar action on an irreducible symmetric space~$G /K$ is less or equal~$\rk(G) + \rk(K)$. For hermitian
symmetric spaces $G / K$ the upper bound on the cohomogeneity can be further improved, see
Proposition~\ref{CoisoDimBound} below. These dimension bounds are essential for our classification,
since they reduce the classification problem to a finite number of cases.


\begin{proposition}\label{CoisoDimBound}
Let $H$ be a compact Lie group acting polarly on a compact K{\"a}hler manifold~$M$. Then the cohomogeneity
of the $H$-action on~$M$ is less or equal~$\rk(H)$.
\end{proposition}

\begin{proof}
By the Equivalence Theorem~\cite{hw}, see also Theorem~1.4 in~\cite{pth2}, the cohomogeneity of the
$H$-action is equal to the difference between the rank of $H$ and the rank of a regular isotropy
subgroup of $H$.
\end{proof}


We have the following lower bounds on the dimension of groups acting polarly on the classical symmetric
spaces.

\begin{proposition}\label{ClassDimBounds}
Let $H$ be a connected compact Lie group acting polarly and non-trivially on a symmetric space~$M$.
Assume $3 \le k \le n - 3$, $2 \le \ell \le n - 2$ and let $d = \dim (H)$.
\begin{enumerate}

\item If $M = \G_{k}(\R^n)$, then $d \ge 2n - 9$.

\item If $M = \G_{\ell}(\C^n)$, then $d \ge 3n - 7$.

\item If $M = \G_{\ell}(\H^n)$, then $d \ge 6n - 16$.

\item If $M = \SO(n) / \U(\frac{n}{2})$, then $d \ge \frac{n^2}{4}
- n$.

\item If $M = \SU(n) / \SO(n)$, then $d \ge \frac{n^2}{2} - n$.

\item If $M = \SU(n) / \Sp(\frac{n}{2})$, then $d \ge
\frac{n^2}{2} - 2n$.

\item If $M = \Sp(n) / \U(n)$, then $d \ge n^2$.

\end{enumerate}
\end{proposition}


\begin{proof}
Follows from Proposition~\ref{CoisoDimBound} in case of the spaces $\G_{\ell}(\C^n)$, $\SO(n) /
\U(\frac{n}{2})$, and $\Sp(n) / \U(n)$, which are Hermitian symmetric, and from Theorem
~\ref{ProductOfSpheres} and Lemma~\ref{ProdSphDim} otherwise.
\end{proof}


\section{Polar subactions}\label{PolSub}


In this section, we will introduce our main tool for classifying polar actions through studying slice
representations. The basic observation is the following maximality property of linear polar actions,
see~\cite{kp}, Theorem~6.


\begin{theorem}\label{MaxPolRep}
Let $G \subset \SO(n)$ be a closed connected subgroup which acts irreducibly on $\R^n$ and
non-transitively on the sphere $\eS^{n-1}\subset\R^n$. Let $H \subseteq G$ be a closed connected
subgroup $\neq\{\e\}$ that acts polarly on $\R^n$. Then the $H$-action and the $G$-action on $\R^n$ are
orbit equivalent.
\end{theorem}


The proof of the above theorem relies on~\cite{simons}. As an immediate consequence of
Theorem~\ref{MaxPolRep} we have the following, cf. \cite{brueck}, Theorem~2.2.


\begin{corollary}\label{PolFixOnHomSp}
Let $X$ be a strongly isotropy irreducible Riemannian homogeneous space. Assume a connected compact Lie
group $H$ acts polarly on~$X$ such that the $H$-action has a one-dimensional orbit $H \cdot p$ or a
fixed point $p \in X$. Then the space $X$ is locally symmetric. Furthermore, $X$ is a rank-one symmetric
space or the action of $H$ is orbit equivalent to the action of the connected component of the isotropy
group of $X$ at~$p$.
\end{corollary}


\begin{proof}
Assume first that $p \in X$ is a fixed point of the $H$-action on $X$. Let $K$ be the connected
component of the isotropy group of $p$. The isotropy representation of $K$ on $\T_{p} X$ restricted to
$H$ is polar by Proposition~\ref{PolarSlice}. If the action of $K$ on the unit sphere in $\T_{p} X$ is
transitive, then the space $X$ is rank-one symmetric. If $K$ does not act transitively on the sphere,
then the linear $H$-action on $\T_p X$ is orbit equivalent to the $K$-action by Theorem~\ref{MaxPolRep},
in particular they have the same cohomogeneity; hence the $K$-action on~$\T_p X$ is polar. It now
follows that the principal orbits of the $H$-action agree with those of the $K$-action on~$X$ and the
orbit equivalence of the two actions follows from Proposition~\ref{PolCrit}, since the principal orbits
of a hyperpolar action determine all other orbits. In case $X$ is compact, the symmetry follows from
\cite{kp}, since then one may assume that $X$ is a homogeneous space of a simple compact Lie group, see
\cite{wolfIrr}, Chapter~I.1. Non-compact strictly isotropy irreducible Riemannian homogeneous spaces are
symmetric by \cite{wolfIrr}.

Now assume $p \in X$ is such that $\dim(H \cdot p) = 1$. If $H \cdot p$ is a regular orbit, it follows
that a section~$\Sigma \subset X$ is a totally geodesic hypersurface and hence $X$ is locally isometric
to a space of constant curvature. Assume now that $H \cdot p$ is a singular orbit, hence the slice
representation of~$H_p$ on~$\N_p (H \cdot p)$ is nontrivial and polar by Proposition~\ref{PolarSlice}.
However, since $\T_p (H \cdot p)$ is one-dimensional, the isotropy representation of~$H_p$ on~$\T_p X =
\T_p (H \cdot p) \oplus \N_p (H \cdot p)$ is polar. It now follows from Theorem~\ref{MaxPolRep} that the
irreducible isotropy representation of $M$ at~$p$ is orbit equivalent to the reducible $H_p$-action on
$\T_p M$, a contradiction.
\end{proof}

In particular, we may restrict our attention to actions without fixed point in the following. The full
classification of connected Lie groups acting polarly with a fixed point on the irreducible symmetric
spaces of higher rank follows immediately from Lemma~\ref{OrbitEqPolarSubGr}. As the proof of
Corollary~6.2 shows, one obtains the same result also under the weaker hypothesis that the linear action
of~$H$ on the tangent space $\T_p X$ is polar.
\par
Let $G$ be a connected compact Lie group acting isometrically on a Riemannian manifold. We say the
action of $G$ on $M$ is {\em polarity minimal} if there is no closed connected subgroup~$H \subset G$
which acts nontrivially and polarly on $M$ and such that the $H$-action is not orbit equivalent to the
$G$-action. Note that a polarity minimal action can be polar or non-polar. We give various sufficient
conditions for an orthogonal representation to be polarity minimal in the following proposition.


\begin{proposition}\label{PolMinCriteria}
Let $\rho \colon G \to \O(V)$ be a representation of the compact connected Lie group~$G$. Then $\rho$ is
polarity minimal if one of the following holds.
\begin{enumerate}

\item[(i)] The representation $\rho$ is irreducible of
cohomogeneity~$\ge 2$.

\item[(ii)] The representation space~$V$ is the direct sum of two
equivalent $G$-modules.

\item[(iii)] The representation space $V$ contains a
$G$-in\-var\-iant sub\-module~$W$ such that the $G$-re\-pre\-sen\-ta\-tion on~$W$ is almost effective,
non-polar, and polarity minimal.

\end{enumerate}
\end{proposition}


\begin{proof}
Part~(i) is a just a reformulation of Theorem~\ref{MaxPolRep}. Assume now $V$ is the direct sum of two
equivalent $G$-modules; then the representation~$\rho$ restricted to any closed connected subgroup $H
\subseteq G$ which acts nontrivially on~$V$ will have two equivalent nontrivial submodules; it then
follows from~\cite{kollross}, Lemma~2.9 that $H$ acts non-polarly on~$V$; this proves part~(ii). To
prove part~(iii), assume there is a closed connected subgroup~$H$ of~$G$ acting polarly on~$V$; since
the $G$-action on the subspace~$W$ is non-polar and polarity minimal, it follows that $H$ acts trivially
on~$W$. But $W$ is an almost effective representation, thus $H$ acts trivially on all of~$V$.
\end{proof}


While we do not have an {\em a priori} proof that polar actions on irreducible compact symmetric spaces
of higher rank are orbit maximal, the following proposition gives various sufficient conditions under
which one can show that certain non-polar actions are polarity minimal. In fact, this is our main tool
to exclude subactions and it will be used frequently in the sequel.

\begin{lemma}\label{PolMinHered}
Let $G$ be compact Lie group and $K \subset G$ be symmetric subgroup such that $M = G / K$ is an
irreducible symmetric space and let $H \subset G$ be a closed subgroup. The action of $H$ on $M$ is
non-polar and polarity minimal if there is a non-polar polarity minimal submodule~$V \subseteq \N_p(H
\cdot p)$ of the slice representation at~$p$ such that one of the following holds.
\begin{enumerate}

\item $M$ is Hermitian symmetric and $\dim(V) > \rk(H)$.

\item $\dim(V) > s(M)$, where $s(M)$ is the maximal dimension of a
totally geodesic submanifold of~$M$ locally isometric to a product of spaces with constant curvature,
cf.\ Lemma~\ref{ProdSphDim}.

\item $V \subseteq \pp = \T_p M$ (where $\g = \k \oplus \pp$ as
usual such that $\k$ is the Lie algebra of~$K = G_p$) contains a Lie triple system corresponding to an
irreducible symmetric space of nonconstant curvature, e.g.\ an irreducible symmetric space of higher
rank.

\item The isotropy group $H \cap K$ acts almost effectively on $V$
and $\rk(H \cap K) = \rk(H)$.

\end{enumerate}
\end{lemma}

\begin{proof}
Assume a closed connected subgroup~$U \subseteq H$ acts polarly on~$M$. Consider the isotropy group at
$U_p$ of the $U$-action on~$M$. Since $U_p \subseteq H_p$, the action of~$U_p$ on the normal space $\N_p
(U \cdot p)$ leaves the subspace~$V$ invariant. By Proposition~\ref{PolarSlice}, the slice
representation of $U_p$ on $\N_p (U \cdot p)$ is polar, in particular, the $U_p$-action on $V$ is polar.
Since the action of~$H_p$ on~$V$ is polarity minimal and non-polar, it follows that the action of the
connected component~$\left ( U_p \right )_0$ on~$V$ is trivial. Hence $V$ is contained in the section of
the polar $U_p$-action on $\N_p (U \cdot p)$ and thus $V$ is tangent to a section $\Sigma$ of the
$U$-action on~$M$; in particular, $\dim(\Sigma) \ge \dim(V)$. Part (i) now follows from
Lemma~\ref{CoisoDimBound}. Parts (ii) and (iii) follow from Theorem~\ref{ProductOfSpheres}. If $\rk(H
\cap K) = \rk(H)$ and $H \cap K$ acts almost effectively on~$V$, then any closed subgroup $U \subseteq
H$ with $\dim U > 0$ will have an intersection $U \cap K \subseteq H \cap K$ of positive dimension
with~$K$; but $U \cap K$ acts on~$V$ non-polarly since $V$ is polarity minimal; this proves~(iv).
\end{proof}


\section{Hermann actions of higher cohomogeneity}
\label{SubHermannHighRk}

In the remaining part of the paper, we will carry out the classification. We begin with subactions of
Hermann actions whose cohomogeneity is~$\ge 2$.

To study actions of reducible groups on the Grassmannians we will need the following technical lemma.
Let us first introduce some notation. Let $G \subseteq \Gl(n,\R)$. Let $V$ be a linear subspace
of~$\R^n$. Then we define the {\em normalizer} of~$V$ in~$G$ as $N_G(V) = \left \{ g \in G \mid g(V) = V
\right \}$, and similarly by $Z_G(V) = \left \{ g \in G \,\mid\, g|_V = \id_V \right \}$ the {\em
centralizer} of~$V$ in~$G$. Clearly, $N_G(V)$ is a subgroup of~$G$ and since the elements of~$N_G(V)$
leave $V$ invariant, the group~$N_G(V)$ acts on~$V$. The kernel of this representation is the normal
subgroup $Z_G(V) \subseteq N_G(V)$. The group $N_G(V)$ is the isotropy subgroup $G_V$ of the $G$-action
on the Grassmannian~$\G_{\dim V}(\R^n)$ of $(\dim V)$-dimensional linear subspaces in~$\R^n$.


\begin{lemma}\label{GrassIsotropy}
Let $H \subset \SO(n)$ be a closed connected proper subgroup.
\begin{enumerate}
\item If for any $8$-dimensional subspace $V \subseteq \R^n$ the
natural action of the connected component of $N_{H}(V)/Z_{H}(V)$ on $V$ is equivalent to the
8-dimensional spin representation of $\Spin(7)$ or the standard representation of $\SO(8)$ then $n = 8$
and $H \cong \Spin(7)$.

\item If for any $7$-dimensional subspace $V \subseteq \R^n$ the
natural action of the connected component of $N_{H}(V)/Z_{H}(V)$ on $V$ is equivalent to the
7-dimensional irreducible representation of $\LG_{2}$ or the standard representation of $\SO(7)$ then
either $n = 7$ and $H \cong \LG_{2}$ or $n = 8$ and $H \cong \Spin(7)$.

\end{enumerate}
\end{lemma}

\begin{proof}
We first show that in both cases the group $H$ acts transitively on the unit sphere in $\R^n$. Let $p$,
$q \in \R^n$ be two unit vectors and let $V$ be linear subspace of~$\R^n$ containing $p$ and $q$ such
that $V$ is $8$- or $7$-dimensional, respectively. Then it follows from the hypothesis that
$N_{H}(V)/Z_{H}(V)$ acts transitively on the unit sphere in the space~$V$, thus there is an element
in~$H$ which maps $p$ to $q$. This shows that $H$ acts transitively on the unit sphere in $\R^n$ and
hence the pair $(H,n)$ is one of the following, see Table~7 in~\cite{oniscik}.
$$
\begin{tabular}{|c||c|c|c|c|c|c|c|c|}
  \hline
  $H$ & $\U(m),$ & $\SU(m),$ & $\Sp(\ell)\cdot\Sp(1),$ & $\Sp(\ell)\cdot\U(1),$
  & $\Sp(\ell),$ & $\Spin(7)$ & $\Spin(9)$ & $\LG_2$ \\
  & $m \ge 2$ & $m \ge 2$ & $\ell \ge 2$ & $\ell \ge 2$ &
  $\ell \ge 2$ & & & \\ \hline
  $n$ & $2m$ & $2m$ & $4 \ell$ & $4 \ell$ & $4 \ell$ & $8$ & $16$ & $7$
  \\ \hline
\end{tabular}
$$
It is easy to see that the first five groups do not have the property described in the hypothesis. For
the groups $\Spin(7)$ and $\LG_2$ the statement is either trivial or follows from the well-known fact
$\eS^7 = \Spin(7) / \LG_2$.
\par
It remains the case of $H = \Spin(9)$ acting on~$\R^{16}$. To prove the assertion for $\dim V = 8$, it
suffices to exhibit an isotropy subgroup of the $\Spin(9)$-action on $\G_8(\R^{16})$ not containing
$\Spin(7)$ as a Lie subgroup. First choose an $8$-dimensional subspace~$\tilde V \subset \R^{16}$ such
that the subgroup~$\Spin(8)$, acting by a representation equivalent to the sum of the two half-spin
representations on~$\R^{16}$, stabilizes~$\tilde V$. Since $\Spin(8) \subset \Spin(9)$ is maximal
connected, it coincides with the connected component of the isotropy group $H_{\tilde V}$ of the
$H$-action on the Grassmannian $\G_8(\R^{16})$. Thus the $H$-orbit through~$\tilde V$ is
$8$-dimensional. We will determine the slice representation of the $H$-action at~$\tilde V$. The group
$\left( H_{\tilde V} \right)_0 \cong \Spin(8)$ acts on the tangent space $\T_{\tilde V}\G_8(\R^{16})$ by
the tensor product of the two half-spin representations of~$\Spin(8)$. By Weyl's dimension formula, this
representation contains an irreducible summand which is $56$-dimensional and must therefore coincide
with the normal space~$\N_{\tilde V} (H \cdot \tilde V)$. This shows that $\tilde V$ lies in a singular
orbit of the~$H$-action on~$\G_8(\R^{16})$. By~\cite{hh}, the principal isotropy subgroups of this slice
representation are finite and we conclude that for generic subspaces~$V \subset \R^{16}$ the group
$N_H(V) / Z_H(V)$ is finite.
\par
Similarly, to prove the assertion for $H = \Spin(9)$ and $\dim V = 7$, choose a $7$-dimensional subspace
$\tilde V \subset \R^{16}$ such that $\Spin(7) \subset \Spin(9)$ stabilizes~$\tilde V$. Using an
analogous argument as in the case $k = 8$ we see that the $48$-dimensional slice representation at
$\tilde V$ of the $H$-action on~$\G_7(\R^{16})$ has finite principal isotropy subgroups and hence for a
generic $7$-dimensional subspace $V \subset \R^{16}$ the group $N_H(V) / Z_H(V)$ is finite.
\end{proof}


\begin{lemma}\label{ReducibleGrass}
Let $H$, $G$, $K$ be as in the following table, where $2\le k,\ell\le\frac{n}{2}$.
\begin{center}
  \begin{tabular}{|c|c|}\hline
\hl \str $H$& $G/K$
\\ \hline \hline
\hl $\SUxU{k}{n-k}$,& $\SU(n) / \SUxU{\ell}{n-\ell}$
\\ \hline
\hl $\SO(k){\times}\SO(n-k)$& $\SO(n) / \SO(\ell){\times}\SO(n-\ell)$
\\ \hline
\hl $\Sp(k){\times}\Sp(n-k)$& $\Sp(n) / \Sp(\ell){\times}\Sp(n-\ell)$
\\\hline
\end{tabular}
\end{center}
Let $U$ be a connected subgroup of $H$. Then the action of $U$ on $G / K$ is polar if either $U$ = $H$
or $U$ is conjugate to one of the following subgroups, where in each case the $U$-action on $G / K$ is
orbit equivalent to the $H$-action. In particular, the $U$-action on $G / K$ is hyperpolar.
\begin{center}
\begin{tabular}{|r|c|l|}
\hline
\hl \str $U$ & $G$ & Range \\
\hline\hline \str $\LG_{2}\times\LG_{2}$ & $\SO(14)$ & $\ell = 2$ \\
\hline \hl $\LG_{2}\times\Spin(7)$ & $\SO(15)$ & $\ell = 2$ \\
\hline \hl $\Spin(7)\times\Spin(7)$ & $\SO(16)$ & $\ell = 2,3$ \\
\hline \hl $\LG_{2}\times\SO(n-7)$ & $\SO(n)$ & $\ell = 2,\,n\ge9$
\\ \hline \hl $\Spin(7)\times\SO(n-8)$ & $\SO(n)$ &
$\ell = 2,3,\,n\ge10$ \\ \hline \hl $\SU(k) \times \SU(n-k)$ & $\SU(n)$ & $(k,\ell) \neq
\left(\frac{n}{2},\frac{n}{2}\right)$
\\\hline
\end{tabular}
\end{center}
\end{lemma}


\begin{proof}
To prove the lemma, we compute certain slice representations, cf.~Section~2.3 in \cite{kollross}. We
assume in the following that the maximal reducible groups are standardly embedded as block diagonal
matrices, cf.\ (\ref{StdEmbedd}). Assume first $k \le \ell$. We compute a slice representation of the
action of $H = \SO(k)\times\SO(n{-}k)$ on $G/K = \SO(n) / \SO(\ell)\times\SO(n{-}\ell)$. The connected
component of the isotropy group is the group $(H \cap K)_0 = \SO(k)\times\SO(\ell-k)\times\SO(n-\ell)$;
it acts on the normal space
$$N_{\e K}(H \cdot \e K) = \left\{ \left. \left( \begin{matrix}
  0 & 0 & M \\
  0 & 0 & 0 \\
  -M^t & 0 & 0 \\
\end{matrix} \right) \right| M \in \R^{k \times n-\ell}
\right\} \subset \so(n) $$ by the tensor product of the two standard representations of the first and
the last factor i.e.\ $\SO(k)\otimes\SO(n-\ell)$. Assume $U \subseteq H$ is a closed subgroup acting
polarly on $G/K$. By Lemma~\ref{OrbitEqPolarSubGr} the connected component of the isotropy group $U \cap
K$ of the $U$-action must contain the product of the first and last factor $L = \SO(k)\times\SO(n-\ell)
\subset H \cap K$, except possibly in cases $k = 2,\,3$ and $n-\ell = 7,\,8$ which will be treated
below. Since this argument also holds for any conjugate subgroup $hUh^{-1}$, $h \in H$, it follows that
$U$ contains $h L h^{-1}$, for all $h \in H$, hence $U$ contains the subgroup generated by $\{ h x
h^{-1} \mid h \in H,\, x \in L\}$, which is the minimal normal subgroup of~$H$ containing~$L$ and we
conclude $H = U$.
\par
Using an analogous argument in the case of $H = \SUxU{k}{n-k}$ acting on $G/K = \SU(n) /
\SUxU{\ell}{n-\ell}$ we see that the only polar subaction is the action of $U = \SU(k)\times\SU(n-k)$
except in the case of $U = \SU(k)\times\SU(k)$ acting on $G/K = \SU(2k)/\SUxU{k}{k}$, where the slice
representation of the $U$-action is non-polar, see Lemma~\ref{OrbitEqPolarSubGr}. The same argument also
works for the actions on the quaternionic Grassmannians and for the case $\ell \le k$.
\par
It remains to study the case where a slice representation of the $U$-action is given by the first three
rows of Lemma~\ref{OrbitEqPolarSubGr}. It follows from Lemma~\ref{GrassIsotropy} that this can only
happen if $U$ is obtained from $H$ by replacing an $\SO(7)$-factor with $\LG_{2}$ or replacing an
$\SO(8)$-factor with $\Spin(7)$. A dimension count shows that the actions obtained in this fashion are
orbit equivalent to the respective $H$-action.
\end{proof}


\begin{table}[h]
\begin{tabular}{|l|c|c|}\hline
\hl \str Action type & Effectivized slice representation & Kernel \\
\hline\hline \str
A\,I-II & $\SU(n) \times \SU(n) / \Delta \SU(n)$ & $\U(1)$
\\\hline \hl
A\,I-III & $\SO(n)/\SO(k)\x\SO(n-k)$ & $$ \\\hline \hl
A\,II-III, $k$ even & $\Sp(n)/\Sp\left(\frac{k}{2}\right)\x\Sp\left(n-\frac{k}{2}\right)$ & $$ \\ \hline
\hl
A\,II-III, $k$ odd  & $\Sp(n-1)/\Sp\left(\frac{k-1}{2}\right)\x\Sp\left(n-\frac{k+1}{2}\right)$ & $\U(1)$ \\
\hline \hl
A\,III-III & $\SU(k+n-\ell)/\SUxU{k}{n-\ell}$ & $\SU(\ell-k)$ \\
\hline \hl
BD\,I-I & $\SO(k+n-\ell)/\SO(k)\x\SO(n-\ell)$ & $\SO(\ell-k)$ \\
\hline \hl
C\,I-II & $\SU(n)/\SUxU{k}{n-k}$ & \U(1) \\ \hline \hl
C\,II-II & $\Sp(k+n-\ell)/\Sp(k)\x\Sp(n-\ell)$ & $\Sp(\ell-k)$ \\
\hline \hl
D\,I-III, $k$ even & $\SU(n) / \eS\left(
\U\left(\frac{k}{2}\right)\x\U\left(n-\frac{k}{2}\right)\right)$ & \U(1) \\\hline \hl
D\,I-III, $k$ odd  & $\SU(n-1) / \eS \left( \U \left( \frac{k-1}{2} \right) \times \U \left(n-\frac{k-1}{2} \right) \right)$ & \U(1) \\
\hline \hl
D\,III-III' & $\SO(4n-4)/\U(2n-1)$ & $$ \\ \hline \hl
D$_4$ I-I' , $k{ = }\ell{ = }3$ & $\LG_2/\SO(4)$ & $$ \\ \hline \hl
E\,I-II    & $\LF_4/\Sp(3)\cdot\Sp(1)$ & $$ \\ \hline \hl
E\,I-III   & $\Sp(4)/\Sp(2)\times\Sp(2)$ & $$ \\ \hline \hl
E\,I-IV    & $\SU(6)/\Sp(3)$ & $\Sp(1)$ \\\hline \hl
E\,II-III  &  $\SU(6)/\eS(\U(2)\x\U(4))$ & $\Sp(1)$ \\ \hline \hl
E\,II-IV   &  $\Sp(4)/\Sp(3)\times\Sp(1)$ & $$ \\ \hline \hl
E\,III-IV  & $\LF_4/\Spin(9)$ & $$ \\ \hline \hl
E\,V-VI    & $\SU(8)/\SUxU{4}{4}$ &  \\ \hline \hl
E\,V-VII   &  $\SU(8)/\Sp(4)$ & \\ \hline \hl
E\,VI-VII  & $\SU(8)/\SUxU{2}{6}$ &  \\ \hline \hl
E\,VIII-IX & $\SO(16)/\U(8)$ & \\ \hline \hl
F\,I-II    & $\Sp(3)/\Sp(2)\times\Sp(1)$ & $\Sp(1)$ \\ \hline
\end{tabular}
\bl\caption{Slice representations of Hermann actions}\label{THermannSliceRep}
\end{table}


\begin{theorem}\label{SubHermannClass}
    Let $G$ be a connected simple compact Lie group and let $H$ and
    $K$ be two non-conjugate connected symmetric subgroups of $G$
    such that the cohomogeneity $r$ of the $H$-action on $G/K$ is
    $\ge 2$.
    Let $U \subseteq H$ be a closed connected nontrivial subgroup acting polarly on~$G/K$.
    Then the action of $U$ on $G/K$ is
    orbit equivalent to the hyperpolar $H$-action on $G/K$.
    \par
    Furthermore, $U \neq H$ if and only if $U$ is as described in
    Lemma~\ref{ReducibleGrass} or the triple $(U,G,K)$ is one of
    $(\SU(2n-2k-1) \times \SU(2k+1),\, \SU(2n),\, \Sp(n))$;
    $(\SU(n),\, \Sp(n),$ $\Sp(k)$ $ \times \Sp(n - k))$;
    $(\SU(n),\, \SO(2n),\, \SO(k) \times \SO(2n-k))$, $k < n$;
    $(\SU(n),\, \SO(2n),$ $\a(\U(n)))$;
    $(\Spin(10),\, \LE_6,\, \SU(6) \cdot \Sp(1))$;
    $(\LE_6,\, \LE_7,\, \SO'(12) \cdot \Sp(1))$.
\end{theorem}


\begin{proof}

To prove the theorem, we use the explicit knowledge of slice representations of Hermann actions as given
in Tables~\ref{THermannActions} and~\ref{THermannSliceRep}. The first table is a list of all Hermann
actions (i.e.\ a list of all pairs $(H,K)$ of non-conjugate symmetric subgroups of the simple compact
Lie groups $G$ up to automorphisms of~$G$); Table~\ref{THermannSliceRep} contains information about one
irreducible slice representation of each action; slice representations of Hermann actions are
s-representations by Lemma~\ref{SubHermann} and so each representation is described by a symmetric space
$G'/K'$ whose isotropy representation $\chi(G',K')$ is equivalent to the slice representation on the Lie
algebra level; in the third column, the (local isomorphism type of the) kernel of the slice
representation is given. It is straightforward to determine these slice representations for actions on
the classical symmetric spaces, for the exceptional symmetric spaces one may use the technique described
in Remark~\ref{MaxRkSubgroups}, cf.\ also \cite{kollross}, Prop.~3.5. Note that actually two different
actions are given in each row of the table, i.e.\ the action of $H$ on $G/K$ and the action of $K$ on
$G/H$; they have the same isotropy subgroups and slice representations.

\par

Assume now that $H$ and $K$ are symmetric subgroups of the simple compact Lie group $G$ and $U \subseteq
H$ is a closed connected subgroup acting polarly on $G/K$ and such that the hyperpolar action of $H$ on
$G/K$ is of cohomogeneity $r \ge 2$. The subactions of the types A\,III-III, BD\,I-I, and C\,II-II were
treated in Lemma~\ref{ReducibleGrass}.

\par

The slice representations given by the table are irreducible and non-transitive on the sphere, since we
assume the cohomogeneity is $\ge 2$, thus we may apply Theorem~\ref{MaxPolRep}. The isotropy group of
the $U$-action at $p = \e K \in G/K$ is $U_{p} = U \cap H_{p}$. The representation of $H_{p} = H \cap K$
on the normal space $N_{p}(H \cdot p)$ restricted to $U_{p}$ occurs as a submodule in the slice
representation of $U_{p}$ on $N_{p}(U \cdot p)$, and is therefore polar. By Theorem~\ref{MaxPolRep}, the
$U_{p}$-action on $V = N_{p}(H \cdot p)$ is either orbit equivalent to the $H_{p}$-action or trivial.

We first show that the slice representation of $H_p$ restricted to~$U_p$ is non-trivial. If it is
trivial, then $V$ is contained in the tangent space of a section through $p$ and we obtain a
contradiction with Theorem~\ref{ProductOfSpheres} since $V\subset\g$ is a Lie triple system
corresponding to a totally geodesic submanifold of $G/K$ which is isometric to an irreducible symmetric
space of higher rank, see the proof of Lemma~\ref{SubHermann}. We may therefore assume the
$U_{p}$-action on $V$ is locally orbit equivalent to the irreducible polar representation of~$H_{p}$.
From Table~\ref{THermannSliceRep} we see that we may assume the slice representation of $H_{p}$ on $V$
is not equivalent to one of the first three items in Lemma~\ref{OrbitEqPolarSubGr}, except in the case
of A\,I-III, which will be treated separately.

Then it follows that the group $U_p$ contains the (component of the) isotropy group $H_p$ or a subgroup
as described in the 4th, 5th, and 6th item of Lemma~\ref{OrbitEqPolarSubGr}. In these cases there exists
a uniquely\footnote{In the first three items of Lemma~\ref{OrbitEqPolarSubGr}, an orbit equivalent
subgroup is only unique up to conjugation.} determined connected subgroup $L \subset H_p$ which is
minimal with respect to the property that the $L$-action on $V$ is orbit equivalent to the $H_p$-action.
Note that this argument actually shows that for any $h \in H$ also $h U h^{-1}$ contains the subgroup $L
\subset H_p$, hence $U$ contains all groups $h^{-1} L h$ conjugate to $L$ in $H$. We conclude that $U$
contains the subgroup $\hat L$ generated by $\left\{ h \ell h^{-1} \mid \ell \in L,\, h \in H \right\}$,
i.e.\ the minimal normal subgroup of $H$ containing $L$.

\paragraph{\em Subgroups of codimension one}

Let us first consider the case where $U \subset H$ is a subgroup of codimension one, i.e.\ $H = U \cdot
\U(1)$, then we have that either $U$ acts transitively on the $H$-orbit through $p$, in which case the
$U$-action and the $H$-action are orbit equivalent, or $U$ acts with cohomogeneity one on the orbit $H
\cdot p$ in which case we arrive at a contradiction since a section $\Sigma$ through $p$ of the
$U$-action contains the flat section $\Sigma_{0}$ of the $H$-action, on whose tangent space
$\T_{p}\Sigma_{0}$ the Weyl group of the irreducible slice representation still acts irreducibly when
restricted to the Weyl group of the $U_p$-representation, so $\Sigma$ would be either flat,
contradicting Proposition~\ref{HyperPMax}, or an irreducible symmetric space of dimension $r+1$ and rank
$r \ge 2$, which does not exist.

\paragraph{\em Subactions of Hermann actions on exceptional symmetric spaces}

Assume the subgroup $U \subset H$ acts polarly on the symmetric space $G/K$. One can see from
Table~\ref{THermannActions} that the group $H$ has either one or two simple factors if it is semisimple
or it is the product $H = H' \cdot \U(1)$ of a one-dimensional abelian and a simple factor. Since $U$
contains the nontrivial normal subgroup $\hat L$ of $H$ it follows that $U = H$ if $H$ is simple; if $H
= H' \cdot \U(1)$ then $U$ contains $H'$ (since $\dim \hat L >1$) and the $U$-action is orbit equivalent
to the $H$-action by the argument above since then $U \subset H$ is a subgroup of codimension one; in
those cases where $H$ is a product of two simple factors, comparison of the Tables~\ref{TSymmSpaces} and
\ref{THermannSliceRep} shows that in each case, except for $H = \SU(6) \cdot \SU(2)$, $G / K = \LE_6 /
\Spin(10) \cdot \U(1)$, the normal subgroup $\hat L$ contains both simple factors of $H$ and it follows
that $H = U$. Consider the action of $H = \SU(6) \cdot \SU(2)$ on $G / K = \LE_6 / \Spin(10) \cdot
\U(1)$; in this case it follows from the data given in Table~\ref{THermannSliceRep} only that $\hat L$
contains the $\SU(6)$-factor of~$H$. An explicit calculation as described in Remark~\ref{MaxRkSubgroups}
shows that the embedding of the connected component of the isotropy group $(H \cap K)_0 = \U(1) \cdot
\SU(4) \cdot \SU(2) \cdot \SU(2)$ into $H = \SU(6) \cdot \SU(2)$ is such that the $\SU(2)$-factor in the
kernel of the slice representation lies in the $\SU(6)$-factor of~$H$, and the other $\SU(2)$-factor of
$(H \cap K)_0$, which acts nontrivially on the slice, coincides with the $\SU(2)$-factor of~$H$. From
this it follows that the actions $\SU(6)$ or $\SU(6) \cdot \U(1)$ on $G / K$ have a slice representation
with two equivalent nontrivial modules and are therefore not polar.

\paragraph{\em Subactions of Hermann actions on classical symmetric spaces}

The cases A\,I-II, A\,II-III, C\,I-II, D\,I-III, D\,III-III', and D$_4$\,I-I' can be handled in a
similar way as the subactions on the exceptional spaces. One can see from Table~\ref{THermannSliceRep}
that $\hat L$ contains every simple factor of $H$. For the case of D$_4$~I-I', i.e.\ subactions of
$\Spin(5) \cdot \Spin(3) \cong \Sp(2) \cdot \Sp(1)$ on $\SO(8) / \SO(3) \times \SO(5)$, the slice
representation was explicitly computed in \cite{kollross}, p.~592-593, and it follows that $\hat L$ is
not contained in one of the simple factors of $H$, thus $H = U$.

It remains to study the case A\,I-III. For the slice representation of this action there are in some
cases orbit equivalent polar subgroups as given in the first three items of
Lemma~\ref{OrbitEqPolarSubGr}; otherwise the argument is as above. Assume $H = \SO(n)$, $G = \SU(n)$, $K
= \SUxU{k}{n-k}$, where $(n,k) = (9,2)$, $(10,2)$ or $(11,3)$. Let us first consider the $H$-action on
$G/K$. The connected component of the isotropy subgroup at~$\eK$ is $\SO(k) \times \SO(n-k)$. It follows
that $U$ must contain the group given in the right column of the table in Lemma~\ref{OrbitEqPolarSubGr}
and it follows from Lemma~\ref{Intermediate} below that either $U = H$ or $U \subseteq H_p$, but in the
latter case the $U$-action on $G / K$ has a fixed point. Finally, consider the $K$-action on $G/H$,
i.e.\ assume a closed connected subgroup $U \subset K = \SUxU{k}{n-k}$ acts polarly on $\SU(n)/\SO(n)$;
it follows by the arguments above that $U$ contains a subgroup $L$ conjugate to $\SO(2)\times \LG_{2}$,
if $K = \SUxU{2}{7}$, or $\SO(k)\times \Spin(7)$, if $K = \SUxU{k}{8}$; all possibilities for the group
$U$ are given by Lemma~\ref{Intermediate}. It follows that the slice representation $V|_{U_p} \oplus
\chi(K,U)|_{U_p}$ of the $U$-action on $G/H$ is non-polar by~\cite{bergmann} if $U$ does not contain
both simple factors of $K$, thus the codimension of $U$ in $K$ is at most one and we conclude that the
$U$-action is orbit equivalent to the $K$-action.
\par
To prove the last part of the theorem, one can easily determine all proper closed subgroups~$U$ of~$H$
whose action on $G / K$ is orbit equivalent to the $H$-action on $G / K$ by using the information from
Table~\ref{THermannSliceRep}.
\end{proof}


For the proof of Theorem~\ref{SubHermannClass} we used the following simple Lemma.
\begin{lemma}\label{Intermediate}
For the following inclusions of compact Lie groups $A \subset B \subset C$, the intermediate subgroups
$B$ are unique in the following sense: If $B' \subset C$ is a closed connected subgroup such that $A
\subsetneqq B' \subsetneqq C$ then $B' = B$.
\begin{align*}
\LG_2\subset\SO(7)&\subset\SU(7);\\
\Spin(7)\subset\SO(8)&\subset\SU(8);\\
\SO(k)\times\LG_2\subset\SO(k)\times\SO(7)&\subset\SO(7+k),\quad k \in \Nat;\\
\SO(k)\times\Spin(7)\subset\SO(k)\times\SO(8)&\subset\SO(8+k),\quad k \in \Nat.
\end{align*}
\end{lemma}
\begin{proof}
It is easily checked in each case that the representation $\chi(C,A)$ splits into the irreducible
modules $\chi(B,A)$ and $\chi(C,B)|_A$. Note that $\chi(\SU(7),\SO(7))|_{\LG_2}$ is equivalent to the
irreducible $27$-di\-men\-sio\-nal representation of~$\LG_2$ and $\chi(\SU(8),\SO(8))|_{\Spin(7)}$ is
equivalent to the irreducible $35$-di\-men\-sio\-nal representation of~$\Spin(7)$, cf.\ Table~1, p.~364
of~\cite{dynkin2}.
\end{proof}


\section{Actions of non-simple irreducible groups}
\label{NonSimpleIrr}

In the following, we will assume that $G$ is a simple classical compact Lie group $G = \SO(n)$, $\SU(n)$
or $\Sp(n)$ and $K$ is a symmetric subgroup such that $\rk(G / K) \ge 2$. We will classify all closed
connected subgroups $H \subset G$ such that $H$ acts polarly on $G / K$.

The symmetric quotient spaces of the simple classical compact Lie groups of rank~$\ge 2$ which are not
locally isometric to one of the Grassmannians $\G_{k}(\R^n)$, $\G_{k}(\C^n)$, $\G_{k}(\H^n)$ are the
following:
\begin{equation}\label{StructSpaces}
\begin{array}{ll}
\SO(2m)/\U(m), & m \ge 5; \\
\SU(m)/\SO(m), & m \ge 3,\quad m \neq 4; \\
\SU(2m)/\Sp(m), & m \ge 3; \\
\Sp(m)/\U(m), & m \ge 3.
\end{array}
\end{equation}
(Note that $\SO(8) / \U(4)$ is locally isometric to~$\G_2(\R^8)$, $\SU(4) / \SO(4)$ is locally isometric
to~$\G_3(\R^6)$ and $\Sp(2) / \U(2)$ is locally isometric to~$\G_3(\R^5)$, cf.~\cite{helgason}, Ch.~X,
\S6.4). In the sequel we will refer to these spaces as {\em ``structure spaces''} since they can be
interpreted as: a space of complex structures on $\R^{2n}$, spaces of real or quaternionic structures on
$\C^n$ or $\C^{2m}$, respectively, and a space of complex structures on $\H^n$.

We will first consider the maximal subgroups $H_{1} \subset G$; for classical groups $G$, they are given
by Propositions~\ref{MaxSubgrSOn}, \ref{MaxSubgrSUn}, and \ref{MaxSubgrSpn}. Note that for $\SO(n)$ the
subgroups (i) and (iii), for $\SU(n)$ the subgroups (i), (ii), (iii) and for $\Sp(n)$ the subgroups (i)
and (ii) are symmetric, thus the actions of these groups are Hermann actions. The remaining types of
subgroups are either given by tensor product representations or are simple irreducible subgroups. We
will also study certain subactions of cohomogeneity one or transitive Hermann actions

Henceforth we will refer to the following maximal connected subgroups of the classical Lie groups
(cf.~Section~\ref{TotGeodSubmf}) as {\em (maximal) tensor product subgroups}:


\begin{equation}
\label{TensorPSubgr}
\begin{array}{lcll}
  H = \SO(p) \otimes \SO(q) &\subset& G = \SO(pq),
  & p \ge 3,\, q \ge 3;\\
  H = \SU(p) \otimes \SU(q)&\subset& G = \SU(pq),
 & p \ge 3,\, q \ge 2;\\
 H = \SO(p) \otimes \, \Sp(q)&\subset& G = \Sp(pq),
 & p \ge 3,\, q \ge 1;\\
 H = \Sp(p)\: \otimes\, \Sp(q)&\subset& G = \SO(4pq),
 & p \ge 2,\, q \ge 1.\\
\end{array} \end{equation}


\begin{proposition}[Tensor product groups on ``structure spaces'']
\label{TensStruct}

Let $G$ be a simple compact classical Lie group, let $H_1$ be a maximal tensor product subgroup of $G$
as in (\ref{TensorPSubgr}), let $K$ be a structure subgroup as in~(\ref{StructSpaces}), and let $H
\subseteq H_1$ be a closed connected subgroup acting nontrivially on~$G / K$. Then the $H$-action on
$G/K$ is not polar.

\end{proposition}


\begin{proof}
There are a few exceptions remaining not excluded by Proposition~\ref{ClassDimBounds}:
\begin{equation}\label{TensorPAct1}
\begin{array}{rcl}
H = \SU(3) \otimes \SU(2) & \mbox{acting on} & \SU(6) / \Sp(3) =  G / K;\\
\SU(4) \otimes \SU(2) & \mbox{acting on} & \SU(8) / \Sp(4);\\
\Sp(3) \otimes \Sp(1) & \mbox{acting on} & \SO(12) / \U(6).\\
\end{array}
\end{equation}
We will apply Lemma~\ref{PolMinHered} to show that none of the actions~\ref{TensorPAct1} can have a
polar subaction. We use the information collected in Table~2 of \cite{kp} to determine slice
representations. In the case of the $H = \SU(3) \otimes \SU(2)$-action on $G / K = \SU(6) / \Sp(3)$, a
slice representation is $\Ad_{\SO(3)}\otimes\Ad_{\SU(2)}$, which is polar. However, an explicit
calculation shows that the normal space to a regular orbit is not a Lie triple system, thus by
Proposition~\ref{PolCrit}, the action is not polar. Let now $U \subset H$ be closed connected proper
subgroup acting polarly on~$G / K$. By Theorem~\ref{MaxPolRep} and Lemma~\ref{OrbitEqPolarSubGr}, the
above slice representation restricted to~$(U \cap K)_0$ is either trivial, leading to a contradiction by
Proposition~\ref{ClassDimBounds} since the slice is $9$-dimensional, or equivalent to the action of $H
\cap K$, in which case $U$ must contain $\SO(3) \cdot \SU(2)$. But since $\SO(3) \subset \SU(3)$ is
maximal connected we have that the $U$-action on $G / K$ has a fixed point and is non-polar by
Corollary~\ref{PolFixOnHomSp}.
\par
Let us consider the $\SU(4) \otimes \SU(2)$-action on $\SU(8) / \Sp(4)$, a slice representation is
$\Ad_{\SO(4)}\otimes\Ad_{\SU(2)}$, which is non-polar~\cite{bergmann} and polarity minimal, hence we may
apply Lemma~\ref{PolMinHered}~(ii). For the $\Sp(3) \otimes \Sp(1)$-action we find a slice
representation $\Rho_{2}({\Sp(3)})\otimes\R^2$ of~$\Sp(3) \otimes \U(1)$, which is also
non-polar~\cite{dadok} and polarity minimal, hence Lemma~\ref{PolMinHered}~(ii) also applies in this
case.
\end{proof}


\begin{proposition}\label{SubHerTensStruct}
Let $H_1 \subset G$ and $G / K$ be as in Table~\ref{TSubHerTensStruct} and let $H \subseteq H_1$ be a
closed connected subgroup acting nontrivially on~$G / K$. Then the $H$-action on $G / K$ is not polar.
\begin{table}[h] \rm
\begin{tabular}{|c|c|c|c|}

\hline \str Type & $H_1$ & $G / K$ & Range \\ \hline \hline

\strh A\,III-I & $\eS( (\U(p) {\otimes} \U(q)) {\times} \U(1))$ &
$\frac{\ts\SU(pq{+}1)}{\ts\SO(pq{+}1)}$ & $p \ge 3, q \ge 2$ \\
\hline

\strh A\,III-II & $\eS( (\U(2r{+}1) {\otimes} \U(2s{+}1)) {\times} \U(1))$ &
$\frac{\ts\SU(4rs{+}2r{+}2s{+}2)}{\ts\Sp(2rs{+}r{+}s{+}1)}$ & $r,s \ge 1$ \\
\hline

\strh A\,III-II & $\eS( (\U(2p) {\otimes} \U(q)) {\times} \U(2))$ &
$\frac{\ts\SU(2pq{+}2)}{\ts\Sp(pq{+}1)}$ & $p \ge 1, q \ge 3$ \\
\hline

\strh A\,III-II & $\eS( (\U(2r{+}1) {\otimes} \U(2s{+}1)) {\times} \U(3))$ &
$\frac{\ts\SU(4rs{+}2r{+}2s{+}4)}{\ts\Sp(2rs{+}r{+}s{+}2)}$ & $r,s \ge 1$ \\
\hline

\strh C\,II-I & $(\SO(p) {\otimes} \Sp(q)) {\times} \Sp(1)$ &
$\frac{\ts\Sp(pq{+}1)}{\ts\U(pq{+}1)}$ & $p \ge 3, q \ge 1$ \\
\hline

\strh D\,I-III & $\SO(2r{+}1) {\otimes} \SO(2s{+}1)$ &
$\frac{\ts\SO(4rs{+}2r{+}2s{+}2)}{\ts\U(2rs{+}r{+}s{+}1)}$ & $r,s \ge 1$ \\
\hline

\strh D\,I-III & $(\SO(2p) {\otimes} \SO(q)) {\times} \SO(2)$ &
$\frac{\ts\SO(2pq{+}2)}{\ts\U(pq{+}1)}$ & $p \ge 2, q \ge 3$ \\
\hline

\strh D\,I-III & $(\SO(2r{+}1) {\otimes} \SO(2s{+}1)) {\times} \SO(3)$ &
$\frac{\ts\SO(4rs{+}2r{+}2s{+}4)}{\ts\U(2sr{+}r{+}s{+}2)}$ &  $r,s \ge 1$\\
\hline

\end{tabular}
\bl\caption{Certain subactions of cohomogeneity one or transitive Hermann actions on ``structure
spaces''} \label{TSubHerTensStruct}
\end{table}

\end{proposition}
\begin{proof}
None of the subgroups~$H_1 \subset G$ fulfills the lower bound on its dimension given in
Proposition~\ref{ClassDimBounds}, except the action of~$\eS((\U(2) \otimes \U(3)) \times \U(2))$ on
$\SU(8) / \Sp(4)$. However, an explicit calculation shows that the cohomogeneity of this action is~$12$,
hence this action is non-polar and polarity minimal by Lemma~\ref{ProdSphDim}.
\end{proof}


\begin{proposition}[Tensor product subgroups on Grassmannians]
\label{TensGrass}

Let $G$, $H$, $K$ be as in (\ref{TensorPGrass}). Assume $n = 2,\dots,\lfloor \frac{pq}{2} \rfloor$ in
cases {\bf (a), (b), (c)}. In case {\bf (d)}, assume that $n = 2,\dots, 2pq $ and $n \ge 3$ if $q = 1$.
Assume further that in case {\bf (d)} $pq \neq 2$. Then the action of $H$ on $G / K$ is non-polar and
polarity minimal, i.e.\ any nontrivial action of a closed connected subgroup $U \subseteq H$ on $G/K$ is
non-polar.
\begin{equation}
\label{TensorPGrass}
\begin{array}{clll}
 \mbox{\bf (a)} &
  H = \SO(p) {\otimes} \SO(q),
 & G/K = \SO(pq)/\SO(n) {{\times}} \SO(pq{-}n),
 & p \ge 3,\, q \ge 3;\\
 \mbox{\bf (b)} &
 H = \SU(p) {\otimes} \SU(q),
 & G/K = \SU(pq)/\eS(\U(n){\times}\U(pq{-}n)),
 & p \ge 3,\, q \ge 2;\\
 \mbox{\bf (c)} &
 H = \SO(p) {\otimes} \Sp(q),
 & G/K = \Sp(pq)/\Sp(n) {{\times}} \Sp(pq{-}n),
 & p \ge 3,\, q \ge 1;\\
 \mbox{\bf (d)} &
 H = \,\Sp(p) {\otimes} \Sp(q), & G/K = \SO(4pq)/\SO(n)
{{\times}} \SO(4pq{-}n), & p \ge 2,\, q \ge 1.\\
\end{array} \end{equation}
\end{proposition}


\begin{proof}
We use the slice representations which were explicitly determined in Section~2.3 of~\cite{kollross}. In
each case, one finds a non-polar, polarity minimal, and almost effective submodule of the slice
representation, by Lemma~\ref{PolMinCriteria}~(iii) this implies that the slice representation is
polarity minimal. But then Lemma~\ref{PolMinHered}~(iii) shows that the $U$-action is non-polar, since
the normal space~$\N_{\eK}(H \cdot \eK)$ contains a Lie triple system corresponding to an irreducible
symmetric space of non-constant curvature in each case, as can be seen from the explicit description of
the normal spaces in~\cite{kollross}.
\end{proof}


\begin{proposition}\label{SubHerTensGrass}
Let $H_1 \subset G$ and $G / K$ be as in Table~\ref{TSubHerTensGrass} and let $H \subseteq H_1$ be a
closed connected subgroup acting nontrivially on~$G / K$. Then the $H$-action on $G / K$ is not polar.
\begin{table}[h] \rm
\begin{tabular}{|c|c|c|c|}

\hline \str Type & $H_1$ & $G / K$ & Range \\ \hline \hline

\strs A\,III-III & $\eS((\U(p) \otimes \U(q)) {\times} \U(1))$ & $\G_k(\C^{pq+1})$ & $p \ge 3,\, q \ge 2,\,2
\le k \le \lfloor \frac{pq+1}{2} \rfloor$ \\ \hline

\strs BD\,I-I & $\SO(p) \otimes \SO(q)$ & $\G_k(\R^{pq+1})$ & $p \ge 3,\, q \ge 3,\,3 \le k \le \lfloor
\frac{pq+1}{2} \rfloor$
\\ \hline

\strs BD\,I-I & $\Sp(p) \otimes \Sp(q)$ & $\G_k(\R^{4pq+1})$ & $p \ge
2,\, q \ge 1,\,3 \le k \le \lfloor \frac{4pq+1}{2} \rfloor$ \\
\hline

\strs C\,II-II & $ (\SO(p) \otimes \Sp(q)) {\times} \Sp(1)$ & $\G_k(\H^{pq+1})$ & $p \ge 3,\, q \ge 1,\,2 \le k
\le \lfloor \frac{pq+1}{2} \rfloor$ \\ \hline \hline

\strs A\,III-III & $\SO(p) \otimes \Sp(q)$ & $\G_{\ell}(\C^{2pq})$ & $p \ge 3,\,q \ge 1,\,\ell = 2,3$
\\ \hline

\strs BD\,I-I & $\U(p) \otimes \U(q)$ & $\G_3(\R^{2pq})$ & $p \ge 3,\, q \ge 2$
\\ \hline

\end{tabular}
\bl\caption{Certain subactions of cohomogeneity one or transitive Hermann actions on Grassmannians}
\label{TSubHerTensGrass}
\end{table}
\end{proposition}

\begin{proof}
Consider the first four items in Table~\ref{TSubHerTensGrass}. By Lemma~\ref{CommInvol}, the
corresponding Hermann action (indicated in the first column) has a totally geodesic orbit isometric to
$\G_k(\C^{pq})$, $\G_k(\R^{pq})$, $\G_k(\R^{4pq})$, and $\G_k(\H^{pq})$, respectively, on which $H$
acts. It follows from Proposition~\ref{TensGrass} and Lemma~\ref{PolarOnTotGeodSubmf} that $H$ acts
non-polarly except if $p = 3$, $q = 1$, $k = 2$ in case of the fourth action; however in this case the
normal space of a principal orbits is not a Lie triple system and $H_1$ acts non-polarly by
Proposition~\ref{PolCrit}; closed proper subgroups of~$H_1$ are excluded by
Proposition~\ref{ClassDimBounds}~(iii).
\par
Let us now consider the last two items of Table~\ref{TSubHerTensGrass}. Assume $H$ acts polarly and
nontrivially on $G / K$. Then the $H$-action has a singular orbit by Corollary~\ref{SingularOrbit}. As
can be seen from Table~\ref{TCohOneTrHermann}, the almost effective slice representation of the
$H_1$-action on~$G / K$ occurs also as a submodule of the isotropy representation~$\chi( H_1, H_1 \cap K
)$ of the $H_1$-orbit $H_1 \cdot \eK$. Thus the $H$-action on~$G / K$ is non-polar by
Proposition~\ref{PolMinCriteria}~(ii).
\end{proof}


\section{Subactions of simple irreducible groups}
\label{SubSimple}

We will now study simple irreducible maximal subgroups of the classical groups acting on the classical
symmetric spaces. We start with actions on the spaces~(\ref{StructSpaces}).

\begin{proposition}[Simple irreducible groups on ``structure spaces'']
\label{SimpleStruct}

Let $G$ be a simple compact classical Lie group $\SO(n)$, $\SU(n)$ or $\Sp(n)$ and let $\rho \colon H
\to G$ be an irreducible representation of corresponding (real, complex or quaternionic) type where $H$
is a simple compact Lie group and such that $\rho(H)$ is a maximal connected subgroup of $G$. Let $K
\subset G$ be a subgroup as in (\ref{StructSpaces}) such that $\rk(G / K) \ge 2$. Then the action of any
closed subgroup of $\rho(H)$ on $G / K$ is non-polar except for the Hermann actions of subgroups
conjugate to $\Spin(7) \subset \SO(8)$ on $\SO(8) / \U(4)$.
\end{proposition}


\begin{proof}
For the spaces $\SU(n) / \SO(n)$ and $\Sp(n) / \U(n)$ this follows directly from Lemmata~2.7, 2.8
of~\cite{kollross} and Proposition~\ref{ClassDimBounds}.
\par
Let us consider the spaces $\SO(n) / \U(\frac{n}{2})$, $n \ge 8$. By Proposition~\ref{ClassDimBounds} we
have that $\dim(H) \ge \frac{n^2}{4} - n$ if $H$ acts polarly on $\SO(n) / \U(\frac{n}{2})$. Hence
$\rho$ is a representation as described in Lemma~2.6 (i),~(iv) of~\cite{kollross} and all possibilities
for~$\rho$ are given in the table of \cite{kollross}, Lemma~2.8~(i). However, all of these
subgroups~$\rho(H) \subset \SO(n)$ are excluded by the dimension bounds given in
Proposition~\ref{ClassDimBounds}, except $\Spin(7) \subset \SO(8)$.
\par
For the spaces $\SU(n) / \Sp(\frac{n}{2})$, $n \ge 3$, it follows from Proposition~\ref{ClassDimBounds}
that $\dim(H) \ge \frac{n^2}{2} - 2n$ for a group $H$ acting polarly on $\SU(n) / \Sp(\frac{n}{2})$.
Thus $\rho$ is a representation as in Lemma~2.6 (ii),~(iv) of~\cite{kollross} and all such
representations~$\rho$ are given by the table in~\cite{kollross}, Lemma~2.8~(ii). However, none of the
simple groups there fulfills the necessary condition on its dimension given by
Proposition~\ref{ClassDimBounds}.
\end{proof}


We also need to consider certain subactions of cohomogeneity one or transitive actions.

\begin{lemma}
\label{SimpleStruct2}

\begin{enumerate}
\item

Let $k \in \{1,2,3\}$. Let $H$ be a simple compact connected Lie group and let $\rho \colon H \to
\SO(2n-k)$ be an irreducible representation of real type such that $\rho(H) \subset \SO(2n-k)$ is
maximal connected. Then any closed subgroup of $\rho(H) \times \SO(k)$ acts non-polarly on $\SO(2n) /
\U(n)$, except if $\rho$ is equivalent to the $7$-dimensional irreducible representation of~$\LG_2$ and
$k = 1$.

\item

Let $k \in \{1,2,3\}$. Let $H$ be a simple compact connected Lie group and let $\rho \colon H \to
\SU(2n-k)$ be an irreducible representation of complex type such that $\rho(H) \subset \SU(2n-k)$ is
maximal connected. Then any closed subgroup of $\eS( (\rho(H) \otimes \U(1) ) \times \U(k) )$ acts
non-polarly on $\SU(2n) / \Sp(n)$.

\item

Let $H$ be a simple compact connected Lie group and let $\rho \colon H \to \SU(n)$ be an irreducible
representation of complex type such that $\rho(H) \subset \SU(n)$ is maximal connected. Then any closed
subgroup of $\eS( (\rho(H) \otimes \U(1) ) \times \U(1))$ acts non-polarly on $\SU(n) / \SO(n)$.

\item

Let $H$ be a simple compact connected Lie group and let $\rho \colon H \to \Sp(n)$ be an irreducible
representation of real type such that $\rho(H) \subset \Sp(n)$ is maximal connected. Then any closed
subgroup of $\rho(H) \times \Sp(1)$ acts non-polarly on $\Sp(n) / \U(n)$.

\end{enumerate}
\end{lemma}

\begin{proof}
The proof is almost literally the same as the proof of Proposition~\ref{SimpleStruct}. For the cases~(i)
and (ii), we may use the tables in parts~(i) and (ii) of Lemma~2.8 in~\cite{kollross}. The only
representation not excluded by this argument are the $8$-dimensional spin representation of~$\Spin(7)$
and the $7$-dimensional representation of~$\LG_2$. However, in case of the actions of $\Spin(7) \times
\SO(2)$ and $\LG_2 \times \SO(3)$ on~$\SO(10) / \U(5)$, the normal space at a principal orbit is not a
Lie triple system and hence these actions are non-polar by Proposition~\ref{PolCrit}. Closed connected
subgroups of these groups can be shown to act non-polarly by the same argument or are excluded by
Proposition~\ref{ClassDimBounds}. The action of $\LG_2$ on $\SO(8) / \U(4)$ is orbit equivalent to the
action of $\Spin(7)$ on $\SO(8) / \U(4)$. The statements (iii) and (iv) follow directly from Lemmata~2.7
and~2.8 of~\cite{kollross}.
\end{proof}


We will now consider the maximal simple irreducible subgroups of the classical groups $\SO(n)$,
$\SU(n)$, $\Sp(n)$, given by irreducible representations of the real, complex, or quaternionic type,
respectively, and their actions on the corresponding Grassmannians $\G_{k}(\K^n)$, $\K = \R$, $\C$ or
$\H$, of higher rank. A necessary condition for polarity on the dimension of these subgroups is given by
Propositions~\ref{ClassDimBounds} and \ref{CoisoDimBound}. The irreducible representations of simple
compact Lie groups whose degrees are sufficiently low can be obtained from Lemma~2.6 of~\cite{kollross},
see also the tables in Lemma~2.8 and the Appendix of~\cite{kollross}. These representations are given by
Table~\ref{TLowDegRep}; the column marked with $k_{\max}$ indicates the maximal rank $k \le \lfloor
\frac{n}{2} \rfloor$ for which the necessary condition for polarity of an action on~$\G_k(\K^n)$ given
by Proposition~\ref{ClassDimBounds} is fulfilled. We only list such representations of complex or
quaternionic type where $k_{\max} \ge 2$ and representations of real type where $k_{\max} \ge 3$, since
polar actions on~$\G_2(\R^n)$ have been classified in~\cite{pth2}. It turns out there are no such
representations of complex or quaternionic type.


\begin{table}
\begin{displaymath}
\begin{tabular}{|c|l|r|c|c|c|} \hline \str
Group & Highest Weight & Degree & Type & Description & $k_{\max}$ \\
\hline \hline

$\LA_{2}$ & $(1,1)$ & $8$ & real & adjoint  & 4 \\\hline

$\LB_{3}$ & $(0,0,1)$ & $8$ & real & $\Spin(7)$  & 4 \\\hline

$\LB_{4}$ & $(0,0,0,1)$ & $16$ & real & F\,II  & 4 \\\hline

$\LC_{3}$ & $(0,1,0)$ & $14$ & real & A\,II  & 3 \\\hline

$\LF_{4}$ & $(1,0,0,0)$ & $26$ & real & E\,IV & 3 \\\hline

$\LG_{2}$ & $(1,0)$ & $7$ & real & $\Aut(\Ca)$ & 3 \\\hline

\end{tabular}
\end{displaymath}
\bl\caption{Representations of low degree} \label{TLowDegRep}
\end{table}

\paragraph{\em Subactions of $\Ad(\SU(3))$ on $\G_3(\R^8)$ and
$\G_4(\R^8)$}

The group $H = \SU(3)$ acts on its Lie algebra $\h$ by the adjoint representation $\Ad \colon H \to
\SO(\h)$ and we obtain a subgroup $\Ad(\SU(3)) \subset \SO(8)$ by identifying $\h$ with $\R^8$. We will
study the actions of this group on the Grassmannians $\G_k(\R^8)$, $k = 3,4$. Any closed connected
proper subgroups of~$\SU(3)$ are of dimension~$\le 4$ by Table~3 of~\cite{kollross} and are thus
excluded by Proposition~\ref{ClassDimBounds}.

The maximal connected subgroup $\SO(3) \subset \SU(3)$ leaves a $3$-dimensional subspace of $\R^8$
invariant, thus it is the connected component of an isotropy subgroup of the $\Ad(\SU(3))$-action on
$\G_3(\R^8)$. The slice representation contains the irreducible $7$-dimensional representation
of~$\SO(3)$ and is hence non-polar~\cite{dadok}.

Consider now the subgroup $\eS(\U(1){\times}\U(2)) \subset \SU(3)$. The action of $H$ on $\h = \R^8$ restricted
to $\eS(\U(1){\times}\U(2))$ leaves the $4$-dimensional linear subspace corresponding to $\suxu{1}{2} \subset
\h$ invariant. Thus the maximal connected $\eS(\U(1){\times}\U(2))$ coincides with the connected component of
the stabilizer of the $H$-action on $\G_4(\R^8)$. Its slice representation contains an $8$-dimensional
irreducible representation of~$\eS(\U(1){\times}\U(2))$ which is non polar~\cite{dadok} and the $H$-action
on~$\G_4(\R^8)$ is polarity minimal by Lemma~\ref{PolMinHered} (ii).

\paragraph{\em Subactions of $\Spin(7)$ on
 $\G_3(\R^{8})$ and $\G_4(\R^{8})$}

The subgroup $\Spin(7) \subset \SO(8)$ gives rise to a Hermann action since its Lie algebra is the fixed
point set of an involution of~$\so(8)$.

\paragraph{\em Subactions of $\Spin(9)$ on
$\G_3(\R^{16})$}

The action on $\G_3(\R^{16})$ was shown not to be polar in~\cite{kollross}, the slice representation
being equivalent to a $16$-dimensional non-polar irreducible representation of $\Sp(1)\cdot\Sp(2)$. Thus
by Lemma \ref{PolMinHered} (ii), no subaction of the $\Spin(9)$-action on~$\G_3(\R^{16})$ is polar.

\paragraph{\em Subactions of $\Spin(9)$ on
$\G_4(\R^{16})$}

Consider the subgroup~$H_0 = \Spin(4) \cdot \Spin(4) \subset \Spin(8) \subset \Spin(9)$, its action
on~$\R^{16}$ leaves a four-dimensional subspace~$V$ invariant. Since $H_0 \subset \Spin(9)$ is a
subgroup of maximal rank, it is easy to check that no other connected subgroup of $\Spin(9)$
containing~$H_0$ leaves $V$ invariant and thus $H_0$ is the connected component of an isotropy subgroup
of the $\Spin(9)$-action on~$\R^{16}$. The slice representation is equivalent to the sum of two
$12$-dimensional irreducible modules and is easily seen to be non-polar and polarity minimal.

\paragraph{\em Subactions of $\Sp(3)$ on~$\G_3(\R^{14})$}

An isotropy subgroup of the $\Sp(3)$-action on~$\G_3(\R^{14})$ is
$$\SO(3) \cdot \U(1) \subset \U(3) \subset \Sp(3).$$ Its
$16$-dimensional slice representation does not contain any trivial submodule and is therefore non-polar
by~\cite{bergmann}. Any proper subgroups of $\Sp(3)$ can be excluded by
Proposition~\ref{ClassDimBounds}.

\paragraph{\em Subactions of $\LF_4$ on~$\G_3(\R^{26})$}

The maximal connected subgroups of maximal rank in~$H_1 = \LF_4$ are, see \cite{oniscikBook}, Chapter~1,
\S~3.11.
\begin{equation}\label{F4MaxSubgroups}
 \Sp(3) \cdot \Sp(1), \quad
 \SU(3) \cdot \SU(3), \quad
 \Spin(9)
\end{equation}
We will determine a slice representation for the $H_1$-action on~$\G_3(\R^{26})$. According
to~\cite{dynkin1}, Table~25, p.~199, the subgroup~$\Spin(9)$ acts on~$\R^{26}$ by the direct sum of the
$9$-dimensional standard representation, the $16$-dimensional spin representation, and a one-dimensional
trivial representation. Thus if we further restrict this representation to the maximal connected
subgroup $\Spin(7) \cdot \SO(2)$ of $\Spin(9)$, a three-dimensional subspace~$W$ is left invariant and
it follows that $\Spin(7) \cdot \SO(2)$ is contained in an isotropy subgroup~$(H_1)_W$ of the
$H_1$-action on~$\G_3(\R^{26})$. The subgroup $\Spin(7) \cdot \SO(2) \subset \LF_4$ is of maximal rank
and it can be deduced from Table~25 of~\cite{dynkin1} that none of the groups in~(\ref{F4MaxSubgroups})
leaves a three-dimensional subspace of~$\R^{26}$ invariant. Hence $\Spin(7) \cdot \SO(2)$ is the
connected component of the isotropy subgroup~$(H_1)_W$ and the slice representation is, by a dimension
count, equivalent to~$\R^7 \oplus 2 \cdot \R^2 \otimes \R^8$, where $\Spin(7)$ acts on~$\R^8$ by the
spin representation, hence the action of~$\LF_4$ on~$\G_3(\R^{26})$ is non-polar and polarity minimal by
Proposition~\ref{PolMinCriteria}~(ii) and Lemma~\ref{PolMinHered}~(ii).

\paragraph{\em Subactions of $\LG_2$ on $\G_3(\R^7)$}

This is a cohomogeneity one action and its subactions will be treated in Section~\ref{SubCohOneTrans}.

\section{Polar actions on the exceptional spaces}
\label{ClassExcept}

In this section we will study those isometric actions on the exceptional symmetric spaces of compact
type which are subactions neither of Hermann actions nor of cohomogeneity one actions. It will turn out
that none of these actions is polar.

The maximal connected subgroups of the simple compact Lie groups were determined in~\cite{dynkin1},
Tables~12 and~12a, p.\ 150--151, and Theorem~14.1, p.\ 231. By Theorem~\ref{ProductOfSpheres} and
Lemma~\ref{ProdSphDim}, the cohomogeneity of a polar action on a symmetric quotient $G / K$ of a simple
Lie group $G$ is at most $\rk (G) + \rk(K)$. By Proposition~\ref{CoisoDimBound}, this estimate can be
further improved for Hermitian symmetric spaces, for which the cohomogeneity is at most~$\rk(G)$. From
this it follows by using the classification of symmetric spaces, see Table~\ref{TSymmSpaces}, that a
group acting polarly on a symmetric quotient $G / K$ with $\rk(G / K) \ge 2$ of one of the simple
exceptional Lie groups $G = \LE_6$, $\LE_7$, $\LE_8$, $\LF_4$, $\LG_2$ is at least of dimension $16$,
$47$, $96$, $20$, $4$, respectively. (We do need not consider the Cayley plane $\LF_4 / \Spin(9)$, since
it is of rank one.) First we would like to recall a method to describe certain subgroups of a
(semi)simple compact Lie group in terms of the root system, which is particularly useful for our
purposes, see~\cite{gg}, \S~8.3 and \cite{oniscikBook}, Ch.~1, \S~3.11.

\begin{remark}[Borel-De Siebenthal theory]
\label{MaxRkSubgroups}

Let $G$ be a connected compact simple Lie group. A subgroup $H \subset G$ is called a {\em subgroup of
maximal rank} if $\rk(H) = \rk(G)$, i.e.\ $H$ contains a maximal torus $T$ of~$G$. Consider the root
space decomposition $\g_{\C} = \g_0 + \sum_{\a \in R} \g_{\a}$, where $\g_0$ is the complexification of
the maximal abelian subalgebra of~$\g$ tangent to~$T$. Since the Lie algebra $\h_{\C}$ contains $\g_0$,
it is a $\g_0$-stable subspace of $\g_{\C}$, and it follows that $\h_{\C} = \g_0 + \sum_{\a \in S}
\g_{\a}$ where $S \subset R$ is a subset of the root system. Conversely, from suitable subsets $S
\subset R$, one may construct the Lie algebra of a subgroup $H \subset G$ of maximal rank, see
\cite{oniscikBook}, Chapter~1, \S~3.11. In particular, one can obtain all maximal connected subgroups of
maximal rank by such a construction. These are obtained by deleting certain vertices from the extended
Dynkin diagram, see \cite{oniscikBook} for details. The classification of all such subgroups up to
conjugation by automorphisms of~$G$ is given in Table~5, p.~64 of~\cite{oniscikBook} or in Table~12,
p.~150 of~\cite{dynkin1}, see also~\cite{gg}.
\par
Now assume $H$ and $K$ are both subgroups of maximal rank in $G$. Then we can use the above description
to obtain information about the $H$-action on the homogeneous space $G / K$, in particular, to compute
an isotropy algebra together with its slice representation. In fact, we may assume by conjugation of $K$
with a suitable element from~$G$ that both $H$ and $K$ contain a maximal torus~$T$ of~$G$. Then $\h_{\C}
= \g_0 + \sum_{\a \in S} \g_{\a}$ and $\k_{\C} = \g_0 + \sum_{\a \in S'} \g_{\a}$ for some subsets
$S,S'$ of the root system~$R$. In particular, the complexified isotropy algebra $(\h \cap \k)_{\C}$ of
the $H$-action on $G / K$ at $\e K$ is spanned by the Cartan algebra $\g_0$ and the root spaces
corresponding to the roots in the intersection $S \cap S'$ and it follows that $H \cap K$ is also a
subgroup of maximal rank in~$G$. On the other hand, the complexified normal space $(\h^{\perp} \cap
\k^{\perp})_{\C}$ of the $H$-orbit through $\e K$ is spanned by the root spaces corresponding to the
roots in $R \setminus \left( S \cup S' \right)$. Since $T$ is also a maximal torus of $H \cap K$, the
roots in $R \setminus \left( S \cup S' \right)$ are exactly the weights of the slice representation of
$H \cap K$ on the normal space $\h^{\perp} \cap \k^{\perp}$. It follows~\cite{jaenich} that the
$H$-orbit through $\e K$ is a singular orbit, since $T$ acts nontrivially on $\h^{\perp} \cap
\k^{\perp}$, in fact, the slice representation does not have any trivial submodules, since the
complexified normal space is spanned by root spaces corresponding to non-zero roots. In the special case
$H = K$, one obtains the isotropy representation $\chi(G, K)$ by this method.
\par
Note that if a subgroup $H \subseteq G$ is a fixed point set of an inner automorphism $\s$ of~$G$, i.e.\
$\s(x) = g x g^{-1}$, it is a subgroup of maximal rank, since the element $g = \exp(X)$, $X \in \g$,
lies in a maximal torus $T$ of $G$ and it follows that $\h_{\C} = \g_0 + \sum_{\{\a \mid X \in \ker
\alpha\}} \g_{\a}$ where $\g_0$ is the complexified Lie algebra of~$T$. (Conversely, if a subgroup of
maximal rank is the fixed point set of an automorphism, then the automorphism is inner.)
\par
Let $H$ and $K$ be two subgroups of maximal rank with common maximal torus $T$ as above. If both groups
are fixed point sets of involutions i.e.\ $H = G^{\s}$, $K = G^{\t}$, then it follows that the
involutions $\s$ and $\t$ commute, since they both act as either plus or minus identity on the root
spaces of $\g$. This shows that if $\s$ and $\t$ are two inner involutions of a simple compact Lie group
$G$, then $\t$ is conjugate to an involution which commutes with $\s$, cf.\ \cite{conlon}.

\end{remark}

\subsection{Symmetric spaces of $\LE_6$}
The maximal connected non-symmetric subgroups of $\LE_6$ of dimension~$\ge 16$ are $\SU(3) \cdot \SU(3)
\cdot \SU(3)$ and $\LG_2^1 \cdot {\LA_2^2}''$, see~\cite{dynkin1}. (The upper indices denote the Dynkin
index of subgroups and the primes are used to distinguish non-conjugate subgroups of the same Dynkin
index). By a dimension count, no closed subgroup of these groups acts polarly on the spaces $\LE_6 /
(\Sp(4)/\{\pm1\})$, $\LE_6 / \SU(6) \cdot \Sp(1)$ or $\LE_6 / \Spin(10) \cdot \U(1)$.
\par
It remains to determine the polar actions on $\LE_6 / \LF_4$. We start with the group $H_1 = \SU(3)
\cdot \SU(3) \cdot \SU(3)$. The subgroup $\SU(3) \cdot \SU(3) \cdot \SU(3)$ is constructed from the
extended Dynkin diagram of $\LE_6$ as follows, cf.\ Remark~\ref{MaxRkSubgroups}.
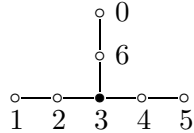
\begin{center}\vspace{20pt}
\begin{minipage}{100pt}
\begin{picture}(100,45)
\multiput(6,17)(16,0){5}{\circle{3}} \put(38,17){\circle*{3}} \multiput(8.5,17)(16,0){4}{\line(8,0){11}}
\multiput(38,19.5)(0,16){2}{\line(0,8){11}} \multiput(38,33)(0,16){2}{\circle{3}} \put(4,5){$1$}
\put(20,5){$2$} \put(36,5){$3$} \put(52,5){$4$} \put(68,5){$5$} \put(44,30){$6$} \put(44,46){$0$}
\end{picture}
\end{minipage}
\end{center}
The vertices numbered $1,\ldots,6$ correspond to the simple roots $\a_1,\ldots,\a_6$ of $\LE_6$ and the
vertex with number~$0$ represents $\a_0$, where $-\a_0$ is the maximal root. Now the group $H_1 = \SU(3)
\cdot \SU(3) \cdot \SU(3)$ arises from the extended Dynkin diagram if one deletes the central
vertex~$3$, i.e.\ it is the regular subgroup whose simple roots are $\a_1$, $\a_2$, $\a_4$, $\a_5$,
$\a_6$, $\a_0$. The subgroup $\LF_4 \subset \LE_6$ is the fixed point set of the diagram automorphism
$\s$ of $\LE_6$ which maps $\a_1 \mapsto \a_5$, $\a_2 \mapsto \a_4$, $\a_4 \mapsto \a_2$, $\a_5 \mapsto
\a_1$ and leaves $\a_3$ and $\a_6$ fixed. This automorphism $\s$ also leaves $\a_0$ fixed, since $- \a_0
= \a_1 + 2 \a_2 + 3 \a_3 + 2 \a_4 + \a_5 + 2 \a_6$, see \cite{oniscikBook}, Chapter 1, \S~3.11. It
follows that $\s$ also acts on $H_1$, i.e.\ trivially on one $\SU(3)$-factor (the one whose simple roots
are $\a_6$, $\a_0$) and by interchanging the other two $\SU(3)$-factors. Thus $H_1 \cap \LF_4$ is the
fixed point set $H_1^{\s}$ and is hence isomorphic to $\SU(3) \cdot \Delta \SU(3)$, where the $\Delta
\SU(3)$-factor is diagonally embedded into two of the $\SU(3)$-factors of $H_1$. Let us determine the
slice representation of the $H_1$-action on $M$, it is a submodule of $\chi(\LE_6, H_1)$ restricted to
$H_1 \cap \LF_4$. The real $54$-dimensional isotropy representation $\chi(\LE_6, H_1)$ is, after
complexification,
\begin{equation}\label{ChiE63SU3}
 \left( \begin{picture}(94,10)
\multiput(6,2)(16,0){6}{\circle{3}} \multiput(26,0)(32,0){2}{$\otimes$}
\multiput(8.5,2)(32,0){3}{\line(8,0){11}} \put(4,7){$1$} \put(36,7){$1$} \put(68,7){$1$}
\end{picture}  \right) \oplus \left(
\begin{picture}(94,10)
\multiput(6,2)(16,0){6}{\circle{3}} \multiput(26,0)(32,0){2}{$\otimes$}
\multiput(8.5,2)(32,0){3}{\line(8,0){11}}  \put(20,7){$1$} \put(52,7){$1$} \put(84,7){$1$}
\end{picture} \right),
\end{equation}
see~\cite{wolfIrr}, Corollary~13.2, i.e.\ the isotropy representation is equivalent to the action of
$\SU(3) \cdot \SU(3) \cdot \SU(3)$ on $\C^3 \otimes \C^3 \otimes \C^3$ by the tensor product of the
standard representations. If we restrict this representation to the subgroup $\SU(3) \cdot \Delta
\SU(3)$, it splits into the irreducible modules $(\C^3 \otimes \Sym^2 \C^3)$ and $(\C^3 \otimes
\Lambda^2 \C^3)$, where the first $\SU(3)$-factor acts on $\C^3$ and the $\Delta \SU(3)$-factor acts on
$\Sym^2 \C^3$ or $\Lambda^2 \C^3$, respectively. A dimension count shows that the real $36$-dimensional
slice representation of the $H_1$-action is equivalent to the first irreducible summand, hence it is
non-polar~\cite{dadok} and polarity minimal. We conclude that the $H_1$-action on~$G / K$ is non-polar
and polarity minimal by Lemma~\ref{PolMinHered}~(ii).
\par
Now consider subactions of $H_1 = \LG_2^1 \cdot {\LA_2^2}''$ on $G / K = \LE_6 / \LF_4$. We determine a
slice representation of the $H_1$-action on~$G / K$. First observe that $\LF_4$ contains the subgroup
$\LG_2^1$, according to Table~39 of~\cite{dynkin1}, p.~233. Since the subgroup $\LF_4 \subset \LE_6$ has
Dynkin index $1$, it follows that $\LG_2^1 \subset \LE_6$ also has Dynkin index~$1$, see~\cite{dynkin1},
Ch.~I, \S~2. By Table~25 of~\cite{dynkin1}, p.~200, there is only one conjugacy class of subgroups
isomorphic to~$\LG_2$ of Dynkin index~$1$ in $\LE_6$ and it follows that an isotropy subgroup $(H_1)_x$
of the $H_1$-action on~$G / K$ contains~$\LG_2^1$.  The homogeneous space $G / H_1$ is strongly isotropy
irreducible, see~\cite{wolfIrr}, Theorem~3.1, p.~66, and its isotropy representation decomposes into
$8$~equivalent $7$-dimensional irreducible modules when restricted to~$\LG_2$. Thus the dimension of the
normal space~$\N_x (H_1 \cdot x)$ to the $H_1$-orbit through~$x$ is a multiple of~$7$. The only
possibility is a $21$-dimensional slice representation which splits into $3$~irreducible $7$-dimensional
modules when restricted to~$\LG_2^1$. By~\cite{dadok} and \cite{bergmann}, such a representation is
non-polar and it is polarity minimal by Lemma~\ref{PolMinCriteria}, part~(i) or~(ii), and hence we can
apply Lemma~\ref{PolMinHered} to show that no closed subgroup~$H \subseteq H_1$ acts polarly on~$G / K$.

\subsection{Symmetric spaces of $\LE_7$}

The only maximal connected non-symmetric subgroup of dimension~$\ge 47$ is $H_1 = \LF_4^1 \cdot
{\LA_1^3}''$, see~\cite{dynkin1}. By Lemma~\ref{ProdSphDim} and Proposition~\ref{CoisoDimBound},
respectively, no closed subgroup of~$H_1$ acts polarly on $\LE_7 / (\SU(8)/\{\pm1\})$ or on $\LE_7 /
\LE_6 \cdot \U(1)$.

Let us determine an isotropy subgroup of the $H_1$-action on~$G / K = \LE_7 / \SO'(12) \cdot \Sp(1)$.
First observe that $\h_1$ contains a subalgebra~$\spin(9) \subset \Lf_4 \subset \h_1$. By Table~25
of~\cite{dynkin1}, p.~201, there is only one conjugacy class of subalgebras isomorphic to~$\spin(9)$ in
$\Le_7$ and it follows that, after conjugation, this subalgebra coincides with the subalgebra~$\spin(9)
\subset \spin(12) \subset K$. Thus there is an isotropy subgroup $(H_1)_x$ of the $H_1$-action on~$G /
K$ whose Lie algebra contains $\spin(9)$ as a subalgebra. The $64$-dimensional isotropy representation
of~$\LE_7 / \SO'(12) \cdot \Sp(1)$ decomposes into $4$ copies of the $16$-dimensional spin
representation when restricted to~$\Spin(9)$. Thus the dimension of an orbit $H_1 \cdot x$ is a multiple
of $16$ and it follows by a dimension count that the Lie algebra of $(H_1)_x$ is isomorphic to~$\spin(9)
\oplus \La_1$. The $32$-dimensional slice representation at the point~$x$ is the sum of two modules
equivalent to the $16$-dimensional spin representation of~$\Spin(9)$, hence non-polar and polarity
minimal and thus by Lemma~\ref{PolMinHered}, the $H_1$-action is non-polar and polarity minimal.
\par
Now consider the action of $H_1 = \LF_4^1 \cdot {\LA_1^3}''$ on $G / K = \LE_7 / \LE_6 \cdot \U(1)$. By
Table~25 of~\cite{dynkin1}, p.~201, there is only one conjugacy class of subgroups isomorphic to~$\LF_4$
in $\LE_7$ and it follows that an isotropy subgroup $(H_1)_x$ of the $H_1$-action on $G / K$ contains
$\LF_4$. The space $G / H_1$ is strongly isotropy irreducible, see~\cite{wolfIrr}, Theorem~3.1, and its
isotropy representation decomposes into $3$~equivalent $26$-dimensional irreducible modules when
restricted to~$\LF_4$. This shows that the dimension of the normal space $\N_x (H_1 \cdot x)$ is a
multiple of $26$ and it follows by a dimension count that the connected component of the isotropy
subgroup $(H_1)_x$ is isomorphic to $\LF_4 \cdot \U(1)$. Thus the slice representation is non-polar and
polarity minimal and we can apply Lemma~\ref{PolMinHered}~(i) to show that no closed subgroup $H
\subseteq H_1$ acts polarly on $G / K$.

\subsection{Symmetric spaces of $\LE_8$} The maximal connected
subgroups of $\LE_8$ whose dimension is at least $96$ are symmetric. Hence any polar action on the
symmetric spaces $\LE_8 / \SO'(16)$ and $\LE_8 / \LE_7 \cdot \Sp(1)$ are subactions of Hermann actions.

\subsection{Actions on $\LF_4 / \Sp(3) \cdot \Sp(1)$}
All maximal connected subgroups of $\LF_4$ whose dimension is at least $20$ are symmetric. Thus any
polar action on the space $\LF_4 / \Sp(3) \cdot \Sp(1)$ is a subaction of a Hermann action.

\subsection{Actions on $\LG_2 / \SO(4)$}
The maximal connected subgroups of $\LG_2$ are
\begin{equation}\label{G2MaxSubGroups}
    \SO(4), \quad \SU(3), \quad \LA_1^{28},
\end{equation}
where $\LA_1^{28}$ is a maximal connected subgroup in~$\LG_2$ of type~$\LA_1$, cf. \cite{dynkin1}. If
the group $H$ acting on $M = \LG_2 / \SO(4)$ is contained in $\SO(4)$, then the action has a fixed
point. The group $\SU(3)$ acts with cohomogeneity one on~$M$. The only closed connected subgroup $H
\subset \SU(3)$ of~dimension~$\ge 4$ is $\SUxU{1}{2} \cong \U(2)$, whose action on~$M$ has a fixed
point. Subgroups of rank one are ruled out by a dimension count.


\section{Subactions of hyperpolar actions}
\label{SubactionsOfHyperpolar}

To complete the classification, it remains to study subactions of cohomogeneity one and transitive
actions. We will need the following lemmata to study subactions of Hermann actions. The first lemma
shows that the slice representations of a Hermann action are s-representations.

\begin{lemma}[Slice representations of Hermann actions]
\label{SubHermann} Let $G$ be a connected simple compact Lie group and let $\s$, $\t$ be two involutive
automorphisms of~$G$. Let $K = G^{\s}_0$ be the connected component of the fixed point set of $\s$ and
let $H_1 = G^{\t}_0$ be the connected component of the fixed point set of~$\t$. Consider the
$H_1$-action on~$G / K$. Then the exponential image $S = \Exp_{\eK} (\N_{\eK} (H_1 \cdot \eK))$ of the
normal space to the orbit through~$\eK$ is a totally geodesic submanifold locally isometric to a
symmetric space $G^{\s \circ \t} / G^{\s} \cap G^{\t}$ whose isotropy representation is on the Lie
algebra level equivalent to the slice representation of the $H_1$-action on $G / K$ at~$\eK$.
\end{lemma}

\begin{proof}
The Lie algebra $\g$ of~$G$ admits the two decompositions
\begin{equation}\label{TwoDecomp}
    \g = \k \oplus \pp = \h_1 \oplus \m_1,
\end{equation}
where $\pp$ and $\m_1$ are the $-1$-eigenspaces of $\s_*$ and $\t_*$, respectively. The normal space
$\N_{\eK} (H_1 \cdot \eK) = \pp \cap \m_1 \subseteq \pp$ is a Lie triple system and $\k \cap \h_1 \oplus
\pp \cap \m_1$ is the Lie algebra generated by $\pp \cap \m_1$; the isotropy algebra at $\eK$ of the
$H_1$-action on $G / K$ is just $\h_1 \cap \k$ and its action on the normal space $\pp \cap \m_1$ agrees
on the Lie algebra level with the isotropy representation of $G^{\s \circ \t} / G^{\s} \cap G^{\t}$,
which is locally isometric to~$S$ by Proposition~\ref{LieTripleSystems}.
\end{proof}

In the special case where the two involutions defining a Hermann action commute (possibly after
conjugation), the action has a totally geodesic orbit. The pairs of involutions on the compact simple
Lie groups for which this is the case have been determined in~\cite{conlon}. By
Lemma~\ref{PolarOnTotGeodSubmf}, a polar subaction acts also polarly on this totally geodesic orbit.

\begin{lemma}[Subactions of Hermann actions with commuting involutions]
\label{CommInvol} Let $G$, $K$, and $H_1$ be as in Lemma~\ref{SubHermann}. Assume in addition that $\s
\circ \t = \t \circ \s$. Let $H \subseteq H_1$ be a closed connected subgroup acting polarly on $M = G /
K$. Then the $H_1$-orbit $H_1 \cdot \e K = H_1 / H_1 \cap K$ through $\e K \in M$ is a totally geodesic
submanifold and $H$ acts polarly on the symmetric space~$H_1 \cdot \e K \cong H_1 / H_1 \cap K$.
\end{lemma}

\begin{proof}
Since $\s$ and $\t$ commute, we have the direct sum decomposition
\begin{equation}\label{CommInvDecomp}
    \g = (\k \cap \h_1) \oplus (\k \cap \m_1) \oplus (\pp \cap
\h_1) \oplus (\pp \cap \m_1).
\end{equation}
Consider now the $H_1$-action on the symmetric space $G / K$. We can identify $\pp$ with the tangent
space $\T_{\e K} G / K$. Then $\h_1 \cap \pp$ is the tangent space of the $H_1$-orbit through the point
$\e K$. Using the Cartan relations for the decompositions (\ref{TwoDecomp}), it is easy to verify that
$\h_1 \cap \pp$ is a Lie triple system. Hence the $H_1$-orbit through $\e K$ is totally geodesic by
Proposition~\ref{LieTripleSystems}. Clearly, the action of~$H$ leaves all $H_1$-orbits invariant and the
polarity of the $H$-action on~$H_1 \cdot \e K$ follows from Lemma~\ref{PolarOnTotGeodSubmf}.
\end{proof}

The following lemma is a just simple reformulation of the criterion for polarity given by
Proposition~\ref{PolCrit} in the special case of subaction of a Hermann action; it is, however, useful
in particular to study polar actions on the exceptional spaces since it enables us to test for polarity
on a subspace.

\begin{lemma}\label{SliceActionTrans}
Let $G$, $K$, $H_1$, $H$ and $M = G / K$ be as in Lemma~\ref{SubHermann}. Assume the group $H$ acts
transitively on the $H_1$-orbit through~$\eK$. Then the action of $H$ on $G / K$ is polar if and only if
the action of the action of~$H \cap K$ on~$S = \Exp_{\eK} (\N_{\eK} (H_1 \cdot \eK))$ is polar and
$[\nu, \nu] \perp \h$, where $\nu \subseteq \m_1 \cap \pp$ is a normal space to a principal orbit of the
slice representation of~$H \cap K$ on~$\m_1 \cap \pp$.
\end{lemma}


\section{Subactions of cohomogeneity one and transitive actions}
\label{SubCohOneTrans}

To finish the proof of our classification result it remains to consider subactions of cohomogeneity one
and transitive actions. Note that polar actions on the real Grassmannians~$\G_2(\R^n)$ of rank two were
completely classified in~\cite{pth2} and we will not consider any spaces locally isometric to them here.
We may also ignore all actions with a fixed point, since they are known to be hyperpolar by
Corollary~\ref{PolFixOnHomSp}, see also~\cite{brueck}, Theorem~2.2.

\begin{proof}
[Proof of Theorems~\ref{MainTheorem} and \ref{MainTheorem2}] \label{MainProofBegin}

Let $G$ be a connected simple compact Lie group and let $K$ be a symmetric subgroup such that $\rk(G/K)
\ge 2$. Assume the closed connected subgroup $H \subset G$ acts polarly on $M = G / K$.
\par
We have already completed the classification in the case where $G$ is an exceptional Lie group, hence it
remains the case where $G$ is one of the classical Lie groups $\SO(n)$, $n \ge 7$, $\SU(n)$, $n \ge 3$,
or $\Sp(n)$, $n \ge 2$. Then $H$ is contained in one of the maximal connected subgroups of~$G$ as
described in Propositions~\ref{MaxSubgrSOn}, \ref{MaxSubgrSUn}, and \ref{MaxSubgrSpn}. Thus at least one
of the following holds:
\begin{itemize}
 \item $H$ is contained in a maximal tensor product subgroup~(\ref{TensorPSubgr}) of~$G$.
 \item $H$ is contained in maximal connected simple irreducible subgroup of~$G$.
 \item $H$ is contained in a symmetric subgroup of~$G$.
\end{itemize}
If the first possibility holds, then the result follows from Propositions~\ref{TensStruct} and
\ref{TensGrass} except if $M$ is a Grassmannian and the tensor product subgroup acts with cohomogeneity
one on~$M$. The second possibility was studied in Section~\ref{SubSimple}, except for subactions of
cohomogeneity one or transitive actions. In the case where $H$ is contained in symmetric subgroup~$H_1$
of~$G$, i.e.\ the $H$-action on $M$ is a subaction of a Hermann action; the result follows from
Theorem~\ref{SubHermannClass}, under the assumption that the cohomogeneity of the $H_1$-action on $M$
is~$\ge 2$. Thus it remains the case where $H$ is a proper closed connected subgroup of $H_1 \subset G$
such that $H_1$ acts on~$M$ with cohomogeneity~$\le 1$.

It will turn out that all polar actions on~$G / K$ are hyperpolar, hence it follows from Corollary~2.12
of~\cite{hptt} that the sections are embedded submanifolds.

\paragraph{\em Subactions of\, {``}exceptional{''} cohomogeneity one and
transitive actions}


\begin{table}[h]
\begin{tabular}{*{4}{|c}|}
\hline \str
$H_1$ & $G$ & $K$ & Cohomogeneity \\
\hline
\hline \hl $\LG_2$ & $\SO(7)$ &  $\SO(5) \times \SO(2)$ & 0 \\
\hline \hl $\LG_2$ & $\SO(7)$ &  $\SO(4) \times \SO(3)$ & 1 \\
\hline \hl $\Spin(9)$ & $\SO(16)$ &  $\SO(14)\times\SO(2)$ & 1\\
\hline
\hl $\Sp(n)\cdot\Sp(1),\,n \ge 2$ & $\SO(4n)$ &  $\SO(4n{-}2)\times\SO(2)$ & 1 \\
\hline \hl $\SU(3)$ & $\LG_2$ &  $\SO(4)$ & 1 \\ \hline
\end{tabular}
\bl\caption{\rm ``Exceptional'' cohomogeneity one and transitive actions}\label{TExcCohOneTypeI}
\end{table}


Let us first consider the case where the $H_1$-action on $M = G / K$ is not of Hermann type. These
cohomogeneity one and transitive actions were determined in \cite{kollross}, Theorem~A and
\cite{oniscik}, respectively. We only consider the cases where $G / K$ is symmetric of rank~$\ge 2$;
these actions are given in Table~\ref{TExcCohOneTypeI}.

We examine these actions case by case. Assume first $H \subset \LG_2$ is acting on~$\G_2(\R^7)$ or
$\G_3(\R^7)$. Then $H$ is contained in one of the maximal connected subgroups (\ref{G2MaxSubGroups}).
Under the $7$-dimensional irreducible orthogonal representation of $\LG_2$, the first two
groups~$\SO(4)$ and $\SU(3)$ are mapped to reducible subgroups of~$\SO(7)$, thus the $H$-action is in
this case a subaction of a Hermann action. If $H$ is contained in $\SU(3)$, then the $H$-action is a
subaction of a cohomogeneity one Hermann action of type~BD\,I-I and will be treated on
page~\pageref{subBDI-I}. The third group can be excluded by Lemma~\ref{ProdSphDim}.

\par

Let us now consider subgroups of $\Spin(9)$, acting on $\G_2(\R^{16})$. The maximal connected subgroups
of $\SO(9)$ are, see Proposition~\ref{MaxSubgrSOn} and~\cite{dynkin1}.
\begin{align*}
\SO(8), \, \SO(7) \times \SO(2), \, \SO(6) \times \SO(3), \, \SO(5) \times \SO(4), \, \SO(3) \otimes
\SO(3), \, \LA_1^{60}.
\end{align*}
Let $H_1 \subset \Spin(9)$ such that $\pi(H_1) \subset \SO(9)$ is one of the above, where $\pi \colon
\Spin(9) \to \SO(9)$ is the double cover. We need to consider the image of $H_1$ under the spin
representation $\delta \colon \Spin(9) \to \SO(16)$. We have $\delta(\Spin(8)) \subset \SO(8) \times
\SO(8)$ and $\delta(\Spin(7) \cdot \Spin(2)) \subset \U(8)$, thus any subgroups of these are contained
in symmetric subgroups of~$\SO(16)$; the remaining subgroups can be excluded by a dimension count, see
Lemma~\ref{ProdSphDim}.

\par

We do not need to consider subactions of the $\Sp(n) \cdot \Sp(1)$-action on $\G_2(\R^{4n}) = \SO(4n) /
\SO(4n-2) \times \SO(2)$, since polar actions on these spaces have been completely classified by Podest{\`a}
and Thorbergsson~\cite{pth2}.

\par

The last item in Table~\ref{TExcCohOneTypeI} was treated in Section~\ref{ClassExcept}.

\paragraph{\em Subactions of cohomogeneity one and transitive Hermann
actions}

\begin{table}[h]\small
\begin{tabular}{*{5}{|c}|c|}
\hline
\hl \str \stru Type & $G^{\sigma}$ & $G$ & $G^{\tau}$ & Coh. & $\frac{G^{\sigma\tau}}{G^{\sigma}\cap G^{\tau}}$  \\
\hline\hline

\stru A\,I-III & $\SO(p+1)$ & $\SU(p+1)$ & $\SUxU{p}{1}$ & $1$ & $\frac{\SO(p+1)}{\SO(p)}$ \\
\hline

\stru A\,III-II & $\eS(\U(2n-1){\times}\U(1))$ & $\SU(2n)$ & $\Sp(n)$ & $0$ &
$\frac{\Sp(n-1){\times}\U(1)}{\Sp(n-1){\times}\U(1)}$ \\\hline

\stru A\,III-II & $\eS(\U(2n-2){\times}\U(2))$ & $\SU(2n)$ & $\Sp(n)$ &
$1$ &   $\frac{\Sp(n)}{\Sp(n-1){\times}\Sp(1)}$ \\
\hline

\stru A\,III-II & $\eS(\U(2n-3){\times}\U(3))$ & $\SU(2n)$ & $\Sp(n)$ &
$1$ &   $\frac{\Sp(n-1){\times}\U(1)}{\Sp(n-2){\times}\U(1){\times}\Sp(1)}$ \\
\hline

\stru A\,III-III & $\eS(\U(a{+}b){\times}\U(1)))$ & $\SU(a{+}b{+}1)$ & $\eS(\U(a){\times}\U(b{+}1))$ &
$1$ &   $\frac{\eS(\U(a+1){\times}\U(b))}{\eS(\U(a){\times}\U(1){\times}\U(b))}$ \\
\hline

\stru BD\,I-I & $\SO(a{+}b)$ & $\SO(a{+}b{+}1)$ &
$\SO(a){\times}\SO(b{+}1)$ & $1$ &   $\frac{\SO(a+1){\times}\SO(b)}{\SO(a){\times}\SO(b)}$ \\
\hline

\stru C\,I-II & $\Sp(p){\times}\Sp(1)$ & $\Sp(p+1)$ & $\U(p+1)$ & $1$ &   $\frac{\U(p+1)}{\U(p){\times}\U(1)}$ \\
\hline

\stru C\,II-II & $\Sp(a+b){\times}\Sp(1)$ & $\Sp(a{+}b{+}1)$ & $\Sp(a){\times}\Sp(b+1)$ & $1$ &
$\frac{\Sp(a+1){\times}\Sp(b)}{\Sp(a){\times}\Sp(1){\times}\Sp(b)}$ \\ \hline

\stru D\,I-III & $\SO(2n{-}1)$ & $\SO(2n)$ & $\U(n)$ & $0$ & $\frac{\U(n-1)}{\U(n-1)}$ \\ \hline

\stru D\,I-III & $\SO(2n{-}2){\times}\SO(2)$ & $\SO(2n)$ & $\U(n)$ & $1$ &   $\frac{\U(n)}{\U(n-1){\times}\U(1)}$ \\
\hline

\stru D\,I-III & $\SO(2n{-}3){\times}\SO(3)$ & $\SO(2n)$ & $\U(n)$ & $1$ &   $\frac{\U(n-1)}{\U(n-2){\times}\U(1)}$ \\
\hline

\stru E\,II-IV & $\SU(6){\cdot}\Sp(1)$&$\LE_6$&$\LF_4$ & $1$ &   $\frac{\Sp(4)}{\Sp(3){\cdot}\Sp(1)}$ \\
\hline

\stru E\,III-IV & $\Spin(10){\cdot}\U(1)$&$\LE_6$&$\LF_4$ & $1$ &   $\frac{\LF_{4}}{\Spin(9)}$ \\
\hline

\stru F\,I-II & $\Sp(3){\cdot}\Sp(1)$&$\LF_4$&$\Spin(9)$& $1$ &   $\frac{\Sp(3){\cdot}\Sp(1)}{\Sp(2){\cdot}\Sp(1){\cdot}\Sp(1)}$ \\
\hline
\end{tabular}
\bl\caption{Cohomogeneity one and transitive Hermann actions.} \label{TCohOneTrHermann}
\end{table}

It now remains to study subactions of cohomogeneity one and transitive Hermann actions. These Hermann
actions are listed in Table~\ref{TCohOneTrHermann}. The column marked
with~$\frac{G^{\sigma\tau}}{G^{\sigma}\cap G^{\tau}}$ indicates (the local isometry type) of the
symmetric space ${G^{\sigma\tau}}/{G^{\sigma}\cap G^{\tau}}$ whose isotropy representation is equivalent
to one slice representation of the $G^{\s}$-action on~$G / G^{\t}$ (and of the $G^{\t}$-action on~$G /
G^{\s}$) by Lemma~\ref{SubHermann}. (The presentation may be non-effective, in particular for transitive
actions the space~${G^{\sigma\tau}}/{G^{\sigma}\cap G^{\tau}}$ is a noneffective presentation of a
zero-dimensional space.) We only have to consider actions on symmetric spaces of rank~$\ge 2$.

\paragraph{\bf A\,III-I}

Consider the action of $H_1 = \SUxU{p}{1}$ on $G / K = \SU(p+1) / \SO(p+1)$, $p \ge 2$. Assume a closed
connected subgroup $H \subset H_1$ acts polarly on~$G / K$. Then $H$ is contained in some maximal
connected subgroup $H_2$ of $H_1$. By Theorem~2.1 of~\cite{kollross}, either $H_2 = \SU(p)$ or $H_2 =
\eS((H_2' \otimes \U(1)) \times \U(1))$ where $H_2'$ is a maximal connected subgroup of~$\SU(p)$, see
Proposition~\ref{MaxSubgrSUn}. In the case of the $\SU(p)$-action on $G / K$, an explicit calculation
shows that the normal space to the orbit at~$\eK$ is not a Lie triple system, thus the $\SU(p)$-action
on $G / K$ is non-polar by Proposition~\ref{PolCrit}. Thus we may restrict our attention to the second
case, where we may further assume that $H_2' \subset \SU(p)$ is irreducible, since otherwise the
$H$-action is a subaction of a Hermann action whose cohomogeneity is~$\ge 2$. Assume first $H_2' =
\SO(p)$. Then the $H$-action on~$G / K$ is non-polar by Corollary~\ref{PolFixOnHomSp}, since the
$H_2$-action has a one-dimensional orbit. Now assume $H_2' = \Sp(p/2)$, $p \ge 4$, then one isotropy
group is $\U(p/2)$ and the slice representation is the adjoint representation of $\SU(p/2)$ plus the
standard representation of~$\U(p/2)$, see Table~\ref{THermannSliceRep}, which is non-polar and polarity
minimal~\cite{bergmann}. The normal space $\N_{\eK} H_2 \cdot \eK$ contains a Lie triple system
corresponding to an irreducible symmetric space of higher rank, thus the $H_2$-action on $G / K$ is
non-polar and polarity minimal by Lemma~\ref{PolMinHered}~(iii). Subgroups of $\eS(\U(p/\ell) \otimes
\U(\ell)) \times \U(1))$ are excluded by Proposition~\ref{SubHerTensStruct}, simple irreducible maximal
connected subgroups $H_2 \subset \SU(p)$ by Proposition~\ref{SimpleStruct2}~(iii).

\paragraph{\bf A\,III-II}

Let $H$ be a closed connected subgroup of $H_1 = \SUxU{2n-k}{k}$, $k = 1,2,3$ acting on $G / K = \SU(2n)
/ \Sp(n)$, $n \ge 3$. It is well known \cite{oniscik} that the action of $H = \SU(2n-1)$ is transitive
on~$G / K$.

We will now study actions of closed connected subgroups~$H$ in $H_2 = \eS( (H_2' \otimes \U(1) ) \times
\U(k))$. The cases where $H_2' \subset \SU(2n-k)$ is a simple irreducible or tensor product subgroup are
excluded by Lemma~\ref{SimpleStruct2} and Proposition~\ref{SubHerTensStruct}. Thus it remains to
consider the case where $H_2'$ is a symmetric subgroup of~$\SU(2n-k)$. Assume $H_2' = \SO(2n-k)$; if $k
= 1$ then $\SO(2n-1)$ acts on the symmetric space $M = \SU(2n) / \Sp(n)$, homogeneously presented as $M
= \SU(2n-1) / \Sp(n-1)$; an isotropy subgroup of this action is $H_2 \cap K = \U(n-1)$, its slice
representation is equivalent to the adjoint representation of~$\SU(n-1)$ plus the standard
representation on~$\C^{n-1} = \R^{2n-2}$, see Table~\ref{THermannSliceRep}. This representation is
non-polar~\cite{bergmann} and polarity minimal by Proposition~\ref{PolMinCriteria} and hence the
$H_2$-action on~$G / K$ is non-polar and polarity minimal by Lemma~\ref{PolMinHered}~(iii), since the
normal space contains an irreducible Lie triple system of higher rank. If $k = 2$ or $k = 3$, then a
slice representation of the $H_2$-action on $G / K$ contains a module equivalent to the isotropy
representation of~$\H \P^{n-1}$ or $\H \P^{n-2}$ restricted to $\U(n-1) \times \Sp(1)$ or $\U(n-2)
\times \Sp(1)$, respectively. This representation contains two equivalent modules and the $H$-action is
thus non-polar and polarity minimal by Lemma~\ref{PolMinHered} (iii), except if $n = k = 3$, a case
which can be handled by explicit calculations using the criterion in Proposition~\ref{PolCrit}.
\par
If $k = 2$ and $H_2' = \Sp(n-1)$, then the $H_2$-action has a one-dimensional orbit and the action of
any closed subgroup $H \subseteq H_2$ on $G / K$ is non-polar by Corollary~\ref{PolFixOnHomSp}. Thus we
are left with the case where $H_2' = \U(2n-k-\ell) \times \U(\ell)$. We may assume $k + \ell \le 3$
since otherwise the $H$-action on $G / K$ is a subaction of a Hermann action of cohomogeneity~$\ge 2$,
which were already treated in Section~\ref{SubHermannHighRk}. If $k = \ell = 1$ then we obtain the
cohomogeneity one actions of $H = \eS(\U(2n-2) \times \U(1) \times \U(1))$ and $\eS(\U(2n-2) \times
\U(1))$ on $G / K$, we have already seen that no further closed proper subgroup of these groups acts
polarly. In case $k + \ell = 3$ we have to consider closed connected subgroups~$H$ of~$H_2 =
\eS(\U(2n-3) \times (H_2'' \otimes \U(1)))$, where $H_2'' \subset \SU(3)$ is a maximal connected
subgroup. Since there are a number of subgroups $H \subseteq H_2$ acting with cohomogeneity two (in
these cases all slice representations are polar), we have to exclude them by explicit calculations
using~Proposition~\ref{PolCrit}.
\par
In case $n = k = 3$ there are additional maximal connected subgroups, i.e.\ $$H_2 = \{
(zA,z^{-1}\alpha(A)) \mid A \in \SU(3),\, z \in \C, \, \left| z \right| = 1 \} \subset \SUxU{3}{3},$$
where $\alpha \in \Aut(\SU(3))$, see Theorem~2.1 of~\cite{kollross}. If $\alpha$ is an outer
automorphism, e.g.\ given by complex conjugation, then the $H_2$-action has a one-dimensional orbit and
is non-polar and polarity minimal by Corollary~\ref{PolFixOnHomSp}. If $\alpha$ is an inner
automorphism, then a stabilizer component is $\U(1) \cdot \SO(3)$, the $9$-dimensional slice
representation is equivalent to $\R^3 \oplus \C^1 \otimes \R^3$ hence non-polar~\cite{bergmann} and
polarity minimal by Proposition~\ref{PolMinCriteria}~(iii). Thus the $H_2$-action on~$G / K$ is
non-polar and polarity minimal by Lemma~\ref{PolMinHered}.

\paragraph{\bf A\,II-III}

Consider the action of~$H_1 = \Sp(n)$ on $G / K = \SU(2n) / \SUxU{2n-2}{2}$, $n \ge 2$. Let $H \subseteq
H_1$ be a closed connected subgroup acting polarly on~$G / K$. The $H_1$-orbit $H_1 \cdot \e K \cong
\Sp(n) / \Sp(n-1) \times \Sp(1)$ is totally geodesic and $H$ acts polarly on this orbit by
Lemma~\ref{CommInvol}, the action being non-transitive by~\cite{oniscik}. The $H$-action on $H_1 \cdot
\e K$ has a singular orbit by Corollary~\ref{SingularOrbit} and we may assume by conjugation of~$H$
in~$H_1$ that $\e K$ lies in a singular orbit. From Table~\ref{TCohOneTrHermann} we read off that the
slice representation of the $H_1$-action on~$G / K$ is equivalent to the isotropy representation of the
symmetric space $H_1 \cdot \e K = H_1 / H_1 \cap K$. Thus the nontrivial slice representation of the
$H$-action on $H_1 \cdot \e K$, which is a submodule of the isotropy representation of~$H_1 / H_1 \cap
K$, also occurs as a submodule of the slice representation of~$H_1$ on $G / K$ restricted to~$H \cap K$.
We conclude that the slice representation of the~$H$-action on~$G / K$ contains two nontrivial
equivalent modules and is hence non-polar by~\cite{kollross}, Lemma~2.9.
\par
Consider the action of~$H_1 = \Sp(n)$ on $G / K = \SU(2n) / \SUxU{2n-3}{3}$, $n \ge 3$. Let $H \subseteq
H_1$ be a closed connected subgroup acting polarly on~$G / K$. Then $H$ is contained in a maximal
connected subgroup~$H_2$ of $H_1 = \Sp(n)$. We may assume that $H_2$ is irreducible, since otherwise the
$H$-action is a subaction of a Hermann action with cohomogeneity~$\ge 2$. If $H_2 = \U(n)$ then $H_2$ is
contained after conjugation in~$\SUxU{n}{n} \subset \SU(2n)$ and the $H$-action is also a subaction of a
Hermann action of cohomogeneity~$\ge 2$. The actions of maximal connected subgroups of type $\SO(q)
\otimes \Sp(p/q)$ have been treated in Proposition~\ref{SubHerTensGrass}. The actions of simple
irreducible subgroups $\rho(H)$, where $\rho \colon H \to \Sp(n)$ is an irreducible representation of
quaternionic type, have been excluded in Section~\ref{SubSimple}.

\paragraph{\bf A\,III-III}

Consider the action of  $H_1 = \SUxU{a + b}{1}$ on the complex Grassmannian $G / K = \SU(a + b + 1) /
\SUxU{a}{b + 1}$, $a \ge b \ge 1$, $a + b \ge 3$. Assume $H \subset H_1$ is a closed connected subgroup
acting polarly on~$G / K$. First note that the action of~$\SU(a + b)$ on~$G / K$ is orbit equivalent to
the $H_1$-action. Now assume $H \subseteq H_2 = \eS( (H_2' \otimes \U(1)) \times \U(1))$ is a closed
connected subgroup acting polarly on~$G / K$, where $H_2' \subset \SU(a + b)$ is a maximal connected
subgroup. We may assume that the standard representation of~$\SU(a + b)$ restricted to~$H_2'$ acts
irreducibly on~$\C^{a + b}$, since otherwise $H_2$ is contained in a subgroup of~$G$ conjugate
to~$\SUxU{k}{a + b + 1 - k}$ for $2 \le k \le a + b - 1$ and the $H_2$-action on $G / K$ is a subaction
of a Hermann action of cohomogeneity~$\ge 2$, which have already been treated in
Section~\ref{SubHermannHighRk}. It follows from Lemma~\ref{CommInvol} that the orbit $H_1 \cdot \eK$ is
a totally geodesic submanifold of~$G / K$ isometric to~$\SU(a + b) / \SUxU{a}{b}$ on which $H$ acts
polarly.
\par
Assume first that $H_2'$ is an irreducible symmetric subgroup of~$\SU(a + b)$, hence conjugate to either
$\SO(a + b)$ or $\Sp(\frac{a + b}{2})$. However, in the first case the isotropy subgroup of the
$H_2$-action at~$\eK$ is $\eS(\O(a) \times \O(b)) \cdot \U(1)$ and its slice representation is
equivalent to~$(\R^a \otimes \R^b) \oplus (\R^a \otimes \C^1)$, thus it is non-polar~\cite{bergmann} and
polarity minimal by Lemma~\ref{PolMinCriteria}. It follows from Table~\ref{TCohOneTrHermann} that the
normal space contains a Lie triple system corresponding to a totally geodesic submanifold isometric
to~$\C \P^a$, thus the $H_2$-action is non-polar and polarity minimal by Lemma~\ref{PolMinHered}~(iii),
since $a \ge 2$. Let us now consider the case where $a + b$ is even and $H_2'$ is conjugate to
$\Sp(\frac{a + b}{2})$, hence $a + b \ge 4$. The group $H$ acts polarly on the totally geodesic
$H_1$-orbit $H_1 \cdot \eK \cong \SU(a + b) / \SUxU{a}{b}$, which is of rank~$b$. If $b \ge 2$, then the
reducible slice representation is non-polar and polarity minimal by \cite{bergmann}, and the $H$-action
is non-polar and polarity minimal by Lemma~\ref{PolMinHered}. If $b = 1$, then $H_2$ acts transitively
on the orbit~$H_1 \cdot \eK$, but an explicit calculation using Proposition~\ref{PolCrit} shows that the
$H_2$-action on $G / K$ is non-polar. Let $H \subset H_2$ be a proper closed subgroup, then $H$ acts
non-transitively on $H_1 \cdot \eK$ by~\cite{oniscik}. If the $H$-action on $G / K $ is polar, then also
the $H$-action restricted to $H_1 \cdot \eK \cong \C \P^a$ is polar by Lemma~\ref{PolarOnTotGeodSubmf}
and it has a singular orbit $H \cdot p$ by Corollary~\ref{SingularOrbit}. The normal space~$\pp \cap
\m_1$, see Lemma~\ref{SubHermann}, of the $H_1$-action on~$G / K$ contains a submodule which is
equivalent to the slice representation at~$p \in H_1 \cdot \eK$ of the $H$-action on~$H_1 \cdot \eK$
after a $\U(1)$-factor is removed from both representations. Since both modules belong to the polar
slice representation of the~$H$-action on~$G / K$, it follows from~\cite{bergmann} that~$H$ is at most
three-dimensional, a contradiction with Proposition~\ref{CoisoDimBound}.

\par

Now assume $H_2'$ is a non-symmetric irreducible maximal connected subgroup of~$\Sp(a + b)$. It follows
from what we have shown so far that this can only happen if $\rk( H_1 \cdot \eK ) = b = 1$. Assume $H_2'
= \SO(p) \otimes \Sp(q)$, then the $H_2$-action on~$H_1 \cdot \eK$ is non-polar and polarity minimal by
Proposition~\ref{TensGrass}. If $H_2'$ is a simple irreducible maximal connected subgroup of~$\Sp(a +
b)$, then it follows from the results of Section~\ref{SubSimple} that the action of $H_2$ on $G / K$ is
non-polar and polarity minimal, since if the action of $H_2'$ on $\G_k(\H^{a + b})$ for $2 \le k \le a +
b - 2$ is excluded by Proposition~\ref{ClassDimBounds}, then also the action of $H_2$ on $\G_k(\H^{a + b
+ 1})$ is excluded by a dimension count.

\paragraph{\bf BD\,I-I}

Let $H_1 = \SO(a + b)$, $G / K = \SO(a + b + 1) / \SO(a) \times \SO(b + 1)$, $a + b \ge 6$, $a \ge b \ge
1$. Assume $H \subset H_1$ is a closed connected subgroup acting polarly on~$G / K$. Without loss of
generality we may assume that $H \subseteq \SO(a + b)$ acts irreducibly on~$\R^{a + b}$, since otherwise
the $H$-action on $G / K$ is a subaction of a Hermann action of cohomogeneity~$\ge 2$, see
Section~\ref{SubHermannHighRk}. By Lemma~\ref{CommInvol}, the $H_1$-orbit $H_1 \cdot \eK$ is a totally
geodesic submanifold isometric to~$\SO(a + b) / \eS(\O(a) \times \O(b))$ and $H$ acts polarly on~$H_1
\cdot \eK$ (the action may be transitive).
\par
We do not need to consider the case where $\rk(H_1 \cdot \eK) = b = 1$, since polar actions
on~$\G_2(\R^{a+2})$ were completely classified in~\cite{pth2}, hence we may assume $b \ge 2$. Let us
first consider the cases where $H \subset H_1$ is not a symmetric subgroup. Then it follows from what we
have shown so far and \cite{oniscik} that we have one of the following:
\begin{itemize}

    \item $H = \LG_2,\quad a + b = 7,\quad b = 2,\,3$;

    \item $H = \Spin(7),\quad a + b = 8,\quad b = 2,\,3,\,4$;

    \item $H = \Spin(9),\quad a + b = 16,\quad b =2$;

    \item $H = \Sp(n) \cdot \Sp(1),\quad a + b = 4n,\quad b =2$;

    \item $H = \U(n),\quad a + b = 2n$.

\end{itemize}
In the case of the $\LG_2$-actions, an explicit calculation using Proposition~\ref{PolCrit} shows that
the actions are non-polar; subgroups of $\LG_2$, see~(\ref{G2MaxSubGroups}), are either reducible or are
excluded by a dimension count. The actions of~$\Spin(7)$ are orbit equivalent to the $\SO(8)$-action in
case $b = 2,\,3$, in case $b = 4$ the action can be shown to be non-polar by an explicit calculation.
Subgroups of~$\Spin(7)$ are either contained in groups treated below or ruled out by a dimension count.
The actions of $\Spin(9)$ and $\Sp(n) \cdot \Sp(1)$ can be excluded by replacing~$K$ with the conjugate
subgroup~$K' = \SO(3) \times \SO(a)$, the actions on $H_1 \cdot \eK'$ have already been shown to be
non-polar and polarity minimal. Assume now $H = \U(\frac{a + b}{2})$ \label{subBDI-I}. The slice
representation of the $H$-action at $\eK$, as can be seen from Tables~\ref{THermannSliceRep} and
\ref{TCohOneTrHermann}, contains a module equivalent to the representation of~$\U \left( \lfloor
\frac{a}{2} \rfloor \right) \times \U \left( \lfloor \frac{b}{2} \rfloor \right)$ on $\C^{\lfloor
\frac{a}{2} \rfloor} \otimes \C^{\lfloor \frac{b}{2} \rfloor} \oplus \C^{\lfloor \frac{a}{2} \rfloor}$,
which is non-polar~\cite{bergmann} and polarity minimal by Lemma~\ref{PolMinCriteria} since $\lfloor
\frac{a}{2} \rfloor \ge 2$.

\paragraph{\bf C\,I-II}

Let $H$ be closed connected subgroup of $H_1 = \Sp(p){\times}\Sp(1)$ acting polarly on $G / K = \Sp(p+1) /
\U(p+1)$, $p \ge 2$. We first observe that the actions of $\Sp(p)$ and $\Sp(p){\times}\U(1)$ are not orbit
equivalent to the $H_1$-action; since the normal space at a regular orbit is not a Lie triple system,
these actions are non-polar by Proposition~\ref{PolCrit}. Now assume $H \subseteq H_2 = H_2' \times
\Sp(1)$, where $H_2'$ is a maximal connected subgroup of~$\Sp(n)$. We may assume $H_2' \subset \Sp(n)$
acts irreducibly on~$\H^n$, since otherwise $H_2$ is a subgroup of~$\Sp(p+1-k) \times \Sp(k)$, $2 \le k
\le p - 1$, see Section~\ref{SubHermannHighRk}. Consider the action of $H_2 = \U(p) \times \Sp(1)$ on $G
/ K$, then the slice representation of the isotropy subgroup $\U(p) \times \U(1)$ is equivalent to the
isotropy representation of~$\Sp(p) / \U(p)$ plus $\C^p \otimes \C^1$. This representation is
non-polar~\cite{bergmann} and the action of~$H_2$ on~$G / K$ is polarity-minimal by
Lemma~\ref{PolMinHered}~(iii). If $H_2 = \rho(H_2') \times \Sp(1)$, where $\rho \colon H_2' \to \Sp(p)$
is an irreducible representation of the simple compact Lie group $H$, then the $H$-action on $G / K$ is
non-polar by Lemma~\ref{SimpleStruct2}. Tensor product subgroups~$H_2'$ have been excluded in
Proposition~\ref{SubHerTensStruct}.

\paragraph{\bf C\,II-II}

Let $H_1 = \Sp(a + b) \times \Sp(1)$ and let $G / K = \Sp(a + b + 1) / \Sp(a) \times \Sp(b + 1)$, $a \ge
b \ge 1$, $a + b \ge 3$. First observe that the action of $\Sp(a + b)$ on~$G / K$ is orbit equivalent to
the $H_1$-action. Now assume $H \subset H_2 = H_2' \times \Sp(1)$ is a closed connected subgroup acting
polarly on~$G / K$, where $H_2' \subset \Sp(a + b)$ is a maximal connected subgroup. We may assume that
$H_2' \subseteq \Sp(a + b)$ acts irreducibly on~$\H^{a + b}$, since otherwise the $H$-action on $G / K$
is a subaction of a Hermann action of cohomogeneity~$\ge 2$, which were examined in
Section~\ref{SubHermannHighRk}. It follows from Lemma~\ref{CommInvol} that the $H_1$-orbit $H_1 \cdot
\eK$ is a totally geodesic submanifold isometric to~$\Sp(a + b) / \Sp(a) \times \Sp(b)$ on which $H$
acts polarly.
\par
Assume $H_2'$ is an irreducible symmetric subgroup of~$\Sp(a + b)$, hence conjugate to~$\U(a + b)$.
However, the isotropy subgroup of the $H_2$-action at~$\eK$ is $\U(a) \times \U(b) \times \Sp(1)$, its
slice representation contains two equivalent modules, thus it is non-polar and polarity minimal
by~Lemma~\ref{PolMinCriteria}, parts~(ii) and~(iii). The normal space contains a Lie triple system
corresponding to a totally geodesic submanifold isometric to~$\H \P^a$, hence the $H_2$-action is
non-polar and polarity minimal by Lemma~\ref{PolMinHered}~(iii).
\par
Now assume $H_2'$ is a non-symmetric irreducible maximal connected subgroup of~$\Sp(a + b)$ acting
polarly on~$H_1 \cdot \eK = \Sp(a + b) / \Sp(a) \times \Sp(b)$. It follows from what we have shown so
far that this can only happen if $\rk( H_1 \cdot \eK ) = 1$, i.e.\ $b = 1$. Assume $H_2' = \SO(p)
\otimes \Sp(q)$, then the $H_2$-action on~$H_1 \cdot \eK$ is non-polar and polarity minimal by
Proposition~\ref{TensGrass}. If $H_2'$ is a simple irreducible subgroup of~$\Sp(a + b)$, then it follows
from the results of Section~\ref{SubSimple} that the action of $H_2$ on $G / K$ is non-polar and
polarity minimal, since if the action of $H_2'$ on $\G_k(\H^{a + b})$ for $2 \le k \le a + b - 2$ is
excluded by Proposition~\ref{ClassDimBounds}, then so is the action of $H_2$ on $\G_k(\H^{a + b + 1})$.

\paragraph{\bf D\,I-III}

Let $H$ be a closed connected subgroup of $\SO(n-k) \times \SO(k)$, $k = 1,2,3$ acting on $G / K =
\SO(2n) / \U(n)$, $n \ge 3$. We first study actions of closed connected subgroups~$H$ in $H_2 = H_2'
\times \SO(k)$. The cases where $H_2' \subset \SO(2n-k)$ is a simple irreducible or tensor product
subgroup were excluded by Lemma~\ref{SimpleStruct2} and Lemma~\ref{SubHerTensStruct}. Let us consider
the case where $H_2'$ is a symmetric subgroup of~$\SO(2n-k)$. If $k = 2$ and $H_2' = \U(n-1)$, then the
$H$-action has a fixed point. Thus it remains the case where $H_2' = \SO(2n-k-\ell) \times \SO(\ell)$.
We may assume $k + \ell \le 3$ since otherwise the $H$-action on $G / K$ is a subaction of a Hermann
action of cohomogeneity~$\ge 2$, see Section~\ref{SubHermannHighRk}. If $k = \ell = 1$ then we obtain
the cohomogeneity one action of $H = \SO(2n-2)$ on $G / K$, we have already seen that no closed proper
subgroup of this group acts polarly. In the case where $k + \ell = 3$, an explicit calculation using
Proposition~\ref{PolCrit} shows that the actions of $\SO(2n-3)$ and $\SO(2n-3) \times \SO(2)$ on $G / K$
are non-polar; we have already excluded any closed subgroups of these two groups.
\par
In case $n = 3$ there is an additional maximal connected subgroup, i.e.\ $\Delta \SO(3) = \{ (g,g) \mid
g \in \SO(3) \} \subset \SO(3) \times \SO(3)$, but this action has a fixed point.

\paragraph{\bf D\,III-I}

Let $H$ be a closed connected subgroup of $H_1 = \U(n)$ acting polarly on the symmetric space~$G / K =
\SO(2n) / \SO(n-3) \times \SO(3)$, $n \ge 3$. It follows from Theorem~2.1 in~\cite{kollross} that the
conjugacy classes of maximal connected subgroups $H_2$ in~$\U(n)$ are exhausted by $\SU(n)$, and $H_2'
\otimes \U(1)$ where $H_2'$ runs through the maximal connected subgroups of~$\SU(n)$; see
Proposition~\ref{MaxSubgrSUn}. We observe first that the action of $\SU(n)$ is orbit equivalent to the
$\U(n)$-action. Now assume $H \subset H_2' \otimes \U(1)$. We do not need to consider reducible
subgroups~$H_2'$ since they lead to subactions of Hermann actions with cohomogeneity~$\ge 2$, which were
treated in Section~\ref{SubHermannHighRk}. Also, if $H_2' = \SO(n)$, then the same argument as in
Proposition~\ref{TensGrass} shows that the action is non-polar and polarity minimal. If $H_2' =
\Sp(n/2)$, then $H_2$ is contained in $\Sp(n/2) \cdot \Sp(1)$ and the action of any closed subgroup $H
\subseteq H_2$ was shown not to be polar in Proposition~\ref{TensGrass}. Assume now $H_2'$ is a tensor
product subgroup $\SU(p) \otimes \SU(n/p)$, then the action of~$H$ is non-polar according to
Proposition~\ref{SubHerTensGrass}. Finally, let $H_2'$ be given by an irreducible representation $\rho
\colon H_2' \to \SU(n)$ where $H_2'$ is a simple compact Lie group; these actions were excluded in
Section~\ref{SubSimple}.

\paragraph{\bf E\,II-IV}

Let $H \subseteq H_1 = \SU(6) \cdot \Sp(1)$ be a closed connected subgroup acting on $G / K = \LE_6 /
\LF_4$. Then $H$ is contained in one of the maximal connected subgroups of $\SU(6) \cdot \Sp(1)$. By
\cite{conlon}, we may assume the involutions of $\LE_6$ corresponding to $H_1$ and $K$ commute such that
the totally geodesic $H_1$-orbit through~$\e K$ is isometric to $\SU(6) \cdot \Sp(1) / \Sp(3) \cdot
\Sp(1) \cong \SU(6) / \Sp(3)$, see Lemma~\ref{CommInvol}; the slice representation is equivalent to the
isotropy representation of $\Sp(4) / \Sp(3) \cdot \Sp(1)$, this shows that the $\Sp(1)$-factor is
inessential for the $H_1$-action. Now assume $H_2' \subseteq \SU(6)$ is a symmetric subgroup and $H
\subseteq H_2 = H_2' \cdot \Sp(1)$ is a closed connected subgroup. If $H_2' = \Sp(3)$ then the
$H$-action on $G / K$ has a fixed point. If $H_2' = \SO(6)$, then the connected component of an isotropy
subgroup is $\U(3) \times \Sp(1)$ and its slice representation is equivalent to the adjoint
representation of $\SU(3)$ plus $\chi( \Sp(4),\, \Sp(3) \times \Sp(1) )|_{\U(3) \times \Sp(1)}$, see
Tables~\ref{THermannSliceRep} and \ref{TCohOneTrHermann}, it is non-polar~\cite{bergmann} and polarity
minimal by Proposition~\ref{PolMinCriteria}. Hence the $H_2$-action on~$G / K$ is non-polar and polarity
minimal by Lemma~\ref{PolMinHered}~(iii).

Now assume $H_2'$ is one of the groups $\SUxU{k}{6-k}$, $k = 1,2,3$. If $k = 3$, then $(H_2 \cap K)_0
\cong \Sp(1) \cdot \U(1) \cdot \Sp(1) \cdot \Sp(1)$ and the slice representation is equivalent to
$$\chi( \Sp(2),\, \Sp(1) \times \Sp(1) ) \oplus \chi( \Sp(4),\,
\Sp(3) \times \Sp(1) )|_{\Sp(1) \cdot \U(1) \cdot \Sp(1) \cdot \Sp(1)},$$ where both modules are polar,
but their sum is non-polar, see~\cite{bergmann}; using the results of ~\cite{bergmann}, it can be
directly verified that this representation is polarity minimal and it follows that the $H_2$-action
on~$G / K$ is non-polar and polarity minimal by Lemma~\ref{PolMinHered}. If $k = 2$, then $(H_2 \cap
K)_0 \cong \Sp(2) \cdot \Sp(1) \cdot \Sp(1)$ and the corresponding slice representation is equivalent to
$$\chi( \Sp(3),\, \Sp(2) \times \Sp(1) ) \oplus \chi( \Sp(4),\, \Sp(3) \times \Sp(1) )|_{\Sp(2) \cdot
\Sp(1) \cdot \Sp(1)}$$ hence non-polar by~\cite{bergmann} and polarity minimal by
Proposition~\ref{PolMinCriteria}~(ii) and the $H_2$-action on~$G / K$ is non-polar and polarity minimal
by Lemma~\ref{PolMinHered}. Finally, if $k = 1$ then $H_2$ acts transitively on the $H_1$-orbit through
$\eK$, and the group $(H_2 \cap K)_0 \cong \U(1) \times \Sp(2) \times \Sp(1)$ acts polarly on
$\Exp_{\eK} (\N_{\eK} (H_2 \cdot \eK) \cong \H \P^3$ by \cite{pth1}, cf.\ Lemma~\ref{SubHermann};
however, an explicit calculation shows that $[\nu, \nu] \not \perp \h_2$, where $\nu \subset \m_1 \cap
\pp$ is the tangent space to a section of the $H_2 \cap K$-action and hence the $H_2$-action and the
orbit equivalent action of $\SU(5) \cdot \Sp(1)$ are non-polar by Lemma~\ref{SliceActionTrans}. The
actions of $\SUxU{1}{5}$ and $\SUxU{1}{5} \cdot \U(1)$ are non-polar since the slice representations at
$\eK$ are non-polar. Any other closed subgroups of $H_2 = \SUxU{1}{5} \cdot \Sp(1)$ are contained in the
groups treated above or excluded by Proposition~\ref{ClassDimBounds}.

\paragraph{\bf E\,IV-II}

Consider now closed connected subgroups $H$ of~$H_1 = \LF_4$ acting polarly on $\LE_6 / \SU(6) \cdot
\Sp(1)$. It follows from Lemma~\ref{ProdSphDim} that $\dim(H) \ge 28$. By \cite{dynkin1}, the only
closed connected subgroup of $\LF_4$ of sufficient dimension is $\Spin(9)$. By conjugation, the subgroup
$H = \Spin(9) \subset \LF_4$ can be chosen such that the connected component of the isotropy group $(H
\cap K)_0$ is $\Sp(2) \cdot \Sp(1)^2 \cong \Spin(5) \cdot \Spin(4)$, see Table~\ref{TCohOneTrHermann}.
From Table~\ref{TCohOneTrHermann} one sees further that the slice representation restricted to $(H \cap
K)_0$ is equivalent to the isotropy representation of $\Sp(3) / \Sp(2) \cdot \Sp(1)$ plus the isotropy
representation of $\Sp(4) / \Sp(3) \cdot \Sp(1)$ restricted to $\Sp(2) \cdot \Sp(1)$, hence non-polar
by~\cite{bergmann}.

\paragraph{\bf E\,III-IV}

Let $H$ be a closed connected subgroup of~$H_1 = \Spin(10) \cdot \U(1)$ acting polarly on $G / K = \LE_6
/ \LF_4$. Let us first show that the action of~$\Spin(10)$ on~$G / K$ is non-polar. Assume the converse,
i.e.\ the action of~$H = \Spin(10)$ on $G / K = \LE_6 / \LF_4$ is polar. Since this action is of
cohomogeneity~two, it follows from Proposition~\ref{HyperPMax} that the sections are locally isometric
to a two-sphere. The $H_1$-orbit through~$\eK$ is totally geodesic by Lemma~\ref{CommInvol} and locally
isometric to~$\eS^9 \times \eS^1$, where the~$\eS^1$-factor is the orbit of the $\U(1)$-factor in $H_1 =
\Spin(10) \cdot \U(1)$, hence totally geodesic in~$G / K$. It follows from the
decomposition~(\ref{CommInvDecomp}) that $\T_{\eK} \eS^1 \perp \T_{\eK} (H \cdot \eK)$, hence $\T_{\eK}
\eS^1 \subset \T_{\eK} \Sigma$, where $\Sigma$ is a section of the $H$-action on~$G / K$
containing~$\eK$. Since the Lie algebra of the $\U(1)$-factor in $H_1 = \Spin(10) \cdot \U(1)$ is
contained in~$\T_{\eK} \Sigma \subset \pp$, it follows by Proposition~\ref{LieTripleSystems} that this
$\U(1)$-factor acts on~$\Sigma$ as a group of transvections. Now, since this $\U(1)$-action commutes
with the $H$-action, it follows that any two points of~$\Sigma$ which lie in the same $\U(1)$-orbit are
of the same orbit type with respect to the $H$-action on~$G / K$. In particular, all singular orbits of
the $H$-action on~$G / K$ intersect $\Sigma$ in the $\U(1)$-orbit which is covered by a great circle
of~$\tilde \Sigma \cong \eS^2$, since all reflection hypersurfaces~$\{P_j\}_{j \in J}$ have to be
invariant under the $\U(1)$-action induced on~$\tilde \Sigma$, see Lemma~\ref{TotGeodHyperWeyl}.
However, the $\U(1)$-action on $\Sigma$, which is isometric to~$\R \P^2$ or $\eS^2$, has at least one
fixed point~$p \in \Sigma$. It follows that this point~$p \in G / K$ lies in a regular orbit of the
$H$-action on~$G / K$, but is left fixed by the $\U(1)$-factor of $H_1 = \Spin(10) \cdot \U(1)$. Hence
the connected component of the isotropy subgroup $(H_1)_p$ is a subgroup~$L \cdot \U(1)$, where $L
\subset \Spin(10)$ is $20$-dimensional, a contradiction, since the principal isotropy subgroup
is~$\Spin(7)$.
\par
Now we may assume that the group $H$ acting polarly on $G / K$ is contained in $H_2' \cdot \U(1)$, where
$H_2'$ is a maximal connected subgroup of~$\Spin(10)$. It follows from Lemma~\ref{ProdSphDim} that
$\dim(H) \ge 16$. This implies that $H_2'$ is one of the following, see Proposition~\ref{MaxSubgrSOn}:
\begin{equation}\label{Spin10U1MaxSub}
  \Spin(9),\,\, \Spin(8) \cdot \SO(2),\,\,
  \Spin(7) \cdot \Spin(3),\,\, \Spin(6) \cdot \Spin(4),\,\, \U(5).
\end{equation}
The actions of these groups are non-polar and polarity minimal by Lemma~\ref{PolMinHered}~(iv). If $H_2'
= \Spin(9)$, then the $H$-action can also be shown to be non-polar by Corollary~\ref{PolFixOnHomSp}.

\paragraph{\bf E\,IV-III}

Assume $H$ is a closed connected subgroup of~$\LF_4$ acting polarly on $G / K = \LE_6 / \Spin(10) \cdot
\U(1)$. Proposition~\ref{CoisoDimBound} implies $\dim(H) \ge 28$. By \cite{dynkin1}, the only maximal
connected subgroups $H \subset \LF_4$ of dimension~$\ge 28$ is $\Spin(9)$. It follows from
Table~\ref{TCohOneTrHermann} that the action of~$\Spin(9)$ leaves a point fixed.

\paragraph{\bf F\,II-I}

Let $H$ be a closed connected subgroup of~$H_1 = \Spin(9)$ acting polarly on $G / K = \LF_4 / \Sp(3)
\cdot \Sp(1)$. Lemma~\ref{ProdSphDim} implies $\dim(H) \ge 20$. The only closed connected subgroups
of~$\Spin(9)$ of dimension~$\ge 20$ are $\Spin(8)$, $\Spin(7) \cdot \SO(2)$ and $\Spin(7)$.

The subgroup $H = \Spin(8) \subset \Spin(9)$ may be chosen such that $(H_1 \cap K)_0 = \Sp(2) \cdot
\Sp(1) \cong \Spin(5) \cdot \Spin(3)$. The group $H$ acts with cohomogeneity one on the orbit $H_1 \cdot
\eK$, which is covered by $\Spin(9) / \Spin(5) \cdot \Spin(4)$. With our choice of the subgroup
$\Spin(8) \subset \Spin(9)$, the slice representation of the $H$-action on $G / K$ at~$\eK$ is
equivalent (on the Lie algebra level) to the representation of $\Sp(2) \cdot \Sp(1) \cong \Spin(5) \cdot
\Spin(3)$ on $\H^2 \otimes_{\H} \H^1 \oplus \R^5$, hence non-polar by~\cite{bergmann}.

The action of $\Spin(7) \cdot \SO(2)$ has an isotropy subgroup whose connected component is isomorphic
to~$\Spin(5) \cdot \SO(2) \cdot \SO(2) \cong \Sp(2) \cdot \U(1) \cdot \U(1)$ and whose slice
representation is $\R^5 \otimes \R^2 \oplus \chi( \Sp(3),\, \Sp(2) \cdot \Sp(1))|_{\Sp(2) \cdot \U(1)}$,
hence non-polar~\cite{bergmann}. This also shows that the $\Spin(7)$-action is non-polar.
\label{MainProofEnd}
\end{proof}


\bibliographystyle{amsplain}

\end{document}